\documentclass{amsart}
\usepackage{amssymb,amscd}
\usepackage[all]{xy}
\usepackage{enumerate}

\numberwithin{equation}{subsubsection}

\swapnumbers
\theoremstyle{plain}
\newtheorem{thm}[subsection]{Theorem}
\newtheorem{prop}[subsubsection]{Proposition}
\newtheorem{lemma}[subsubsection]{Lemma}
\newtheorem{cor}[subsection]{Corollary}

\theoremstyle{definition}
\newtheorem{defn}[subsubsection]{Definition}
\newtheorem{defns}[subsubsection]{Definitions}

\theoremstyle{remark}
\newtheorem{rem}[subsubsection]{Remark}

\newcommand{\Curve}{\mathcal{C}}
\renewcommand{\O}{\mathcal{O}}
\newcommand{\FF}{\mathcal{F}}

\newcommand{\LL}{\mathcal{L}}

\newcommand{\F}{\mathbb{F}}

\newcommand{\Fq}{{\mathbb{F}_q}}
\newcommand{\Z}{\mathbb{Z}}
\newcommand{\Q}{\mathbb{Q}}
\newcommand{\R}{\mathbb{R}}
\newcommand{\C}{\mathbb{C}}
\newcommand{\A}{\mathbb{A}}
\renewcommand{\P}{\mathbb{P}}

\newcommand{\m}{\mathfrak{m}}
\newcommand{\n}{\mathfrak{n}}

\newcommand{\G}{\backslash G}

\newcommand{\sha}{{\hbox to 10pt{\rlap{\hskip2.8pt\vrule
height6pt\hskip1.6pt\vrule height6pt\hskip1.6pt
\vrule height6pt}\hskip1pt\vrule height0.8pt width 8pt\hskip1pt}}}

\newcommand{\into}{\hookrightarrow}
\newcommand{\onto}{\twoheadrightarrow}
\newcommand{\isoto}{\tilde{\to}}
\newcommand{\tensor}{\otimes}
\newcommand{\compose}{\circ}
\def\nodiv{\mathrel{\mathchoice{\not|}{\not|}{\kern-.2em\not\kern.2em|}{\kern-.2em\not\kern.2em|}}}

\newcommand{\SL}{\mathrm{SL}}
\newcommand{\GL}{\mathrm{GL}}
\newcommand{\Sp}{\mathrm{Sp}}
\newcommand{\SO}{\mathrm{SO}}
\newcommand{\Or}{\mathrm{O}}

\DeclareMathOperator{\im}{Im}

\DeclareMathOperator{\res}{Res}
\DeclareMathOperator{\ind}{Ind}

\DeclareMathOperator{\cond}{Cond}
\DeclareMathOperator{\ord}{ord}
\DeclareMathOperator{\rk}{Rank}

\DeclareMathOperator{\Hom}{Hom}
\DeclareMathOperator{\aut}{Aut}

\DeclareMathOperator{\gal}{Gal}
\DeclareMathOperator{\spec}{Spec}

\DeclareMathOperator{\sgn}{sgn}

\renewcommand{\vector}[1]{\mathbf{#1}}
\renewcommand{\v}{\vector{v}}

\newcommand{\<}{\langle}
\renewcommand{\>}{\rangle}

\newcommand{\pmat}[1]{\begin{pmatrix}#1\end{pmatrix}}

\newcommand{\Fqbar}{{\overline{\mathbb{F}}_q}}

\newcommand{\Fqn}{{\mathbb{F}}_{q^n}}

\newcommand{\Qbar}{{\overline{\mathbb{Q}}}}
\newcommand{\Ql}{{\mathbb{Q}_\ell}}
\newcommand{\Zl}{{\mathbb{Z}_\ell}}
\newcommand{\Qlbar}{{\overline{\mathbb{Q}}_\ell}}

\newtheorem*{thm*}{Theorem}
\newcommand{\GG}{{\mathcal{G}}}
\newcommand{\HH}{{\mathcal{H}}}
\newcommand{\XX}{{\mathcal{X}}}
\newcommand{\DD}{{\mathcal{D}}}
\newcommand{\uE}{\underline{E}}
\newcommand{\Xbar}{\overline{X}}
\newcommand{\Cbar}{\overline{\Curve}}
\renewcommand{\ggg}{\mathfrak{g}}
\newcommand{\hh}{\mathfrak{h}}
\renewcommand{\G}{\mathbb{G}}
\newcommand{\Ubar}{\overline{U}}

\newcommand{\sdp}{{\rtimes}}

\def\hyp#1{{\advance\hsize by -1in\parindent=0pt
    \vtop{#1}}}

\begin{document}
\title[Geometric non-vanishing]{Geometric non-vanishing}
\author{Douglas Ulmer}
\address{Department of Mathematics \\
  University of Arizona \\ Tucson, AZ 85721}
\email{ulmer@math.arizona.edu} 
\thanks{This paper is based upon work
  supported by the National Science Foundation under Grant No.
  DMS 0070839} 
\date{May 13, 2004}
\begin{abstract}
  We consider $L$-functions attached to representations of the Galois
  group of the function field of a curve over a finite field.  Under
  mild tameness hypotheses, we prove non-vanishing results for twists
  of these $L$-functions by characters of order prime to the
  characteristic of the ground field and more generally by certain
  representations with solvable image.  We also allow local
  restrictions on the twisting representation at finitely many places.
  Our methods are geometric, and include the Riemann-Roch theorem, the
  cohomological interpretation of $L$-functions, and monodromy
  calculations of Katz.  As an application, we prove a result which
  allows one to deduce the conjecture of Birch and Swinnerton-Dyer for
  non-isotrivial elliptic curves over function fields whose
  $L$-function vanishes to order at most 1 from a suitable
  Gross-Zagier formula.
\end{abstract}
\maketitle

\section{Introduction}\label{s:intro}

Non-vanishing results have long played an important role in the
application of $L$-functions to arithmetic, beginning with Dirichlet's
1837 proof of the infinitude of primes in an arithmetic progression.
The area remains active and there is a vast literature.  We refer to
\cite{BFH}, \cite{Murty}, \cite{Goldfeld}, and their bibliographies
for an overview of some recent work in the area.  

Over number fields, one typically considers {\it automorphic\/}
$L$-functions, since only these are known to have good analytic
properties.  Here, proofs of non-vanishing results necessarily use
automorphic methods such as modular symbols, Fourier coefficients of
half-integral weight forms, metaplectic Eisenstein series, or average
value computations based on character sum estimates.  Over function
fields, similar automorphic ideas can be applied (see for instance
\cite{HoffRosen} and \cite{DuttaGupta}), but the theory is much less
developed.

On the other hand, in the function field case, one has a much better
understanding of {\it motivic\/} $L$-functions, i.e., those attached
to Galois representations, because of Grothendieck's analysis of
$L$-functions.  This powerful cohomological interpretation allows one
to apply geometric methods to the study of these $L$-functions.

The goal of this paper is to use geometric methods to prove a very
general non-vanishing result for twists of motivic $L$-functions over
a function field.  Because Lafforgue has proven the Langlands
correspondence for $\GL_n$ over function fields \cite{Lafforgue}, our
results apply to many automorphic $L$-functions as well.

To state the result more precisely, let $\Curve$ be a smooth, proper,
geometrically irreducible curve over a finite field $\Fq$ of
characteristic $p$, $F=\Fq(\Curve)$, and $\overline{F}$ a separable
closure of $F$.  Let $\Fqn\subset\overline F$ be the subfield of $q^n$
elements, $\Fqbar=\cup_{n\ge1}\Fqn$, and set $F_n=\Fqn(\Curve)$
($n\ge1$) and $F_{\infty}=\Fqbar(\Curve)$.  Let $\rho$ be a continuous,
absolutely irreducible $\ell$-adic representation of
$\gal(\overline{F}/F)$ for some $\ell\neq p$.  We assume that $\rho$
is unramified outside a finite set of places of $F$ and that it is
geometrically absolutely irreducible, i.e., that it is absolutely
irreducible when restricted to $\gal(\overline{F}/F_{\infty})$.  We write
$L(\rho,F,s)$ for the $L$-function attached to $\rho$ (see
\ref{sss:L-funs} for the definition) and $L(\rho,K,s)$ for the
$L$-function of $\rho|_{\gal(\overline{F}/K)}$ for any finite
extension $K$ of $F$ contained in $\overline{F}$.

Fix a positive integer $d$ not divisible by $p$ and a complex number
$s_0$.  We seek elements $f\in F^\times$ such that $F(f^{1/d})$ has
degree $d$ over $F$ and the ratio 
\begin{equation*}
\frac{L(\rho,F(f^{1/d}),s)}{L(\rho,F,s)}
\end{equation*}
is non-vanishing at $s=s_0$.  We can find such $f$ if we first
replace $F$ with $F_n$ for sufficiently large $n$.  More
precisely, here is the statement of a very weak version of our main
result:

\begin{thm}\label{thm:intro-main}
  Assume that $d|q-1$ and that $\rho$ is everywhere at worst tamely
  ramified or that $p>\deg\rho+2$.  Then for infinitely many
  integers $n$, there exists an element $f\in F_n^\times$ such that
  the extension $F_n(f^{1/d})$ of $F_n$ has degree $d$ and
\begin{equation*}
\frac{L(\rho,F_n(f^{1/d}),s)}{L(\rho,F_n,s)}
\text{ does not vanish at $s=s_0$.}
\end{equation*}
\end{thm}

Before discussing the strengthenings of this result which are our
goal, let us remark on the difference between it and what one might
expect from analogy with the classical case.  Fix $F$ as above and
consider extensions of the form $K=F(f^{1/d})$ partially ordered by
the degree of their conductors.  Then one might expect that for
sufficiently large conductor, there exists an extension $K$ of this
type such that the non-vanishing conclusion of the theorem holds.
More optimistically, one might hope that as the degree of the
conductor goes to infinity, the proportion of the extensions $K$ that
satisfy the non-vanishing conclusion is positive and bounded away from
0.  This may well be true, but the methods of this paper lead to
a slightly different point of view (for reasons explained in
Section~\ref{s:easy-case}). Namely, we consider extensions
$K=F_n(f^{1/d})$ of {\it bounded conductor\/} for varying $n$.  We
show that for large $n$ there exist extensions for which the
non-vanishing conclusion holds.  Our methods also show that the
density of extensions for which we have non-vanishing is positive and
bounded away from 0 as $n\to\infty$.  (We do not, however, state
explicitly the densities.  If needed, they may easily be extracted
from the proofs of Proposition~\ref{prop:good-pts-density} and
Corollary~\ref{cor:non-vanishing-density}.)

The first strengthening of Theorem~\ref{thm:intro-main} concerns the
hypothesis $d|q-1$.  Because of it, the ratio in the conclusion of the
theorem is a product of twists $L(\rho\tensor\chi,F_n,s)$ where $\chi$
runs through the non-trivial characters of
$\gal(F_n(f^{1/d})/F_n)\cong\Z/d\Z$.  Thus Theorem~\ref{thm:intro-main}
is about the non-vanishing of abelian twists of $L(\rho,F_n,s)$.  In
our main theorem, we drop the condition that $d|q-1$ and so the
extension $F_n(f^{1/d})/F_n$ may not be Galois.  This means
that we have to consider twists of $\rho$ by certain non-abelian
representations of $\gal(\overline{F}/F_n)$.

The second strengthening is that we are able to impose local
conditions (splitting, inertness, ramification) on the extension
$F_n(f^{1/d})/F_n$ at finitely many places.  

The third strengthening is that we make a statement for all
sufficiently large $n$.  It turns out that for 
certain data ($\rho$, $d$, local conditions, and points $s_0$), the ratio
$L(\rho,F_n(f^{1/d}),s)/L(\rho,F_n,s)$ vanishes at $s=s_0$ for
arbitrarily large $n$ and all $f\in F_n^\times$ satisfying the local
conditions.  (Think for example of a situation where the local
conditions force the sign in a functional equation to be $-1$.)  In
these ``exceptional situations'' our result will assert simple
vanishing, rather than non-vanishing.  The analysis of the exceptional
situations is somewhat intricate.  From a monodromy point of view, their
cause is clear enough (it is related to the fact that every odd
dimensional orthogonal matrix has 1 or $-1$ as an eigenvalue), but we
have gone to some pains to describe the exceptional situations in
terms of {\it easily computable\/} (essentially local) data, like
local root numbers and conductors.  This yields criteria which are
well-suited to applications.  The precise result is stated as
Theorem~\ref{thm:main}.

Another strengthening
is that we allow the point $s_0$ to vary with $n$.  I do not know of
any application of this generalization, but it is natural from a certain 
point of view and it does not make the proof any harder.

The case of the theorem where $d|q-1$ and we do not impose local
conditions follows fairly easily from the monodromy calculations
\cite{Katz} of Katz.  The motivation for considering degrees $d$ that
do not divide $q-1$ and local conditions comes from an application to
elliptic curves which was the genesis of this project.  The result
says roughly that any non-isotrivial elliptic curve over $F$ whose
$L$-function vanishes to order $\le1$ can be put into position to
apply a Gross-Zagier formula.  More precisely:

\begin{thm}\label{thm:intro-elliptic}  
  Assume that $F=\Fq(\Curve)$ has characteristic $p>3$ and let $E$ be
  an elliptic curve over $F$ with $j(E)\not\in\Fq$.  Then there exists
  a finite separable extension $F'$ of $F$ and a quadratic extension
  $K$ of $F'$ such that the following conditions hold:
\begin{enumerate}
\item[(a)] $E$ is semi-stable over $F'$.
\item[(b)] There is a place of $F'$, call it $\infty$, where $E$ has split
      multiplicative reduction.
\item[(c)] The place $\infty$ of $F'$ is not split in $K$.
\item[(d)] Every other place of $F'$ where $E$ has bad reduction is split
      in $K$.
\item[(e)] $\ord_{s=1}L(E/F',s)=\ord_{s=1}L(E/F,s)$ and
      $\ord_{s=1}L(E/K,s)$ is odd and $\le\ord_{s=1}L(E/F',s)+1$.  In
      particular, if $\ord_{s=1}L(E/F,s)\le1$, then\hfil\break
      $\ord_{s=1}L(E/K,s)=1$.
\end{enumerate}
\end{thm}

As we have explained elsewhere \cite[3.8]{UlmerAnalogies}, this result
together with a suitably general Gross-Zagier formula implies that the
conjecture of Birch and Swinnerton-Dyer holds for elliptic curves $E$ over
function fields $F$ of characteristic $p>3$ with
$\ord_{s=1}L(E/F,s)\le1$. 

The plan of the paper is as follows.  In the next section we consider the
simplest case of Theorem~\ref{thm:intro-main}, in which we take
$\Curve=\P^1$, $\rho$ the trivial representation, $d=2$, and
$s_0=1/2$.  The result in this case can easily be proven by elementary
methods, but we give a proof which already contains the main ideas of
the general case.  This section is meant for motivation and none of
the rest of the paper relies on it.  In Sections~\ref{s:L-funs} and
\ref{s:forced-zeroes} we discuss some preliminaries on the
factorization of the ratio $L(\rho,F_n(f^{1/d}),s)/L(\rho,F_n,s)$ into
twists of $L(\rho,F_n,s)$ and on local root numbers and conductors and
then use them to analyze the exceptional situations mentioned above.
Then we are ready to state the main theorem in
Section~\ref{s:statement}.  The main body of the proof begins in
Section~\ref{s:local-conds} where we define a variety $X$
parameterizing extensions $F_n(f^{1/d})/F_n$ and study the set of
points of $X$ satisfying local conditions of splitting, inertness, and
ramification.  In Section~\ref{s:cohomology} we review the
cohomological interpretation of $L$-functions and construct a sheaf
$\GG$ on $X$ whose stalks give the twisted $L$-functions we are
studying.  In Section~\ref{s:monodromy} we calculate the monodromy
groups of $\GG$, using crucially the results of \cite{Katz}.  In
Sections~\ref{s:equidistribution}-\ref{s:end-of-proof} we apply a
variant of Deligne's equidistribution theorem and the monodromy
calculations to prove our non-vanishing results.  The application to
elliptic curves is given in Section~\ref{s:elliptic-app}.

This paper relies heavily on the difficult work of Katz \cite{Katz}.
Fortunately, we are able to treat his results as a ``black box'' for
most of the argument (one important exception being the proof of
Proposition~\ref{prop:GGis-distinct}.)  We hope that this paper may
serve as an introduction to some of the powerful ideas in \cite{Katz}.

{\it Acknowledgements:\/} It is a pleasure to thank Nick Katz for
making a preliminary version of \cite{Katz} available to me and for
some helpful remarks at an early stage of the project.  I also thank
Minhyong Kim for encouraging me to think about the problem in its
natural generality and the referee for making several comments and
corrections.

\section{The simplest case}\label{s:easy-case}
In this section we consider the simplest case of
Theorem~\ref{thm:intro-main}, namely that where $\Curve=\P^1$, $\rho$
is the trivial representation, $d=2$, and $s_0=1/2$.  (For brevity, we
use certain notational conventions which are not spelled out until
later, but which are standard and should be clear.)
Since we assume as always that $p\nodiv d$, we have $p>2$.  If $f\in
F_n^\times$ is not a square, then on one hand,
$L(\rho,F_n(\sqrt{f}),s)$ is the zeta function of the hyperelliptic
curve $\Curve_f$ over $\Fqn$ with function field $F_n(\sqrt{f})$, and
on the other hand,
\begin{equation*}
L(\rho,F_n(\sqrt{f}),s)=L(\rho,F_n,s)L(\rho\tensor\chi_f,F_n,s)=
\frac{L(\rho\tensor\chi_f,F_n,s)}{(1-q^{-ns})(1-q^{-n(1-s)})}
\end{equation*}
where $\chi_f$ is the quadratic character of $\gal(\overline{F}/F_n)$
associated to the extension $F_n(\sqrt{f})/F_n$.
This means that 
\begin{equation*}
\frac{L(\rho,F_n(\sqrt{f}),s)}{L(\rho,F_n,s)}=
L(\rho\tensor\chi_f,F_n,s)
\end{equation*}
 is the numerator of the zeta
function of $\Curve_f$.  

Thus Theorem~\ref{thm:intro-main} asserts that for infinitely many
$n$, there exists a hyperelliptic curve over $\Fqn$ whose zeta
function does not vanish at the center point of the functional
equation, namely at $s_0=1/2$.  This in fact holds for all
sufficiently large $n$ and it is possible to give elementary proofs of
this fact, but we need a proof that will work in a much more general
situation.  In the rest of this section we give such a proof in order to
illustrate the main ideas of the proof of
Theorem~\ref{thm:main}.

The first point is to note that 
\begin{equation*}
L(\rho\tensor\chi_f,F_n,s)
=\det\left(1-Fr^n\,q^{-ns}
   \left|H^1(\Curve_f\times\spec\Fqbar,\Ql\right.\right)
\end{equation*}
by Grothendieck's analysis of $L$-functions.  Here $Fr$ is the
endomorphism of $H^1(\Curve_f\times\spec\Fqbar,\Ql)$ induced by the
identity on $\Curve_f$ and the geometric ($q^{-1}$-power) Frobenius on
$\Fqbar$.  Thus we need to study the distribution of eigenvalues of
Frobenius on $H^1$ of hyperelliptic curves and in particular to find
an $f$ such that $q^{n/2}$ is not an eigenvalue of $Fr^n$ on
$H^1(\Curve_f\times\spec\Fqbar,\Ql)$.

To that end, we construct a large family of hyperelliptic curves.
More precisely, fix an odd integer $D\ge3$.  Let $X$ be the variety
over $\Fq$ whose $\Fqn$ points are the monic polynomials of degree $D$
over $\Fqn$ with distinct roots.  I.e., $X$ is obtained from affine space
$\A^D$ by removing a discriminant hypersurface.  Over $X$ we construct
a family $\pi:Y\to X$ of hyperelliptic curves of genus $g=(D-1)/2$ in
such a way that the fiber over $f\in X(\Fqn)$ is the curve
$\Curve_f$.  Explicitly, we view polynomials as rational functions on
$\P^1$.  We have a rational function $f_{univ}$ on $\P^1\times X$,
namely $f_{univ}=x^D+a_1x^{D-1}+\cdots+a_D$ where $x$ is the standard
coordinate on $\P^1$ and $a_1,\dots,a_D$ are the natural coordinates
on $X$. Taking the square root of $f_{univ}$ gives a surface $Y$ with
a map $Y\to\P^1\times X\to X$ with the desired property.

Next we consider the sheaf $\GG=R^1\pi_*\Ql$ on $X$ which is lisse
because $\pi$ is smooth and proper.  The stalk of $\GG$ at a geometric
point over $f\in X(\Fqn)$ is canonically isomorphic to
$H^1(\Curve_f\times\Fqbar,\Ql)$ and so we have united the cohomology
groups we wish to study in one object.  Let $\overline{\eta}$ be a
geometric generic point of $X$ and consider the natural monodromy
representation of $\pi_1(X,\overline{\eta})$ on the stalk
$\GG_{\overline{\eta}}$, which is a $2g$-dimensional $\Ql$ vector
space.  (See Section~\ref{s:cohomology} for more on lisse sheaves and
monodromy representations.)  Let us assume for convenience that $q$ is
a square in $\Ql$ and fix a square root.  Then we twist the
representation of $\pi_1(X,\overline{\eta})$ by the unique character
which sends a geometric Frobenius element at a place $v$ of $X$ to
$q^{-\deg(v)/2}$.  Call the resulting representation $\tau$.

The key input is a calculation of the monodromy group of $\tau$.  More
precisely, write $\pi_1^\text{arith}$ for $\pi_1(X,\overline{\eta})$
and $\pi_1^\text{geom}$ for $\pi_1(X\times\Fqbar,\overline{\eta})$.
Then we define $G^\text{arith}$, the arithmetic monodromy group of
$\tau$, as the Zariski closure of $\tau(\pi_1^\text{arith})$ in
$\GL(\GG_{\overline{\eta}})$.  Similarly, $G^\text{geom}$ is the
Zariski closure of $\tau(\pi_1^\text{geom})$.  By
\cite[1.3.9]{WeilII}, $G^\text{geom}$ is a (not necessarily connected)
semisimple algebraic group over $\Ql$.  Note that $\GG$ carries a
natural alternating form (the cup product on cohomology) with values
in $\Ql(-1)$.  This form is respected by the action of
$\pi_1^\text{arith}$ and so the arithmetic monodromy group lies {\it a
  priori\/} in a symplectic group.  (This is why we introduced the
twist by $Fr_v\mapsto q^{-\deg v/2}$; otherwise, $\pi_1^\text{arith}$
would act by symplectic similitudes.)  Theorem~10.1.18.3 of
\cite{KatzSarnak} is a calculation of this monodromy group.  Namely,
Katz and Sarnak show that $G^\text{geom}$ is the full symplectic group
$\Sp_{2g}$, and therefore so is the {\it a priori\/} larger group
$G^\text{arith}$.

At this point we could apply Deligne's equidistribution result, which
says roughly that Frobenius elements are equidistributed in the
monodromy group.  (This is what we will do in the general case.)  But
in the current simple context, it is more efficient to proceed as
follows.  Let $E_1\subset\Sp_{2g}(\Ql)\subset\GL(\GG_{\overline\eta})$
be the subset of matrices which have $1$ as an eigenvalue.  This is a
proper Zariski closed subset and so there exists an element
$c\in\pi_1^\text{arith}$ such that $\tau(c)\not\in E_1$.

Since $\pi_1^\text{arith}$ is compact, choosing a suitable basis, we may
assume that the image of $\tau$ lies in $\Sp_{2g}(\Zl)$ and then we
may form the reduced representations
$\tau_m:\pi_1^\text{arith}\to\Sp_{2g}(\Z/\ell^m\Z)$.  For large enough $m$
we have that $\det(1-\tau_m(c))\neq0$.  If $f\in X(\Fqn)$ we write
$Fr_{n,f}\in\pi_1^\text{arith}$ for the corresponding geometric Frobenius
element (induced by the map
$\spec\Fqn\to X$); it is well-defined up to conjugacy.  By the
Cebotarev density theorem, for all sufficiently large $n$ there exist
elements $f\in X(\Fqn)$ such that $\tau_m(Fr_{n,f})$ and $\tau_m(c)$
are in the same conjugacy class.  This implies that 1 is not an
eigenvalue of $Fr_{n,f}$ on $\GG_{\overline{\eta}}$ and so $q^{n/2}$
is not an eigenvalue of $Fr^n$ on
$H^1(\Curve_f\times\Fqbar,\Ql)$.  Therefore $s=1/2$ is not a
zero of $L(\rho\tensor\chi_f,F_n,s)$ which is the desired result.

It is clear from this argument why we need to use the extensions $F_n$
in the main theorem.  Indeed, if we consider extensions of $F_n$ (for
varying $n$) of the form $F_n(f^{1/d})$ and of bounded conductor, then
there is a scheme $X$ of finite type whose $\Fqn$ points parameterize
the extensions under consideration and there is a lisse sheaf $\GG$ on
$X$ whose stalks are the cohomology groups related to twisted
$L$-functions.  On the other hand, if we were to consider only extensions
of $F$ of the form $F(f^{1/d})$ then the set of extensions under
consideration would naturally be the $\Fq$ points of an inductive
limit of schemes of finite type, with components of arbitrarily large
dimension.  Moreover, the relevant sheaf on this ind-scheme would have
stalks of arbitrarily large rank.  It is not at all clear how to
handle this situation.

\section{Preliminaries on $L$-functions}\label{s:L-funs}

\subsection{Input data and hypotheses}\label{ss:data}

The notation and hypotheses the following paragraphs
(\ref{sss:data-begin} through \ref{sss:data-end}) will be in force for
the rest of the paper. 

\subsubsection{}\label{sss:data-begin}
Let $\Curve$ be a smooth, proper, geometrically irreducible curve over
the finite field $\Fq$ of characteristic $p$ and let $F=\Fq(\Curve)$
be its field of functions.  Choose an algebraic closure $F^{\rm alg}$
of $F$ and let $\overline{F}\subset F^{\rm alg}$ be the separable
closure of $F$.  Let $G=\gal(\overline{F}/F)$ be the absolute Galois
group of $F$.  For each place $v$ of $F$ we choose a decomposition
group $D_v\subset G$ and we let $I_v$ and $Fr_v$ be the corresponding
inertia group and geometric Frobenius class.  We write $\deg v$ for
the degree of $v$ and $q_v=q^{\deg v}$ for the cardinality of the
residue field at $v$.

For positive integers $n$ we write $\Fqn$ for the subfield of
$\overline{F}$ of cardinality $q^n$, $F_n$ for the compositum $\Fqn
F$, and $G_n\subset G$ for $\gal(\overline{F}/F_n)$.  We write $F_{\infty}$
for $\Fqbar F$ and $G_{\infty}$ for $\gal(\overline{F}/F_{\infty})$.

\subsubsection{}\label{sss:root-q}
Fix a prime $\ell\neq p$ and let $\Qlbar$ be an algebraic closure
of $\Ql$, the field of $\ell$-adic numbers.  Fix also imbeddings
$\overline\Q\into\C$ and $\overline\Q\into\Qlbar$ and a compatible
isomorphism $\iota:\Qlbar\to\C$.  Whenever a square root of $q$ is
needed in $\Qlbar$, we take the one mapping to the positive square
root of $q$ in $\C$.  Having made this choice, we can define Tate
twists by half integers.

\subsubsection{}\label{sss:rho}
Fix a continuous, absolutely irreducible representation
$\rho:\gal(\overline{F}/F)\to\GL_r(E)$ where $E$ is a finite extension
of $\Ql$ in $\Qlbar$.  (We may extend the coefficient field $E$ as
necessary below.)  We assume that $\rho$ remains absolutely
irreducible when restricted to $G_{\infty}$ and that it is unramified outside
a finite set of places, so that it factors through
$\pi_1(U,\overline{\eta})$ for some non-empty open subscheme
$j:U\into\Curve$.  (Here $\overline{\eta}$ is the geometric point of
$\Curve$ defined by the fixed embedding $F\into F^{\rm alg}$.)  By
\cite[VII.6]{Lafforgue} and \cite[1.2.8-10]{WeilII}, $\rho$ is
$\iota$-pure of some weight $w$, i.e., for every place $v$ where
$\rho$ is unramified, each eigenvalue $\alpha$ of $\rho(Fr_v)$
satisfies $|\iota(\alpha)|=q_v^{w/2}$.  For convenience, we assume
that $w$ is an integer.

\subsubsection{}\label{sss:self-duality}
We say that a representation $\tau$ of $G_n$ is self-dual if it is isomorphic
to its contragredient.  This is equivalent to saying that there is a
non-degenerate $G_n$-equivariant bilinear pairing on the underlying
space.  If this pairing is symmetric, we say $\tau$ is
``orthogonally self-dual'' and that $\tau$  ``has sign
$+1$''.  If it is alternating, we say $\tau$ is
``symplectically self-dual'' and that $\tau$  ``has sign
$-1$''.  Schur's lemma implies that if an irreducible representation is
self-dual, then it is either orthogonally or symplectically self-dual.
Also, if a representation is symplectically self-dual, then its
degree is even.

If $\tau$ is
self-dual, then the weight $w$ of $\tau$ is 0.  To generalize
slightly, we say that $\tau$ is symplectically (orthogonally)
self-dual of weight $w$ if $\tau$ has weight $w$ and the Tate twist $\tau(w/2)$ is
symplectically (orthogonally) self-dual.  (Here $\tau(w/2)$ is
characterized by the equation $\tau(w/2)(Fr_v)=\tau(Fr_v)q_v^{-w/2}$.)

\subsubsection{}\label{sss:rho-duality-hyp}
Since $\rho$ is absolutely irreducible when restricted to $G_{\infty}$,
Schur's lemma implies that if $\rho$ is self-dual when restricted to
$G_{\infty}$, then for any integer $w$ there is a character of
$G/G_{\infty}\cong\gal(\Fqbar/\Fq)$ such that $\rho\tensor\chi$ is self-dual
of weight $w$.

We always assume that if $\rho$ is self-dual when restricted to $G_{\infty}$,
then it is already self-dual of some integer weight $w$ as a
representation of $G$.  In light of the above, this is not a serious
restriction.

\subsubsection{}
For each place $v$ of $F$ we write $\cond_v\rho$ for the exponent of
the Artin conductor of $\rho$ at $v$.  (See~\cite[Chap.~VI]{SerreLF}
for definitions.)  We let $\n=\sum_v({\cond_v\rho})[v]$ be the global
Artin conductor of $\rho$, viewed as an effective divisor on $\Curve$.
We write $|\n|$ for the support of $\n$, i.e., for the set of places
of $F$ where $\rho$ is ramified.

\subsubsection{}\label{sss:L-funs}
Attached to $\rho$ we have an $L$-function, defined formally by the
product 
\begin{align*}
L(\rho,F,T)&=\prod_v\det\left(1-\rho(Fr_v)T^{\deg v}
                  \left|(E^r)^{\rho(I_v)}\right.\right)^{-1}.
\end{align*}

The Grothendieck-Lefschetz trace formula implies that $L(\rho,F,T)$ is
actually a rational function of $T$.  More precisely, if $\rho$ is the
trivial representation, $L(\rho,F,T)$ is just the $Z$-function of $F$
(so $L(\rho,F,q^{-s})=\zeta(\Curve,s)$) and if $\rho$ is geometrically
non-trivial, i.e., non-trivial when restricted to $G_{\infty}$, then
$L(\rho,F,T)$ is a polynomial in $T$ of degree
$N=(2g_\Curve-2)(\deg\rho)+\deg\n$.

Writing the numerator of $L(\rho,F,T)$ as $\prod(1-\beta_iT)$, we call
the $\beta_i$ the inverse roots of $L(\rho,F,T)$.  Deligne's purity
result \cite[3.2.3]{WeilII} says that the inverse roots of
$L(\rho,F,T)$ have $\iota$-weight $w+1$, i.e.,
$|\iota(\beta_i)|=q^{w+1}$ for all $i$.

If $K$ is a finite extension of $F$ contained in $\overline{F}$, we
abbreviate $L(\rho|_{\gal(\overline{F}/K)},K,T)$ to $L(\rho,K,T)$.

Using the embedding $E\into\Qlbar\cong\C$ we may view $L(\rho,F,T)$ as
a rational function in $T$ with complex coefficients.  Then the
$L$-function appearing in the Introduction is $L(\rho,F,q^{-s})$.

\subsubsection{}\label{sss:local-conditions}
Fix a positive integer $d$ prime to $p$.  We let $a=[\Fq(\mu_d):\Fq]$
where $\mu_d$ denotes the $d$-th roots of unity.  
If necessary, we expand the coefficient field $E$ so that it contains
the $d$-th roots of unity and a square root of $q$.

Fix also three finite sets of places of $F$ called $S_s$, $S_i$,
$S_r$, which are pairwise disjoint.

We will be considering extensions of $F_n$ of the form
$K=F_n(f^{1/d})$ where $f\in F_n^\times$ and where the places of $F_n$
over $S_s$, $S_i$, and $S_r$ are split, inert, or ramified in $K$.
More precisely:

\begin{defn}\label{defn:local-conds}
  We say that {\it $f$ satisfies the local conditions\/} or {\it
    $K=F_n(f^{1/d})$
    satisfies the local conditions\/} if the
    following hold:
\begin{enumerate}
\item[(a)] for every place $v$ of $F_n$ over $S_s$, there is a place of $K$
  over $v$ unramified and of residue degree 1;
\item[(b)] for every place $v$ of $F_n$ over $S_i$, there is a place of $K$
  over $v$ unramified and of largest possible residue degree, namely
  $\gcd(d,q_v-1)$;
\item[(c)] every place of $F_n$ over $S_r$ is totally ramified in $K$; and
\item[(d)] every place of $F_n$ over $|\n|\setminus S_r$ is unramified in
  $K$.
\end{enumerate}
\end{defn}

\subsubsection{}\label{sss:alphas}\label{sss:data-end}
The last piece of data we need is a sequence of algebraic numbers
$\alpha_n$, indexed by positive integers $n$, which we view as
elements of $\Qlbar$ via the fixed embedding $\Qbar\into\Qlbar$.  We
assume that the image of $\alpha_n$ under 
$\iota:\Qlbar\to\C$ has absolute value $q^{n(w+1)/2}$.

\subsection{Base change and twisting}\label{ss:twisting}
Our main theorem is a statement about the existence of $f\in
F_n^\times$ such that $L(\rho,F_n(f^{1/d}),T)$ has no higher order of
zero at $T=\alpha_n^{-1}$ than $L(\rho,F_n,T)$ does.  If $F_n$
contains the $d$-th roots of unity then $F_n(f^{1/d})$ is a Kummer
extension of $F_n$ and we have a factorization
\begin{equation*}
L(\rho,F_n(f^{1/d}),T)=\prod_{i=0}^{d-1}L(\rho\tensor\chi_f^i,F_n,T)
\end{equation*}
where $\chi_f$ is a character of order $d$ of $G_n$ trivial on
$\gal(\overline{F}/F_n(f^{1/d}))$.  (Here and elsewhere we write
$\rho\tensor\chi^i_f$ for what should properly be denoted
$\rho|_{G_n}\tensor\chi^i_f$.)  Thus in this case the main theorem is a
non-vanishing statement for abelian twists of $\rho$.  The purpose of
this subsection is to record a similar factorization valid without the
assumption that $F_n$ contains the $d$-th roots of unity.  This will
relate the main theorem to a statement about non-vanishing of certain
non-abelian twists.

\subsubsection{}\label{sss:base-change}
Let $f$ be an element of $F_n^\times$ which is not an $e$-th power for
any divisor $e>1$ of $d$ and choose a $d$-th root $f^{1/d}$ of $f$ in
$\overline F$.  Set $G_n=\gal(\overline{F}/F_n)$,
$H_{n,f}=\gal(\overline{F}/F_n(f^{1/d}))$,
$G_{na}=\gal(\overline{F}/F_n(\mu_d))$, and
$L_{n,f}=\gal(\overline{F}/F_n(\mu_d,f^{1/d}))$. Here is the diagram
of fields:
\begin{equation*}
\xymatrix{&F_n(\mu_d,f^{1/d})=F_{na}(f^{1/d})\ar@{-}[dl]\ar@{-}[dr]&\\
F_n(f^{1/d})\ar@{-}[dr]&&F_n(\mu_d)=F_{na}\ar@{-}[dl]\\
&F_n&}
\end{equation*}
and the corresponding diagram of Galois groups:
\begin{equation*}
\xymatrix{&L_{n,f}\ar@{-}[dl]\ar@{-}[dr]&\\
H_{n,f}\ar@{-}[dr]&&G_{na}\ar@{-}[dl]\\
&G_n&}
\end{equation*}
Clearly $G_{na}$ and $L_{n,f}$ are normal subgroups of $G_n$ and
$H_{n,f}$ is a (possibly non-normal) subgroup of $G_n$ of index $d$.
Moreover, $G_n/L_{n,f}$ is a semi-direct product:
\begin{equation*}
 G_n/L_{n,f}=G_{na}/L_{n,f}\sdp H_{n,f}/L_{n,f}\cong
G_{na}/L_{n,f}\sdp G_n/G_{na}.
\end{equation*} 
Fix an isomorphism $\mu_d(\overline{F})\isoto\mu_d(E)$
and let $\chi_f$ be the $E$-valued character of $G_{na}$ of order $d$
given by the natural isomorphism
$G_{na}/L_{n,f}\isoto\mu_d(\overline{F})$ (namely
$\sigma\mapsto\sigma(f^{1/d})/f^{1/d}$) followed by
$\mu_d(\overline{F})\isoto\mu_d(E)$.  Note that $\chi_f^i$ is in fact
well-defined on the possibly larger group
$\gal(\overline{F}/F_n(\mu_{d/\gcd(d,i)}))$.

\begin{lemma}\label{lemma:Frob-twist}
If $\Phi\in G$ lies over the geometric Frobenius in $\gal(\Fqbar/\Fq)$ and
$(\chi^i_f)^\Phi$ is defined by $(\chi^i_f)^\Phi(h)=\chi^i_f(\Phi h\Phi^{-1})$,
then $(\chi^i_f)^\Phi=\chi^{iq}_f$.
\end{lemma}

\begin{proof}
This is an easy consequence of the definitions.
\end{proof}

\subsubsection{}
Our notational convention in \ref{sss:L-funs} says that
\begin{equation*}
L(\rho,F_n(f^{1/d}),T)=L(\res^{G_n}_{H_{n,f}}\rho,F_n(f^{1/d}),T)
\end{equation*}
and by standard properties of $L$-functions (e.g.,
\cite[3.8.2]{DeligneConstants}),
\begin{equation*}
L(\res^{G_n}_{H_{n,f}}\rho,F_n(f^{1/d}),T)=
L(\ind^{G_n}_{H_{n,f}}\res^{G_n}_{H_{n,f}}\rho,F_n,T).
\end{equation*}
Also, $\ind^{G_n}_{H_{n,f}}\res^{G_n}_{H_{n,f}}\rho
\cong\rho\tensor\ind^{G_n}_{H_{n,f}}\boldsymbol1$ where we write
$\boldsymbol1$ for the trivial representation with coefficients in $E$
of $H_{n,f}$ (\cite[3.3 Example 5]{SerreLR}).  It is well-known that
$\ind^{G_n}_{H_{n,f}}\boldsymbol1$ is the linear representation
associated to the permutation action of $G_n$ on the coset space
$G_n/H_{n,f}$.  We need to know how this representation factors into
irreducibles.

\begin{lemma}\label{lemma:reps}
Let $\sigma_f=\ind^{G_n}_{H_{n,f}}\boldsymbol1$.  Then the irreducible
constituents of $\sigma_f$ are in bijection with the orbits of
multiplication by $q^n$ on $\Z/d\Z$.  For each orbit $o\subset
\Z/d\Z$, set $d_o=d/gcd(d,i)$ for any $i\in o$ and set
$a_o=[F_n(\mu_{d_o}):F_n]=\#o$.  Then the representation
$\sigma_{o,f}$ corresponding to the orbit $o$ has dimension $a_o$ and
the restriction of $\sigma_{o,f}$ to
$G_{na_o}=\gal(\overline{F}/F_n(\mu_{d_o}))$ splits into lines; more
precisely, $\sigma_{o,f}|_{G_{na_o}}\cong\oplus_{i\in o}\chi_f^i$.
\end{lemma}

For example, when $\mu_d\subset F_n$, i.e., $q^n\equiv1\pmod d$, all
the orbits $o$ are singletons and we have
$\sigma_f\cong\oplus_{i\in\Z/d\Z}\chi_f^i$ as representations of
$G_n$.

\begin{proof}
  This is a standard exercise in representation theory.  Indeed, by
  \cite[7.3]{SerreLR}, the restriction of $\sigma_f$ to $G_{na}$ is
  $\ind_{L_{n,f}}^{G_{na}}\boldsymbol1$ which is easily seen to be
  $\oplus_{i\in\Z/d\Z}\chi_f^i$.  (To apply \cite{SerreLR}, we should
  note that all the representations in question are trivial on
  $L_{n,f}$ and so we are really working with subgroups of the finite
  group $G_n/L_{n,f}$.)  The factors $\chi_f^i$ are permuted by $G_n$
  and Lemma~\ref{lemma:Frob-twist} shows that under the right action
  $(\chi_f^i)^g(h)=\chi^i_f(ghg^{-1})$, we have
  $(\chi_f^i)^g=\chi_f^{iq^{n(g)}}$ where $q^{n(g)}$ is the image of
  $g$ under the natural map $G_n\to G_n/G_{na}\subset(\Z/d\Z)^\times$;
  the image of this map is the cyclic subgroup of $(\Z/d\Z)^\times$
  generated by $q^n$.  This proves that $\sigma_{o,f}=\oplus_{i\in
    o}\chi_f^i$ is an irreducible constituent of $\sigma_f$ and it is
  clear that $\sigma_{o,f}$ splits into lines when restricted to
  $G_{na_o}$.
\end{proof}

Thus, the general analogue of the factorization at the beginning of
this section is
\begin{equation}\label{eqn:L-factors}
L(\rho,F_n(f^{1/d}),T)=\prod_{o\subset\Z/d\Z}L(\rho\tensor\sigma_{o,f},F_n,T)
\end{equation}
where the product is over the orbits of $q^n$ on $\Z/d\Z$.  If we
assume that $n$ is relatively prime to $a=[\Fq(\mu_d):\Fq]$ then the
orbits for multiplication by $q^n$ are the same as the orbits for
multiplication by $q$.

It will be useful to know that $\sigma_{o,f}$ is itself induced.
Recall that $d_o=d/\gcd(i,d)$ for any $i\in o$ and
$a_o=[F_n(\mu_{d_o}):F_n]=\#o$. 

\begin{lemma}\label{lemma:induced}
 $\sigma_{o,f}\cong\ind_{G_{na_o}}^{G_n}\chi^i_f$ for any
$i\in o$. 
\end{lemma}

\begin{proof}
This is immediate from the facts that
$\sigma_{o,f}|_{G_{na_o}}\cong\oplus_{i\in o}\chi_f^i$ and that
$G_n/G_{na_o}$ permutes the factors $\chi^i_f$ simply transitively.
\end{proof}

Here is a criterion for $\sigma_{o,f}$ to be self-dual.

\begin{lemma}\label{lemma:sigma-self-dual}
$\sigma_{o,f}$ admits a non-degenerate $G_n$-invariant pairing if and
only if $-o=o$, i.e., if and only if $\{-i\,|\,i\in o\}=o$.
In this case, the pairing is symmetric.
\end{lemma}

\begin{proof}
If $-o=o$ it is easy to write down explicitly a $G_n$-invariant
symmetric pairing.  Indeed, this is obvious if $o=\{d/2\}$ and so
$\sigma_{o,f}=\chi^{d/2}_f$.  Otherwise, $a_o=\#o$ is even and 
by the previous lemma $\sigma_{o,f}$ can be realized as
\begin{equation*}
\{\phi:G_n\to E\,|\,\phi(gh)=\chi_f^i(h)\phi(g)\text{ for
all }h\in G_{na_o}\}
\end{equation*}
for any fixed $i\in o$.  The $G_n$ action is given by
$(g\phi)(g')=\phi(g^{-1}g')$.  Let $g$ be a generator of
$\gal(F_n(\mu_{d_o})/F_n)\cong\Z/a_o\Z$.  A suitable
pairing is then given by
\begin{equation*}
\<\phi_1,\phi_2\>=\sum_{j=0}^{a_o-1}\phi_1(g^j)\phi_2(g^{j+a_o/2}).
\end{equation*}
(To check the $G_n$-invariance, one uses that
$\gal(F_n(\mu_{d_o},f^{i/d})/F_n)$ is the semi-direct product 
\begin{equation*}
\gal(F_n(\mu_{d_o},f^{i/d})/F_n(\mu_{d_o})) \sdp
\gal(F_n(\mu_{d_o})/F_n)
\end{equation*}
and that $g^{a_o/2}$
acts by inversion on
$\gal(F_n(\mu_{d_o},f^{i/d})/F_n(\mu_{d_o}))=G_{na_o}/\ker(\chi^i_f)$.)
Since $\sigma_{o,f}$ is irreducible, any other $G_n$-invariant pairing
is a scalar multiple of this one, so is symmetric.

Conversely, if $-o\neq o$, then $\sigma_{o,f}$ is visibly not
self-dual when restricted to $G_{na_o}$ (where it is isomorphic to
$\oplus_{i\in o}\chi^i_f$).  Thus it cannot be self-dual as a
representation of $G_n$.
\end{proof}

\section{Forced zeroes}\label{s:forced-zeroes}

\subsection{Functional equations and forced zeroes}\label{ss:funl-eqns}

As a step toward stating a more precise version of the main theorem,
we review some well-known facts about the functional equations
satisfied by $L(\rho,F,T)$ and its twists.

\subsubsection{}
Let $\tau:G_n\to\GL_r(E)$ be an absolutely irreducible representation
of $G_n$ which is unramified outside a finite set of places and is
$\iota$-pure of some weight $w$.  (In the applications, $\tau$ will be $\rho$
or one of its twists $\rho\tensor\sigma_{o,f}$.)

\subsubsection{}\label{sss:funl-eqn}
The $L$-functions $L(\tau,F_n,T)$ satisfy functional equations.  If
$\tau$ has weight $w$, then we have
\begin{equation}\label{eqn:funl-eqn}
L(\tau,F_n,T)=W(\tau,F_n)q^{n\frac{(w+1)}{2}N}T^N
L(\tau\check{\ },F_n,(q^{n(w+1)}T)^{-1})
\end{equation}
where $W(\tau,F_n)$ is an algebraic number of weight 0,
$N=(2g_\Curve-2)(\dim\tau)+\deg\cond(\tau)$, and $\tau\check{\ }$ is the
contragredient of $\tau$.  (If $n(w+1)N$ is odd, we take the positive
square root of $q$, or more precisely, the square root of $q$ in $E$
which maps to the positive square root of $q$ under $\iota$.  Although
we have omitted it from the notation, in general $W(\tau,F_n)$ depends
on $\iota$, via the choice of square root of $q$.) 

If $\tau$ is the trivial representation, $W(\tau,F_n)=1$.  If $\tau$
is geometrically non-trivial and the inverse roots of $L(\tau,F_n,T)$
are $\beta_1,\dots,\beta_N$, then
$W(\tau,F_n)q^{n\frac{(w+1)}{2}N}$ can be described succinctly (and
independently of $\iota$) as $\prod_{i=1}^N(-\beta_i)$.  This implies
that $W(\tau,F_{nm})=(-1)^{N(m+1)}W(\tau,F_n)^m$.

\subsubsection{}
As is well-known, functional equations sometimes force zeroes of
$L$-functions at certain values of $s$ or $T=q^{-s}$.  In the
remainder of this subsection, we explain how this works out in the
function field situation (where the functional equation has {\it
  two\/} fixed points).

As usual, let $\tau$ be an absolutely irreducible representation of
$G_n$ of weight $w$ and consider the functional equation
\ref{eqn:funl-eqn}.  Note that the involution $T\mapsto
(q^{n(w+1)}T)^{-1}$ has two fixed points, namely $T=\pm
q^{-n(w+1)/2}$.  Thus, when $\tau$ is self-dual the functional
equation may force $\pm q^{n(w+1)/2}$ as inverse roots of
$L(\tau,F_n,T)$.  Here is the precise statement, which we leave as a
simple exercise for the reader.

\begin{lemma}\label{lemma:quad-zeroes}
  Suppose that $\tau$ is geometrically non-trivial and self-dual of
  weight $w$ and let $N$ be the degree, as a polynomial in $T$, of
  $L(\tau,F_n,T)$.
\begin{enumerate}
\item If $N$ is even and $W(\tau,F_n)=-1$, then $\pm q^{n(w+1)/2}$ are
both inverse roots of $L(\tau,F_n,T)$.
\item If $N$ is odd, then $-W(\tau,F_n)q^{n(w+1)/2}$ is an
inverse root of $L(\tau,F_n,T)$.
\end{enumerate}
\end{lemma}

The lemma applies only if $\tau$ is symplectically self-dual.  Indeed,
when $\tau$ is orthogonally self-dual, $W(\tau,F_n)=1$ and $N$ is even
(see~\ref{sss:coh-and-Ls}) and when $\tau$ is not self-dual, the
functional equation does not force any inverse roots since in that
case the functional equation relates two different $L$-functions.

\subsubsection{}
Let $\Psi_n\subset F_n^\times$ be the set of functions $f$ such that
the quadratic extension $F_n(\sqrt{f})/F_n$ satisfies the local conditions
(i.e., it is split at places of $F_n$ over $S_s$, inert at places of
$F_n$ over $S_i$, totally ramified at places of $F_n$ over $S_r$, and is
unramified at all places of $F_n$ over $|\n|\setminus S_r$).  Note
that we make no restrictions at places not over $S_s\cup S_i\cup
S_r\cup |\n|$ and so $\Psi_n$ is an infinite set.

If $f\in F_n^\times$ we write $\psi_f$ for the character of
$\gal(\overline{F}/F_n)$ corresponding to the quadratic extension
$F_n(\sqrt{f})/F_n$.  (In the notation of the previous section, this would
be $\sigma_{\{d/2\},f}=\chi_f^{d/2}$.)

As we have seen, functional equations can force certain numbers
$\alpha$ to be inverse roots of the twisted $L$-functions
$L(\rho\tensor\psi_f,F_n,T)$ for many $\psi_f$.  In order to control
this situation, we need to analyze when the signs
$W(\rho\tensor\psi_f,F_n)$ are fixed as $f$ varies over $\Psi_n$.  To
do so, we need to collect some facts about the signs $W$.

\subsubsection{}
It will be important for us that the sign $W(\tau,F_n)$ in the
functional equation admits an expression as a product of local
factors.  (See \cite[9.9]{DeligneConstants}, \cite[3.4]{Tate}, or
\cite[3.2.1.1]{Laumon} for more details; \cite{Laumon} treats the case
where the representation $\rho$ need not {\it a priori\/} be part of a
compatible family.)  In general, one must make auxiliary choices of a
measure and an additive character to define these local factors, but
in case $\tau$ is symplectically self-dual, the local factors are
independent of these choices.  Since this is the only case we need, we
assume for the rest of this subsection that $\tau$ is symplectically
self-dual.

Under that assumption, there are local factors $W_v(\tau,F_n)=\pm1$
which depend only on the restriction of $\tau$ to $D_v$ (and $\iota$)
and which are 1 wherever $\tau$ is unramified.  The global sign is
then given by $W(\tau,F_n)=\prod_vW_v(\tau,F_n)$.  We need to know how
these local factors behave under quadratic twists.

\begin{lemma}\label{lemma:signs}
Suppose that $p=char(F)$ is odd, $\tau$ is a symplectically self-dual
representation of $G$ of some weight $w$, and
$\psi$ is a quadratic character of $G_n$.
\begin{enumerate}
\item If $\tau$ is unramified at $v$ and $\psi$ is ramified at $v$,
then 
\begin{align*}
W_v(\tau\tensor\psi,F_n)&=(-1)^{(q_v-1)(\dim\tau)/4}\\
&=(-1)^{(\deg v)(q^n-1)(\dim\tau)/4}.
\end{align*}
\item If $\tau$ is ramified at $v$ and $\psi$ is unramified and
non-trivial at $v$, then
\begin{equation*}
W_v(\tau\tensor\psi,F_n)=(-1)^{\cond_v(\tau)}W_v(\tau,F_n).
\end{equation*}
\end{enumerate}
\end{lemma}

\begin{proof}
In case 1, $\tau\tensor\psi|_{D_v}$ is a direct sum of 1-dimensional
representations and we can compute the value of $W_v$ using classical
results on Gauss sums.  We leave the details as an exercise.

In case 2, \cite[5.5.1]{DeligneConstants}  or \cite[3.4.6]{Tate} says that
\begin{equation*}
W_v(\tau\tensor\psi,F_n)=
\psi(\pi_v)^{\cond_v(\tau)}W_v(\tau,F_n).
\end{equation*}
(Here we use that $\tau$ is symplectic and so $\deg(\tau)$ is even.) 
But our assumptions imply that $\psi(\pi_v)=-1$.
\end{proof}

\subsubsection{}
For the rest of this subsection, we assume:
\begin{equation*}
\hyp{$\rho$ is symplectically self-dual of some weight
$w$}
\end{equation*}

This hypothesis implies that the dimension of $\rho$ is
even and using this, it is not hard to check that the parity of the
degree of $\cond(\rho\tensor\psi_f)$ is the same for all $f\in\Psi_n$.
Since the degree in $T$ of $L(\rho\tensor\psi_f,F_n,T)$ is
$N=(2g_\Curve-2)(\deg\rho)+\deg\cond(\rho\tensor\psi_f)$, the parity of
$N$ is independent of the choice of $f\in\Psi_n$.

We now discuss a (local) hypothesis which determines whether the sign
$W(\rho\tensor\psi_f,F_n)$ is the same for all $f\in\Psi_n$ or whether
it varies.  Note that for all $v$ over $|\n|$, $\cond_v(\rho\tensor\psi_f)$
is independent of the choice of $f\in\Psi_n$.
\begin{equation}
\hyp{For every place $v$ of $F_n$ over $|\n|\setminus(S_s\cup S_i)$,
$\cond_v(\rho\tensor\psi_f)$ is even for one (and thus every)
$f\in\Psi_n$.}\label{eqn:cond1}
\end{equation}

\begin{lemma}\label{lemma:global-signs}
  If hypothesis \ref{eqn:cond1} is satisfied then the signs
  $W(\rho\tensor\psi_f,F_n)$ for $f\in\Psi_n$ are all the same.  On
  the other hand, if hypothesis \ref{eqn:cond1} fails 
  then $W(\rho\tensor\psi_f,F_n)$ takes both
  values $\pm1$ as $f$ varies through $\Psi_n$.
\end{lemma}

\begin{proof}
Recall that we have assumed that $\rho$ is symplectically
self-dual.  If $f\in\Psi_n$ then for places $v$ of $F_n$ not over
$|\n|$, part 1 of Lemma~\ref{lemma:signs} tells us that
\begin{equation*}
W_v(\rho\tensor\psi_f,F_n)=\begin{cases}
1&\text{if $\psi_f$ is unramified at $v$}\\
(-1)^{(\deg v)(q-1)(\deg\rho)/4}&\text{if $\psi_f$ is ramified at $v$.}
\end{cases}
\end{equation*}
But the sum of $\deg v$ over places of $F_n$ which are not over $|\n|$
and where $\psi_f$ is ramified has fixed parity independent of $f$.
Indeed 
\begin{equation*}
\sum_{v\text{ over }|{\cond(\psi_f)}|\setminus|\n|}\deg v\equiv
\sum_{v\text{ over }|{\n}|\cap S_r}\deg v\pmod2
\end{equation*}
since $\cond(\psi_f)$ has even degree.  Thus
the sign $\prod_{v\text{ not over }|\n|}W_v(\rho\tensor\psi_f,F_n)$ is
independent of $f\in\Psi_n$.

Now take $f$ and $f'$ in $\Psi_n$.  If $v$ is over $|\n|$, then
$\psi''=\psi_{f'}/\psi_f$ is unramified at $v$ and it is trivial on
$D_v$ if $v$ is over $S_s\cup S_i$.  Applying part 2 of
Lemma~\ref{lemma:signs} at those places $v$ over
$|\n|\setminus(S_s\cup S_i)$ where $\psi''$ is non-trivial (with
$\tau$ replaced by $\rho\tensor\psi$ and $\psi$ replaced by $\psi''$),
we conclude that
\begin{equation*}
\frac{W(\rho\tensor\psi_f,F_n)}{W(\rho\tensor\psi_{f'},F_n)}=
\prod_{\substack{v\text{ over }|\n|\setminus(S_s\cup S_i)\\
                    \psi''\text{ non-trivial on }D_v}}
(-1)^{\cond_v(\rho\tensor\psi_f)}
\end{equation*}

Hypothesis \ref{eqn:cond1} implies that this quantity is 1.  If
\ref{eqn:cond1} fails, by the Riemann-Roch theorem, we can choose $f$
and $f'$ in $\Psi_n$ so that this quantity takes both values $\pm1$.
(See Section~\ref{s:local-conds} for more details and a
quantitative statement about the density of such $f'$.) 
\end{proof}

\subsubsection{}\label{sss:quad-zeroes}
A common situation where $\tau=\rho\tensor\sigma_{o,f}$ is
symplectically self-dual is when $\rho$ is symplectically self-dual
and $-o=o$, so that $\sigma_{o,f}$ is orthogonally self-dual and the
tensor product is symplectically self-dual.  In particular, the
results of the preceding subsections are relevant to the special case
where $\sigma_{o,f}$ is a quadratic character $\psi_f$, i.e., when
$o=\{d/2\}$.

Recall that $N$, the degree of $L(\rho\tensor\psi_f,F_n,T)$ as a
polynomial in $T$, is given by 
\begin{equation*}
N=(2g_\Curve-2)(\deg\rho)+\deg\cond(\rho\tensor\psi_f).
\end{equation*}
Lemma~\ref{lemma:global-signs} and
Lemma~\ref{lemma:quad-zeroes} imply that $L(\rho\tensor\psi_f,F_n,T)$
has certain predictable inverse roots, as $f$ varies over the set
$\Psi_n$, in the following two situations:
\begin{enumerate}[(i)]
\item if $\rho$ is symplectically self-dual, the hypothesis
\ref{eqn:cond1} is satisfied, $N$ is even, and
$W(\rho\tensor\psi_f,F_n)=-1$ for one (and thus all) $f\in\Psi_n$,
then $\alpha=\pm q^{n(w+1)/2}$ are both inverse roots of
$L(\rho\tensor\psi_f,F_n,T)$
\item if $\rho$ is symplectically self-dual, the hypothesis
\ref{eqn:cond1} is satisfied, and $N$ is odd, then
$\alpha=-W(\rho\tensor\psi_f,F_n)q^{n(w+1)/2}$ is an inverse root of
$L(\rho\tensor\psi_f,F_n,T)$ for all $f\in\Psi_n$.
\end{enumerate}

\subsection{Zeroes forced by induction}\label{ss:induced-zeroes}
It turns out that there is another source of forced inverse roots of
$L$-functions, not visible via functional equations, coming from the
fact that $\sigma_{o,f}$ is an induced representation.  Here is the
precise statement:

\begin{prop}\label{prop:induced-zeroes}
Let $F$, $\rho$, $d$, and $n>0$ be as in \ref{ss:data}.  Fix an orbit
$o$ of multiplication by $q^n$ on $\Z/d\Z$ and set as usual
$d_o=d/\gcd(d,i)$ for any $i\in o$ and
$a_o=\#o=[F_n(\mu_{d_o}):F_n]$.  Assume that $-o=o$ and $a_o>1$.  Fix
$f\in F_n^\times$ and assume that the degree of
$\cond(\rho\tensor\chi^i_f)$ is odd for one (and thus every) $i\in
o$.  Then
\begin{enumerate}
\item if $\rho$ is symplectically self-dual of weight $w$, then
$1-\left(Tq^{n\frac{w+1}2}\right)^{a_o}$ divides
$L(\rho\tensor\sigma_{o,f},F_n,T)$.
\item if $\rho$ is orthogonally self-dual of weight $w$, then
$1+\left(Tq^{n\frac{w+1}2}\right)^{a_o}$ divides
$L(\rho\tensor\sigma_{o,f},F_n,T)$.
\end{enumerate}
\end{prop}

Note that the asserted inverse roots of the $L$-function are not
fixed points of the involution $T\mapsto (q^{n(w+1)}T)^{-1}$ in part
(2), and not all of them are fixed points in part (1) as soon
as $a_o>2$.  

We delay the proof of the proposition until
\ref{sss:induced-zeroes-proof} below, where it can be most naturally
explained in terms of cohomology.  For the moment, we just check the
assertion that $\cond(\rho\tensor\chi^i_f)$ is odd for all $i\in o$ if
it is so for one $i\in o$.  In fact, if $\Phi$ denotes an element of
$\gal(\overline F/F)$ which induces the $q$-power Frobenius on
$\Fqbar$ and $(\rho\tensor\chi^i_f)^{\Phi^n}$ is defined by
\begin{equation*}
(\rho\tensor\chi^i_f)^{\Phi^n}(g)=(\rho\tensor\chi^i_f)(\Phi^{n}g\Phi^{-n})
\end{equation*}
then
\begin{equation*}
\cond(\rho\tensor\chi^{iq^{n}}_f)=
\cond\left((\rho\tensor\chi^i_f)^{\Phi^n}\right)
\end{equation*}
and so $\cond(\rho\tensor\chi^i_f)$
and $\cond(\rho\tensor\chi^{iq^{n}}_f)$ have the
same degree.

\subsubsection{}
Here is an example of forced zeroes ``in nature'' which can be treated
by elementary means.

Let $q$ be a prime power with $q\equiv2\pmod3$ and let $X$
be a smooth projective curve over $\Fq$ given as a 3-fold cover of
the projective line by the equation
\begin{equation*}
y^3=f(x)
\end{equation*}
where $f(x)$ is a rational function on the line.  Suppose that $X$ has
odd genus.   (This can be arranged, for example, by assuming that $f$
has $d$ simple zeroes and $d-1$ poles, one of which is double, the
others simple.)    Note that as we vary $f$ we get a large family of
curves.

The claim then is that the numerator of the zeta function (or rather
$Z$-function) $Z(X,T)$ is divisible by $1+qT^2$, i.e., it has
$\pm\sqrt{-q}$ as inverse zeroes.  (In terms of
$\zeta(X,s)=Z(X,q^{-s})$, we are claiming that there are zeroes at
$s=\frac12+\frac{\pi i}{2\ln q}$ and $s=\frac12+\frac{3\pi i}{2\ln q}$.)
Note that these inverse zeroes are {\it not\/} at fixed points of the
functional equation.

The claim can be seen by an elementary argument: observe that since
$q\equiv2\pmod3$, every element of $\Fq$ has a unique cube root and so
the number of points on $X$ over $\Fq$ is $q+1$.  A similar statement
applies for all odd degree extensions of $\Fq$.  This implies that the
set of inverse roots of the numerator of the $Z$-function is
invariant under $\alpha\mapsto-\alpha$.  It is also invariant under
$\alpha\mapsto q/\alpha$ and the product of the inverse roots is
$q^g$.  Since there are $2g$ inverse roots and $g$ is odd, it follows
that for some inverse root $\alpha$ we have $\alpha=-q/\alpha$, as
claimed.  (Thanks to Mike Zieve for supplying this argument.)

\subsubsection{}
We now return to the general analysis of forced zeroes.
We want to give a simple local criterion which determines whether the
condition ``$\cond(\rho\tensor\chi_f^i)$ is odd for $i\in o$'' holds
for a fixed $o$ and all $f$ satisfying the local conditions. 
Consider the following hypothesis: 
\begin{equation}
\hyp{
For all places $v$ of $F_n$ over $|\n|\cap S_r$, 
$\cond_v(\rho\tensor\chi_v)$ has fixed parity as 
$\chi_v$ varies over totally ramified characters of $D_v$ of order
exactly $d_o=d/\gcd(i,d)$.  Moreover,
\begin{equation*}
\sum_{v\in|n|\cap S_r}\cond_v(\rho\tensor\chi_v)\deg v +
\sum_{v\in|n|\setminus S_r}\cond_v(\rho)\deg v
\end{equation*}
is odd for some (and thus any) choice of totally ramified local
characters $\chi_v$.
}\label{eqn:cond2}
\end{equation}

\begin{prop}\label{prop:avoiding-induced-zeroes}
  Fix an orbit $o\neq\{d/2\}$, an $i\in o$, and an integer $n$ prime
  to $a_o$.  Then $\deg\cond(\rho\tensor\chi_f^i)$ is odd for all
  $f\in F_n^\times$ satisfying the local conditions
  if and only if $\rho$ is even dimensional and hypothesis
  \ref{eqn:cond2} is satisfied.
\end{prop}

\begin{proof}
The degree of the Artin conductor is 
\begin{equation*}
\sum_{v\text{ over }|\n|}\cond_v(\rho\tensor\chi_f^i)\deg v 
+\sum_{v\text{ not over }|\n|}\cond_v(\rho\tensor\chi_f^i)\deg v 
\end{equation*}
and we have
\begin{multline*}
\sum_{v\text{ over }|\n|}\cond_v(\rho\tensor\chi_f^i)\deg v \\
=\sum_{v\text{ over }|\n|\cap S_r}\cond_v(\rho\tensor\chi_f^i)\deg v 
+\sum_{v\text{ over }|\n|\setminus S_r}\cond_v(\rho)\deg v 
\end{multline*}
and
\begin{equation*}
\sum_{v\text{ not over }|\n|}\cond_v(\rho\tensor\chi_f^i)\deg v 
=\sum_{\substack{v\text{ not over }|\n|\\v(f)\not\equiv0\pmod{d_0}}}
\deg\rho\deg v. 
\end{equation*}
It is thus clear that if $\rho$ is even dimensional and hypothesis
\ref{eqn:cond2} is satisfied, then $\deg\cond(\rho\tensor\chi^i_f)$
is odd for all $f$ satisfying the local conditions.

For the converse, first assume that $\deg\rho$ is even and hypothesis
\ref{eqn:cond2} fails.  Choose local characters $\chi_v$ of order
exactly $d_o$ at each $v$ over $|\n|\cap S_r$ so that the sum
appearing in \ref{eqn:cond2} is even.  Then there is an element $f\in
F_n^\times$ satisfying the local conditions such that the local
component at $v$ of $\chi^i_f$ is the fixed $\chi_v$ for all $v$ over
$|\n|\cap S_r$.  (This is an easy consequence of the Riemann-Roch
theorem; we just need to fix (modulo $d_o$) the valuation of $f$ at
$v$ over $|\n|\cap S_r$. See Section~\ref{s:local-conds} below for a
more precise version of this result.)  For such $f$, it is clear that
$\deg\cond(\rho\tensor\chi^i_f)$ is even.

Finally, assume that $\deg\rho$ is odd.  Fix a divisor $D$ which is
the sum of the places over $S_r$, each with multiplicity one, plus a
sum of places of odd degree not over $S_s\cup S_i\cup S_r\cup|\n|$
also with multiplicity one; we insist that there should be at least 2
such places and that the degree of $D$ be sufficiently large, namely
greater than $2g-2$ plus the sum of the degrees of all places over
$S_s\cup S_i$.  Let $f$ be an element of $F_n^\times$ satisfying the
local conditions and such that the divisor of $f$ is $-D$ plus an
effective square free divisor.  (I.e., $f$ has polar divisor $D$ and
its zeroes are distinct.)  The existence of such an $f$ again follows
easily from the Riemann-Roch theorem.  We note that
$\deg\cond(\rho\tensor\chi^i_f)$ only depends on $D$, not on the
specific $f$ chosen.  If, for this $D$,
$\deg\cond(\rho\tensor\chi^i_f)$ is even, we are finished.  If not,
modify $D$ as follows: drop one place (of odd degree) not over
$S_s\cup S_i\cup S_r\cup|\n|$ and change the coefficient of another
place not over $S_s\cup S_i\cup S_r\cup|\n|$ from 1 to 2.  Calling the
resulting divisor $D'$, choose $f'\in F_n^\times$ with polar divisor
$D'$ and distinct zeroes which satisfies the local conditions.  Then
we have removed one term $\deg\rho\deg v$ from the last displayed sum
and not changed anything else (here we use that $i\neq d/2$) and so
\begin{equation*}
\deg\cond(\rho\tensor\chi^i_{f'})=
\deg\cond(\rho\tensor\chi^i_{f})-\deg\rho\deg v
\end{equation*}
 which is even.
\end{proof}

\begin{rem} 
In general it is not possible to find one $f$ which makes
$\deg\cond(\rho\tensor\chi^i_f)$ even for all $i$ in several different
orbits $o$.  By the proposition, we can always arrange this for one
orbit. If $\deg\rho$ is odd, we can find one $f$ which makes this true
for at least half of the orbits in any fixed collection of orbits (for
a fixed $d$).
\end{rem}

\subsubsection{}\label{sss:induced-zeroes}
The upshot of this subsection is that
$L(\rho\tensor\sigma_{o,f},F_n,T)$ has certain predictable inverse
roots for all $f$ satisfying the local conditions in the following two
situations:
\begin{enumerate}[(i)]
\item if $\rho$ is symplectically self-dual of weight $w$, $o=-o$,
      $a_o=\#o>1$, and hypothesis \ref{eqn:cond2} is satisfied, then
      the solutions $\alpha$ of $\alpha^{a_o}=q^{n(w+1)/2}$ are
      inverse roots of $L(\rho\tensor\sigma_{o,f},F_n,T)$.
    \item if $\rho$ is orthogonally self-dual of weight $w$ and of
      even degree, $o=-o$, $a_o=\#o>1$, and hypothesis \ref{eqn:cond2}
      is satisfied, then the solutions $\alpha$ of
      $\alpha^{a_o}=-q^{n(w+1)/2}$ are inverse roots of
      $L(\rho\tensor\sigma_{o,f},F_n,T)$.
\end{enumerate}

\section{Statement of the main technical theorem}\label{s:statement}

\subsection{Exceptional situations}
It will turn out that \ref{sss:quad-zeroes} and
\ref{sss:induced-zeroes} exhaust the supply of ``forced
zeroes'' in our situation.  The following definitions give a
convenient terminology for when forced zeroes occur:

\begin{defns}
  We say that $\rho$, $d$, $S_s$, $S_i$, $S_r$, $n$, $o\subset\Z/d\Z$,
  and $\alpha_n$ are ``exceptional'' (or more briefly ``$n$ is
  exceptional'') if one of the following conditions holds:
\begin{enumerate}[(i)]
\item $\rho$ is symplectically self-dual of weight $w$, $o=\{d/2\}$,
  the hypothesis \ref{eqn:cond1} is satisfied,
  $\deg\cond(\rho\tensor\chi_f)$ is even and $W(\rho\tensor\chi_f,F_n)=-1$
  for one (and thus all) $f\in\Phi_n$, and $\alpha_n=\pm q^{n(w+1)/2}$
\item $\rho$ is symplectically self-dual of weight $w$, $o=\{d/2\}$,
  the hypothesis \ref{eqn:cond1} is satisfied,
  $\deg\cond(\rho\tensor\chi_f)$ is odd and
  $\alpha_n=-W(\rho\tensor\chi_f,F_n)q^{n(w+1)/2}$ for one (and thus
  all) $f\in\Phi_n$
\item $\rho$ is symplectically self-dual of weight $w$, $-o=o$ and
  $a_o=\#o>1$, hypothesis~\ref{eqn:cond2} holds for $i\in o$, and
  $\alpha_n^{a_o}=q^{n\frac{w+1}2a_0}$
\item $\rho$ is even dimensional and orthogonally self-dual of weight
  $w$, $-o=o$ and $a_o=\#o>1$, hypothesis~\ref{eqn:cond2} holds for
  $i\in o$, and $\alpha_n^{a_o}=-q^{n\frac{w+1}2a_0}$
\item $\rho$ is symplectically self-dual of weight $w$, $-o=o$ and
  $a_o=\#o>1$, and $\alpha_n^{a_o}=q^{n\frac{w+1}2a_0}$
\item $\rho$ is orthogonally self-dual of weight $w$, $-o=o$ and
  $a_o=\#o>1$, and $\alpha_n^{a_o}=-q^{n\frac{w+1}2a_0}$
\end{enumerate}
\end{defns}

The exceptional situation (i) and (ii) arise in the context of
elliptic curves and ``Heegner conditions'' as we will see in
Section~\ref{s:elliptic-app} below.  Situations (iii) and (iv) are
related to the ``exotic'' forced zeroes of
Subsection~\ref{ss:induced-zeroes} and are also needed for the
application to elliptic curves.

Note that exceptional situations (iii) and (iv) are subsets of
situations (v) and (vi) respectively.  When we consider several orbits
$o$ at once, we will find that there is always a forced zero (of
multiplicity one) in situations (iii) and (iv), whereas in situations
(v) and (vi), we will only be able to assert that the
multiplicity of a zero at $\alpha_n$ is at most one.

We can now state the main theorem: 

\begin{thm}\label{thm:main}
Suppose that $F$, $\rho$, $d$, $S_s$, $S_i$, $S_r$, $(\alpha_n)$
satisfy the hypotheses of \ref{ss:data}.  Suppose also either that
$\rho$ is at worst tamely ramified at every place $v$ of $F$ or that
$p\ge\deg\rho+2$ and $\rho$ is tamely ramified at all places
$v\in|\n|\cap S_r$.
\begin{enumerate}
\item Fix an orbit $o\subset\Z/d\Z$ for multiplication by $q$ and set
$d_o=d/\gcd(d,i)$ for any $i\in o$ and $a_o=\#o=[F(\mu_{d_o}):F]$.
Then for all sufficiently large $n$ relatively prime to $a_o$, there
exists $f\in F_n^\times$ such that every place of $F_n$ over $S_s$
(resp.~$S_i$, $S_r$) splits (resp. is ``as inert as possible'', is totally
ramified) in $F_n(f^{1/d})$ and
\begin{itemize}
\item if $n$ is exceptional of type (i)-(iv),
$L(\rho\tensor\sigma_{o,f},F_n,T)$ has $\alpha_n$ as a simple inverse
root 
\item in all other cases, $\alpha_n$ is not an inverse root of 
$L(\rho\tensor\sigma_{o,f},F_n,T)$
\end{itemize}
\item Set $a=[F(\mu_d):F]$.  Then for all sufficiently large $n$
relatively prime to $a$, there exists $f\in F_n^\times$ such that 
every place of $F_n$ over $S_s$ (resp.~$S_i$, $S_r$) splits
(resp. is ``as inert as possible'', is totally ramified) in $F_n(f^{1/d})$ and 
for each orbit $o\subset\Z/d\Z$ for multiplication by $q$:
\begin{itemize}
\item if $n$ is exceptional of type (i)-(iv),
$L(\rho\tensor\sigma_{o,f},F_n,T)$ has $\alpha_n$ as a simple inverse
root
\item if $n$ is exceptional of type (v) or (vi),
$L(\rho\tensor\sigma_{o,f},F_n,T)$ has $\alpha_n$ as an inverse
root of multiplicity at most 1
\item in all other cases, $\alpha_n$ is not an inverse root of 
$L(\rho\tensor\sigma_{o,f},F_n,T)$
\end{itemize}
\end{enumerate}
\end{thm}

It is possible to get somewhat better control of the type (v) and (vi)
exceptional situations in various contexts.  For example, if
$\deg\rho$ is odd, we can find an $f$ so that
$L(\rho\tensor\sigma_{o,f},F_n,T)$ does not vanish at $\alpha_n$ for
at least half of the orbits $o$ of type (vi).  These improvements do
not seem likely to be of much use, so we omit them.

\subsection{Twisting}
Note that the truth of the theorem is invariant under twisting in the
following sense: the theorem holds for $F$, $\rho$, $d$, $S_s$, $S_i$,
$S_r$, $(\alpha_n)$ if and only if it holds for $F$, $\rho(t)$, $d$,
$S_s$, $S_i$, $S_r$, $(q^{-tn}\alpha_n)$.  (Here as in
\ref{sss:self-duality}, $\rho(t)$ is the Tate twisted representation,
characterized by $\rho(t)(Fr_v)=\rho(Fr_v)q_v^{-t}$.)   Thus by
twisting we may assume that $\rho$ has weight $w=-1$ and the
$\alpha_n$ all have $\iota$-weight 0, i.e., satisfy
$|\iota\alpha_n|=1$.  We make this assumption for the rest of the
paper.

                \section{Local conditions}\label{s:local-conds}

\subsection{Notational conventions}
The rest of this article will use more algebraic geometry.  We set the
following notations and conventions.

All schemes considered will be of finite type over $\spec\Fq$.  If $X$
is such a scheme and $k$ is an extension field of $\Fq$, we write
$X\times k$ for $X\times_{\spec\Fq}\spec k$.  Let $\Fqbar$ denote an
algebraic closure of $\Fq$.  We will often use a bar to denote the
base change to $\Fqbar$, so for example $\Cbar=\Curve\times\Fqbar$.

We write $Fr$ for the geometric ($q^{-1}$-power) Frobenius of $\Fqbar$
and its subfields, and also for the automorphism of
$\overline{X}=X\times\Fqbar$ which is the identity on $X$ and $Fr$ on
$\Fqbar$, and for its action on cohomology.

Suppose that $X$ is reduced and irreducible and let $\overline{\eta}$
be a geometric generic point of $X$ with residue field
$\kappa(\overline\eta)$.  To fix ideas, we take $\overline\eta$ to be
the spectrum of an algebraic closure of the field of rational
functions on $X$.  Let $\pi_1(X,\overline\eta)$ be the fundamental
group of $X$ with base point $\overline\eta$.  (See
\cite[Exp.~V]{SGA1}.)

Let $k$ be a finite extension of $\Fq$ and $x\in X(k)$ be a $k$-valued
point of $X$, i.e., a morphism $\spec k\to X$.  Choosing an algebraic
closure $\overline k$ of $k$ yields a geometric point
$\overline{x}:\spec\overline k\to X$ over $x$, from which we deduce an
embedding
\begin{equation*}
\gal(\overline k/k)\cong\pi_1(\spec k,\spec\overline k)
\into\pi_1(X,\overline x)\cong\pi_1(X,\overline\eta).
\end{equation*}
where the last isomorphism is a non-canonical ``path'' isomorphism.
We write $Fr_{k,x}$ for the image of the geometric Frobenius.  The
conjugacy class of $Fr_{k,x}$ is well-defined independently of the
choices.  When $k$ is a field with $q^n$ elements, we also write
$Fr_{n,x}$ for $Fr_{k,x}$.

Similarly, if $x$ is a closed point of $X$ with residue field
$\kappa(x)$, we may view $x$ as a $\kappa(x)$-valued point of $X$ and
form a Frobenius element $Fr_x=Fr_{\kappa(x),x}\in\pi_1(X,\overline
\eta)$ which is well-defined up to conjugation.  We will use this
notation mostly in the case where $X$ is a curve and $x$ is the closed
point associated to a place of the function field of $X$.

\subsection{The parameter space $X$}
We now introduce an effective $\Fq$-rational divisor $D$ on $\Curve$,
say $D=\sum_va_v[v]$ where $v$ runs over places of $F$ and the
coefficients $a_v$ are non-negative.  As usual, $\deg(D)=\sum_va_v\deg
v$ denotes the degree of $D$ and $|D|$ denotes the support of $D$,
i.e., the set of places where $a_v\neq0$.  We consider $D$ as a
divisor on the curves $\Curve\times\F_{q^n}$ in the natural way.  In
the course of the discussion we may enlarge $D$ so that its degree is
``sufficiently large'' in a sense which will be made precise as
needed.

Let $L$ be the scheme representing the functor on $\Fq$-algebras
\begin{equation*}
R\mapsto H^0(\Curve\times_{\spec\Fq}\spec R,\O(D))
=H^0(\Curve,\O(D))\tensor_\Fq R.
\end{equation*}
In concrete terms, this just means that $L$ is an
affine space over $\spec\Fq$ of dimension $\dim_\Fq
H^0(\Curve,\O(D))$.  Note as well that the set of $\Fq$ points
$L(\Fq)$ is what would classically be denoted $L(D)$.

Now let $X$ be the scheme which represents the functor $R\mapsto$
``the set of elements of $L(R)$ whose zeroes (as section of $\O(D)$)
are distinct and disjoint from $|D|\cup|\n|\cup S_s\cup S_i$.''  It is
clear what the quoted phrase means when $R$ is a field; the precise
meaning for a general scheme and a very detailed proof of the
existence of $X$ is explained in \cite[5.0.6, 6.0, and 6.1]{Katz}.
Among other things it is proven there that $X$ is an open subscheme of
$L$.  (Essentially, $X$ is obtained from $L$ by removing the
hyperplanes corresponding to sections of $\O(D)$ vanishing at some
point in $|D|\cup|\n|\cup S_s\cup S_i$ and a discriminant locus
corresponding to sections with multiple zeroes.)

\subsection{General local conditions}\label{ss:local-conditions}
In this subsection we make some general definitions which will allow
us to identify those points of $X(\Fqn)$ which satisfy various local
conditions needed in the proof of the main theorem.

Let us fix for each $n$ a finite set of places $S_n$ of $F_n$ and for
each place $w\in S_n$ a non-empty subset of $C_{n,w}\subset
F_{n,w}^\times/F_{n,w}^{\times d}$.  We define the degree of $S_n$ by
$\deg(S_n)=\sum_{w\in S_n}\deg(w)$ and we say that the collection
$(S_n,C_{n,w})$ is {\it compatible with\/} $D$ if the following condition
is satisfied: for every every $w$ in $S_n$, there
exists an element $f\in F_n^\times$ such that the order of pole
$-w(f)$ is equal to the coefficient of $w$ in $D$ and the class of $f$
in $F_{n,w}^\times/F_{n,w}^{\times d}$ lies in the subset $C_{n,w}$.

We say that $f\in F_n^\times$ {\it satisfies the local conditions
  imposed by $(S_n,C_{n,w})$\/} if for every $w\in S_n$, the class of
$f$ in $F_{n,w}^\times/F_{n,w}^{\times d}$ lies in the subset
$C_{n,w}$.  It is a consequence of the Riemann-Roch theorem that if
the degree of $D$ is sufficiently large (namely
$>2g_\Curve-2+\deg(S_n)$) and $(S_n,C_{n,w})$ is compatible with $D$,
then there are elements of $L(\Fqn)\subset F_n$ which satisfy the
local conditions imposed by $(S_n,C_{n,w})$.  The next proposition
tells us that the set of such elements which also lie in $X(\Fqn)$ has
a positive density, bounded away from 0.

Fix an effective divisor $D$ and a set of local conditions $(S_n,C_{n,w})$
for each $n$ which are compatible with $D$.  Define
\begin{equation*}
Y_n=\{f\in X(\Fqn)|\text{ $f$ satisfies the local conditions imposed
by $(S_n,C_{n,w})$}\}.
\end{equation*}

\begin{prop}\label{prop:good-pts-density}
  Assume that $\deg(D) > 2g_\Curve-2+\deg(S_n)$ for all $n$.  Then
  there exists a constant $C>0$, independent of $n$, such that
  \begin{equation*}
    \frac{\#Y_n}{\#X(\Fqn)}>C
  \end{equation*}
for all sufficiently large $n$.
\end{prop}

\begin{proof}
Let $B=\sup_n\deg(S_n)$, which is finite by hypothesis.
For each $n$, introduce an auxiliary effective divisor defined by
$D'_n=\sum_{w\in S_n}[w]$.  Note that $\deg(D'_n)=\deg(S_n)\le B$.  

Define $L(\Fqn)^{\text{good}}$ to be those elements of $L(\Fqn)$ which
satisfy the local conditions imposed by $(S_n,C_{n,w})$.  Whether
an element $f\in L(\Fqn)$ lies in $L(\Fqn)^{\text{good}}$ is
determined by the leading terms in the expansion of $f$ at places in
$|D_n'|=S_n$.  More precisely, note that for each $w\in S_n$ there is
a well-defined map
\begin{equation*}
  \frac{\O(D)_w}{\O(D-D'_n)_w}\setminus0\to F_{n,w}^\times/F_{n,w}^{\times d}
\end{equation*}
where $\O(D)_w$ and $\O(D-D'_n)_w$ are the stalks of $\O(D)$ and
$\O(D-D'_n)$ at $w$.
This map is not surjective, but its image does meet $C_{n,w}$ (this is the
definition of compatible) and its non-empty fibers all have
cardinality 
\begin{equation*}
\frac{(q_w-1)}{\gcd(q_w-1,d)}\ge\frac{(q_w-1)}{d}.
\end{equation*}
Let $C'_{n,w}$ be the subset of $\O(D)_w/\O(D-D'_n)_w$ consisting of
non-zero elements which map to $C_{n,w}$.  Then $C'_{n,w}$ is
non-empty and its ``density'' (i.e., its cardinality divided by that
of $\O(D)_w/\O(D-D'_n)_w$) is positive and bounded away
from 0 for all $n\gg0$.  (Indeed, it is bounded below by
$(q_w-1)/dq_w=1/d-1/dq_w$.)  Now let $C'_n$ be the subset $\prod_{w\in
  S_n}C'_{n,w}$ of
\begin{equation*}
  H^0(\Curve\times\Fqn,\O(D)/\O(D-D'_n))=
\prod_{w\in S_n}\O(D)_w/\O(D-D'_n)_w.
\end{equation*}
Again $C'_n$ has positive density which is bounded away from 0 for all
$n$.  Moreover, $f\in L(\Fqn)$ is in $L(\Fqn)^{\text{good}}$ if
and only if its image under
the natural homomorphism
\begin{equation*}
L(\Fqn)=H^0(\Curve\times\Fqn,\O(D))\to 
H^0(\Curve\times\Fqn,\O(D)/\O(D-D'_n))
\end{equation*}
lies in $C'_n$.

By the Riemann-Roch theorem, this homomorphism is surjective because
\begin{equation*}
\deg D>2g_\Curve-2+B\ge2g_\Curve-2+\deg D'_n.  
\end{equation*}
Also, the fibers of
this homomorphism all have the same cardinality, so the density 
of $L(\Fqn)^{\text{good}}$ in $L(\Fqn)$ is bounded away from 0 for all
$n$: there is an explicit constant $C'>0$ such that
\begin{equation*}
\frac{\#L(\Fqn)^{\text{good}}}{\#L(\Fqn)}>C'
\end{equation*}
for all $n$.

On the other hand, $X$ contains the complement of a
hypersurface in $L$ and so there is a constant $C''$ such that we have
a Lang-Weil type estimate
\begin{equation*}
\frac{\#(L(\Fqn)\setminus X(\Fqn))}{\#L(\Fqn)}<\frac{C''}{q^{n/2}}.
\end{equation*}
Thus
\begin{align*}
\frac{\#Y_n}{\#X(\Fqn)}&\ge
\frac{\#Y_n}{\#L(\Fqn)}\\
&=\frac{\#\left(L(\Fqn)^{\text{good}}\cap X(\Fqn)\right)}{\#L(\Fqn)}\\
&\ge\frac{\#L(\Fqn)^{\text{good}}
           -\#\left(L(\Fqn)\setminus X(\Fqn)\right)}{\#L(\Fqn)}\\
&>C'-\frac{C''}{q^{n/2}}
\end{align*}
and this proves the proposition.
\end{proof}

\subsection{Typical $D$ and local conditions}\label{ss:typical-local}
Now we discuss the local conditions to be used in the proof of the
main theorem.  In this subsection we give the ``typical'' conditions,
then in the next subsection we explain how they should be
modified in certain special circumstances.  The point is that the proofs of
\ref{lemma:global-signs} and \ref{prop:induced-zeroes} give conditions
under which certain zeroes can be avoided and we need to
insure that these conditions are satisfied.

Here are the typical conditions on $D$.  We require that the
effective divisor $D=\sum_va_v[v]$ satisfies:
\begin{enumerate}
  \item[(a)] $a_v$ is relatively prime to $d$ for all $v\in S_r$.
  \item[(b)] $a_v=0$ for all $v\in S_s\cup S_i\cup(|\n|\setminus S_r)$.
  \item[(c)] $a_v=1$ for at least one $v\not\in S_s\cup S_i\cup S_r$.
  \item[(d)] $\deg(D)>2g_\Curve-2+\deg(S_s\cup S_i)$.
  \item[(e)] $\deg(D)>\max(12g_\Curve+9,
        6\deg(\n)+11,72\deg(\rho)-(2g_\Curve-2))$. 
\end{enumerate}
The reason for the requirements on $\deg D$ will become clear later in
the proof.  Less stringent requirements are needed in many cases, but
we have chosen to simplify by making a uniform hypothesis.  

Our typical local conditions are as follows: $S_n$ will be the set of
places of $F_n$ over $S_s\cup S_i$.  If $w\in S_n$ lies over $S_s$,
then $C_{n,w}=\{1\}\subset F_{n,w}^\times/F_{n,w}^{\times d}$.  If $w\in
S_n$ lies over $S_i$, then $C_{n,w}\subset F_{n,w}^\times/F_{n,w}^{\times
  d}$ is the set of generators of the cyclic subgroup
  $\O_{n,w}^\times/\O_{n,w}^{\times d}$.

It is clear that the local conditions $(S_n,C_{n,w})$ are compatible with
$D$ and that an element $f\in F_n^\times$ which satisfies the local
conditions imposed by $(S_n,C_{n,w})$ satisfies the
local conditions in the sense of \ref{sss:local-conditions}.

\subsection{Three special situations}\label{ss:special-local}
First suppose we are considering part (1) of the main theorem,
$o=\{d/2\}$, $\rho$ is symplectically self-dual (of weight $w=-1$),
and for some $n$ the condition \ref{eqn:cond2} fails and $\alpha_n=\pm
1$.  Then we need to impose additional local conditions $C_{n,w}$ at
places $w$ over $|\n|\setminus(S_s\cup S_i)$.  So we replace condition
(d) above with
\begin{enumerate}
  \item[(d)] $\deg(D)>2g_\Curve-2+\deg(S_s\cup
        S_i)+\deg(|\n|\setminus(S_s\cup S_i))$.
\end{enumerate}
Note that the right hand side of this inequality is independent of $n$
and so we may fix one $D$ which works for all $n$.  The proof of
\ref{lemma:global-signs} shows that by imposing local conditions at
places over $|\n|\setminus(S_s\cup S_i)$, we may fix the sign in the
functional equation of $L(\rho\tensor\psi_f,F_n,T)$ where $\psi_f$ is
the quadratic character of $G_n$ corresponding to the extension
$F_n(f^{1/2})/F_n$.  We choose such local conditions so that the sign
is $+1$ when the degree of the $L$-function is even and so that the
sign is equal to that of $\alpha_n$ when the degree of the
$L$-function is odd.  (Cf.~\ref{lemma:quad-zeroes}.)  It is clearly
possible to do this in such a way that the new local conditions
$(S_n,C_{n,w})$ are still compatible with $D$.

The second special situation is when we consider part (2) of the main
theorem, $d$ is even (so that one of the orbits considered is
$o=\{d/2\}$), $\rho$ is symplectically self-dual (of weight $w=-1$),
and for some $n$ the
condition \ref{eqn:cond2} fails and $\alpha_n=\pm 1$.
Again we replace condition (d)
above with
\begin{enumerate}
  \item[(d)] $\deg(D)>2g_\Curve-2+\deg(S_s\cup
        S_i)+\deg(|\n|\setminus(S_s\cup S_i))$
\end{enumerate}
and we add local conditions at places $w$ over
$|\n|\setminus(S_s\cup S_i)$ to force the sign in the functional
equation of the quadratic twist $L(\rho\tensor\psi_{f},F_n,T)$
to take a certain value depending on $N$ and $\alpha_n$.  It is
clearly possible to do this in such a way that the new local
conditions $(S_n,C_{n,w})$ are still compatible with $D$.

The third special situation is when we consider part (1) of the main
theorem, $o=-o$, $a_o=\#o>1$, and either (a) $\rho$ is symplectically
self-dual (of weight $w=-1$) and for some $n$ $\alpha_n^{a_o}=1$; or
(b) $\rho$ is orthogonally self-dual (of weight $w=-1$) and for some
$n$ $\alpha_n^{a_o}=-1$.  By
Propositions~\ref{prop:avoiding-induced-zeroes} and
\ref{prop:induced-zeroes} we have forced zeroes if $\deg\rho$ is even
and hypothesis \ref{eqn:cond2} holds.  If one of these conditions
fails, the proof of Proposition~\ref{prop:avoiding-induced-zeroes}
tells us how to avoid forced zeroes and we must build this into the
definition of $D$.  More precisely, if $\rho$ is odd-dimensional we
choose $D$ so that the sum of the degrees of places in $|D|$ and not
over $S_s\cup S_i\cup S_r\cup|\n|$ is either odd or even, as required
to make the degree of $\cond(\rho\tensor\chi^i_f)$ even.  On the other
hand, if $\rho$ is even dimensional and hypothesis \ref{eqn:cond2}
fails, then we choose the coefficients $a_v$ of $D$ at places
$v\in|\n|\cap S_r$ so as to make the conclusion of \ref{eqn:cond2}
false.  Note that fixing the integer $a_v$ (modulo $d$) is the same as
fixing the local character of inertia $\chi^i_f$.  Note also that the
conditions on $D$ are independent of $n$ so there is one $D$ which
works for all $n$.  In this third special situation, the new
conditions are all on the coefficients of $D$ away from the places in
$S_n$, so the new $D$ and the local conditions $(S_n,C_{n,w})$ are
clearly compatible.

\subsection{Summary}
For the rest of the paper, we fix a divisor $D$ and compatible local
conditions $(S_n,C_{n,w})$ according to the recipe in
Subsections~\ref{ss:typical-local}-\ref{ss:special-local}.  (This data
depends of course on the data $F$, $\rho$, $d$, $S_s$, $S_i$, $S_r$, and
$(\alpha_n) $ fixed in Subsection~\ref{ss:data} and, when we are
considering part (1) of the theorem, on a fixed orbit
$o\subset\Z/d\Z$ for multiplication by $q$.)

The divisor $D$ determines a parameter space $X$ of functions, and the
conditions $(S_n,C_{n,w})$ determine a subset $Y_n\subset X(\Fqn)$ of
functions which satisfy the local conditions imposed by
$(S_n,C_{n,w})$.  Because we assumed the degree of $D$ is large
(specifically, because of the first hypothesis on $\deg(D)$ in
\ref{ss:typical-local} above), the density of $Y_n$ in $X(\Fqn)$ is
positive and bounded away from 0 for all sufficiently large $n$.
Moreover, by our choice of local conditions $(S_n,C_{n,w})$
the functions $f\in Y_n$ satisfy the local conditions in the sense of
Definition~\ref{defn:local-conds}.

\section{Twisted $L$-functions and sheaves on $X$}\label{s:cohomology}

In this section, we relate the twisted $L$-functions
$L(\rho\tensor\sigma_{o,f},F_n,T)$ to certain sheaves on the parameter
space $X$.  We assume familiarity with the basic formalism and
techniques of \'etale sheaves and their cohomology, as explained for
example in \cite[{[Arcata] and [Rapport]}]{SGA4.5} or \cite{Milne},
and in much more detail in \cite{SGA4} and \cite{SGA5}.

\subsection{$L$-functions and cohomology}
We begin in this subsection by reviewing Grothendieck's
cohomological expression for the $L$-functions
$L(\rho\tensor\sigma_{o,f},F_n,T)$, ``one $f$ at a time."

\subsubsection{}
Let $\tau:G\to\GL_r(E)$ be a continuous Galois representation such
that there exists a non-empty Zariski open subset $j:U\into\Curve$
with $\tau$ unramified at all places in $U$, i.e., such that
$\tau$ factors through $\pi_1(U,\overline{\eta})$.  Then there is a
twisted constant constructible (i.e., lisse) sheaf of $E$ vector
spaces $\FF_U$ on $U$ corresponding to $\tau$.  (Briefly, since $G$
and $\pi_1(U,\overline\eta)$ are compact we may conjugate $\tau$ so
that its image lies in $\GL_r(\O_E)$.  If $\m$ denotes the maximal
ideal of $\O_E$, reducing modulo powers of $\m$ gives representations
$\pi_1(U,\overline\eta)\to\GL_r(\O_E/\m^n)$ into finite groups.  These
correspond to \'etale sheaves of $\O_E/\m^n$-modules, free of rank
$r$.  For varying $n$, these finite sheaves collate into a $\m$-adic
system and tensoring with $E$ gives $\FF_U$.  Here of course we are
using the standard abuse of terminology, according to which a ``lisse
sheaf of $E$ vector spaces'' is actually an inverse system of twisted
constant, constructible sheaves of $\O_E/\m^n $-modules, up to
torsion.)

Conversely, given a lisse sheaf of $E$ vector spaces on some non-empty
open subset $U$ of $\Curve$, taking the stalk at $\overline\eta$
yields a continuous representation of $G$.  These constructions set up
an equivalence of categories between lisse sheaves of $E$ vector
spaces on $U$ and continuous representations of
$\pi_1(U,\overline\eta)$ on finite dimensional $E$ vector spaces.  (We
refer to \cite[{[Rapport] \S2}]{SGA4.5} or \cite[I.5, II.1, and
V.1]{Milne} for more details, and \cite[VII, VIII, IX]{SGA4} plus
\cite[V, VI]{SGA5} for many more details.)

Given $\tau$ as above, form $\FF_U$ and set $\FF_\tau=j_*\FF_U$.  Note
that $j^*\FF_\tau=\FF_U$.  If $j':V\into U$ is a smaller Zariski open
set, then it follows easily from the definitions that
$j'_*\FF_V\cong\FF_U$ and so $\FF_\tau$ is independent of the choice
of $U$.

\subsubsection{}
A ``middle extension'' sheaf of $E$ vector spaces on $\Curve$ is a
constructible sheaf $\FF$ of $E$ vector spaces for the \'etale topology
such that: (i) there exists a non-empty Zariski open $j:U\into\Curve$
such that $j^*\FF$ is lisse; and (ii) for one (and thus any)  such $U$,
$j_*j^*\FF\cong\FF$.  The preceding subsection describes a functor
from the category of finite dimensional continuous representations of $G$ on
vector spaces over $E$ ramified only at a finite set of places to the
category of middle extension sheaves of $E$ vector spaces on $\Curve$.
This functor is an equivalence of categories whose quasi-inverse
sends a sheaf $\FF$ to its geometric generic stalk
$\FF_{\overline\eta}$ equipped with the natural action of $G$.

\subsubsection{}\label{sss:middle-ext-pullback}
Suppose $\Curve'\to \Curve$ is an \'etale Galois cover and $\FF$ is
the middle extension sheaf on $\Curve'$ corresponding to a
representation $\tau$ of the fundamental group of $\Curve'$.  If
$g\in\pi_1(\Curve,\overline\eta)$, then $g$ induces an automorphism
$g:\Curve'\to \Curve'$.  We have that $g^*\FF$ is the middle extension
sheaf corresponding to the representation $\tau^g$, defined by
$\tau^g(h)=\tau(ghg^{-1})$.  We will apply this remark below in the
case where $\Curve'=\Curve\times\Fqn$ and $g$ is a lift of the
geometric Frobenius under
$\pi_1(\Curve,\overline\eta)\onto\gal(\Fqbar/\Fq)$.

\subsubsection{}
Some caution is required when applying standard constructions of
linear algebra (such as $\tensor$ and $\Hom$) in the category of
middle extension sheaves.  For example, it is not true in general that 
$\FF_{\tau_1\tensor\tau_2}\cong\FF_{\tau_1}\tensor\FF_{\tau_2}$.  What
is true is that if 
$j:U\into\Curve$ is a Zariski open such that $\tau_1$ and
$\tau_2$ factor
through $\pi_1(U,\overline{\eta})$, then 
$\FF_{\tau_1\tensor\tau_2}\cong
j_*\left(j^*(\FF_{\tau_1})\tensor j^*(\FF_{\tau_2})\right)$. 
In what follows we will be explicit about constructions like this
one.

\subsubsection{}
If $\FF$ is a constructible $\ell$-adic sheaf on a scheme $X$ we write
$H^i(X,\FF)$ and $H_c^i(X,\FF)$ for the cohomology and cohomology with
compact supports of $X$ with coefficients in $\FF$.  There is a
natural ``forget supports'' morphism $H^i_c(X,\FF)\to H^i(X,\FF)$
which in general is neither injective nor surjective.  

If $\FF_\tau$ is the middle extension sheaf on $\Curve$ attached to a
representation $\tau$ of $G$ unramified over $U$ and
$\Ubar=U\times\Fqbar$, then one has that
$H^0(\Cbar,\FF_\tau)=H^0(\Ubar,\FF_U)$ is the $G_{\infty}=\gal(\overline
F/\Fqbar F)$-invariants in the representation space of $\tau$.  If $U$
is affine (i.e., a proper subset of $X$), then $H^2(\Ubar,\FF_U)=0$.
By Poincar\'e duality, $H^0_c(\Ubar,\FF_U)=0$ and $H^2_c(\Ubar,\FF_U)$
is the $G_{\infty}$-coinvariants of $\tau$, with $\gal(\Fqbar/\Fq)$ action
twisted by $E(-1)$.

The following lemma is well-known but I know of no convenient
reference for the proof.

\begin{lemma}\label{lemma:IH}
Let $\tau$ be a representation of $G$ as above which is unramified
over the open $j:U\into \Curve$.  Form the sheaves $\FF_U$ and
$\FF_\tau=j_*\FF_U$.  Then $j^*:H^1(\Cbar,\FF_\tau)\to
H^1(\Ubar,\FF_U)$ is injective and
\begin{equation*}
j^*\left(H^1(\Cbar,\FF_\tau)\right)=
\im\left(H^1_c(\Ubar,\FF_U)\to H^1(\Ubar,\FF_U)\right).
\end{equation*}
\end{lemma}

\begin{proof}
Consider the Leray spectral sequences for $j$ with and without compact
supports.  The exact sequences of low degree terms and the ``forget
supports'' maps yield a commutative diagram with exact rows
\begin{equation*}
\xymatrix{0\ar[r]&H^1(\Cbar,\FF_\tau)\ar[r]^{a=j^*}
&H^1(\Ubar,\FF_U)\ar[r]&H^0(\Cbar,R^1\FF_\tau)\\
0\ar[r]&H^1_c(\Cbar,R^0_c\FF_\tau)\ar[r]^{d=j^*}\ar[u]^{b}
&H^1_c(\Ubar,\FF_U)\ar[r]\ar[u]_{c}&H^0_c(\Cbar,R^1_cj_*\FF_U).}
\end{equation*}
In particular,  $j^*:H^1(\Cbar,\FF_\tau)\to
H^1(\Ubar,\FF_U)$ is injective.  Since $j$ is quasi-finite and separated,
$R^1_cj_*\FF_U=0$ and so $d$ is an isomorphism.  On the other hand, we have
an exact sequence of constructible sheaves on $\Cbar$
\begin{equation*}
0\to j_!\FF_U\to \FF_\tau \to i_*i^*\FF_\tau\to0
\end{equation*}
where $i:Z\into \Curve$ is the complement of $U$.  Since $j$ is \'etale,
$R^0_cj_*\FF_U=j_!\FF_U$ and so taking cohomology with compact supports
yields an exact sequence
\begin{equation*}
H^1_c(\Cbar,j_!\FF)\to H^1_c(\Cbar,\FF_\tau)\to H^1_c(\Cbar,i_*i^*\FF_\tau).
\end{equation*}
But $i_*i^*\FF_\tau$ is a skyscraper sheaf and so
$H^1_c(\Cbar,i_*i^*\FF_\tau)=0$.  This shows that $b$ is surjective.  Thus
we have $\im(a)=\im(ab)=\im(cd)=\im(c)$, as desired.
\end{proof}

\subsubsection{}
It follows from the lemma and Poincar\'e duality that if $\tau$ is
self-dual of some weight $w$ then $H^1(\Cbar,\FF_\tau)$
is self-dual of weight $w+1$, i.e., we have a perfect pairing of
representations of $\gal(\Fqbar/\Fq)$
\begin{equation*}
H^1(\Cbar,\FF_\tau)\times H^1(\Cbar,\FF_\tau)\to E(-w-1).
\end{equation*}
If $\tau$ is orthogonally (resp.~symplectically) self-dual, then
$H^1(\Cbar,\FF_\tau)$ is symplectically (resp.~orthogonally)
self-dual.

\subsubsection{}
Let $\FF=\FF_\rho$ be the middle extension sheaf on $\Curve$
corresponding to the representation $\rho$ fixed in
Subsection~\ref{ss:data}.  The Grothendieck-Lefschetz trace formula
computes the $L$ function of the representation $\rho$ in terms of the
cohomology of the sheaf $\FF$.  More precisely,
\begin{equation*}
L(\rho,F,T)=\prod_{i=0}^2
\det\left(1-Fr\,T|H^i(\Cbar,\FF_\rho)\right)^{(-1)^{i+1}}
\end{equation*}
where as usual $Fr$ is the geometric ($q^{-1}$-power on $\Fqbar$)
Frobenius endomorphism of $\Cbar$.  The cohomology groups
are finite dimensional $E$ vector spaces and so the $L$-function is a
rational function in $T$.  When $\rho$ is irreducible and
geometrically non-trivial
(or more generally a direct sum of geometrically non-trivial
irreducibles), the groups $H^i(\Cbar,\FF_\rho)$ vanish
for $i=0,2$ and the $L$-function is a polynomial in $T$.

\subsubsection{}\label{sss:coh-and-Ls}
Now assume that $o$ is an orbit of multiplication by $q^n$ on
$(\Z/d\Z)^\times$ and $f\in X(\Fqn)\subset F_n^\times$.  Then
$\rho\tensor\sigma_{o,f}$ is semisimple as a representation of $G_n$
and also semisimple when restricted to $G_{\infty}$.
By our choice of $D$ defining $X$, $\sigma_{o,f}$ is totally ramified at at
least one place where $\FF$ is lisse and so $\rho\tensor\sigma_{o,f}$
does not contain the trivial representation when restricted to $G_{\infty}$.
Thus we have
\begin{equation*}
L(\rho\tensor\sigma_{o,f},F_n,T)=
\det\left(1-Fr^{n}\,T|
H^1(\Cbar,\FF_{\rho\tensor\sigma_{o,f}})\right).
\end{equation*}

As we saw in Lemma~\ref{lemma:sigma-self-dual}, if $o=-o$, then
$\sigma_{o,f}$ is orthogonally self-dual and so if $\rho$ is
self-dual, then so is $\rho\tensor\sigma_{o,f}$, with the same sign as
$\rho$.  In particular, if $\rho$ is orthogonally self-dual (of weight
$w=-1$), then $H^1(\Cbar,\FF_{\rho\tensor\sigma_{o,f}})$
is symplectically self-dual and so the $L$-function has even degree
in $T$ and satisfies the functional equation
\begin{equation*}
L(\rho\tensor\sigma_{o,f},F_n,T)=L(\rho\tensor\sigma_{o,f},F_n,1/T)
\end{equation*}
in other words, the root number $W(\rho\tensor\sigma_{o,f},F_n)=1$.

\subsubsection{}  
Recall from Lemma~\ref{lemma:reps} that if $f\in F_n^\times$, $o$ is
an orbit of multiplication by $q^n$ on $(\Z/d\Z)^\times$, and
$a_o=\#o$, then when restricted to $G_{na_o}$,
$\sigma_{o,f}\cong\oplus_{i\in o}\chi_f^i$.  This implies that as
sheaves on $\Curve\times\F_{q^{na_o}}$
\begin{equation*}
  \FF_{\rho\tensor\sigma_{o,f}}\cong\oplus_{i\in o}
\FF_{\rho\tensor\chi_f^i}
\end{equation*}

Similarly, since $\sigma_{o,f}\cong\ind^{G_n}_{G_{na_o}}\chi_f^i$ for
any $i\in o$ (Lemma~\ref{lemma:induced}), we have
\begin{equation*}
  \FF_{\rho\tensor\sigma_{o,f}}\cong b_*\FF_{\rho\tensor\chi_f^i}
\end{equation*}
where $b:\Curve\times\F_{q^{na_o}}\to\Curve$ is the natural
projection.  

By the remark in \ref{sss:middle-ext-pullback} and Lemma~\ref{lemma:Frob-twist}, 
if $\Phi$ is an
element of $G$ lifting the geometric Frobenius of $\gal(\Fqbar/\Fq)$,
then 
\begin{equation*}
\Phi^*\FF_{\rho\tensor\chi^i_f}\cong\FF_{(\rho\tensor\chi^i_f)^\Phi}
\cong\FF_{\rho\tensor\chi^{iq}_f}.
\end{equation*}

\subsubsection{}\label{sss:induced-zeroes-proof}
We can now give the proof of Proposition~\ref{prop:induced-zeroes}.
In light of \ref{sss:coh-and-Ls} and our assumption that $\rho$ has
weight $w=-1$, what is to be proven is that all of the $a_o$-th roots
of $-\sgn(\rho)$ appear as eigenvalues of $Fr^n$ on
$H^1(\Cbar,\FF_{\rho\tensor\sigma_{o,f}})$.  (Here
$\sgn(\rho)$ is $1$ if $\rho$ is orthogonally self-dual and $-1$ if it is
symplectically self-dual.)  Since
\begin{equation*}
H^1(\Cbar,\FF_{\rho\tensor\sigma_{o,f}})\cong
\oplus_{i\in o}
H^1(\Cbar,\FF_{\rho\tensor\chi^i_{f}})
\end{equation*}
and $(\chi_f^i)^{Fr^n}=\chi_f^{iq^n}$, the matrix of $Fr^n$ is
a block permutation matrix, i.e., has the form 
\begin{equation}\label{eqn:block-permutation}
  \begin{pmatrix}
0&0&0&\cdots&A_{iq^{n(a_o-1)}}\\
A_i&0&0&\cdots&0\\
0&A_{iq^n}&0&\cdots&0\\
\vdots&\vdots&\vdots&\ddots&\vdots\\
0&0&0&\cdots&0
\end{pmatrix}
\end{equation}
where $A_j$ is the matrix of
$Fr^n:H^1(\Cbar,\FF_{\rho\tensor\chi^j_{f}})\to
H^1(\Cbar,\FF_{\rho\tensor\chi^{jq^n}_{f}})$.  This
implies that the eigenvalues of $Fr^n$ are all of the $a_o$-th
roots of the eigenvalues of $Fr^{na_o}$.  (I.e., if
$P_{Fr^n}(T)$ and $P_{Fr^{na_o}}(T)$ are the characteristic
polynomials of $Fr^n$ and $Fr^{na_o}$ on
$H^1(\Cbar,\FF_{\rho\tensor\sigma_{o,f}})$ and
$H^1(\Cbar,\FF_{\rho\tensor\chi^i_{f}})$, then
$P_{Fr^n}(T)=P_{Fr^{na_o}}(T^{a_o})$.)  Thus we must show
that $-\sgn(\rho)$ is an eigenvalue of $Fr^{na_o}$
on $H^1(\Cbar,\FF_{\rho\tensor\chi^i_{f}})$.

We assumed that $\rho$ is self-dual of weight $w=-1$.
Since $\sigma_{o,f}$ is
orthogonally self-dual,
$H^1(\Cbar,\FF_{\rho\tensor\sigma_{o,f}})$ is literally
self-dual (i.e., self-dual of weight 0),
of sign opposite to that of $\rho$.  Moreover, the subspaces
$H^1(\Cbar,\FF_{\rho\tensor\chi^i_{f}})$ and
$H^1(\Cbar,\FF_{\rho\tensor\chi^{-i}_{f}})$ are put in
duality by the restriction of the form.  

Let us fix bases of each
$H^1(\Cbar,\FF_{\rho\tensor\chi^i_{f}})$ such that for
all $i$ the chosen basis of
$H^1(\Cbar,\FF_{\rho\tensor\chi^i_{f}})$ is dual to that
of $H^1(\Cbar,\FF_{\rho\tensor\chi^{-i}_{f}})$.  Then in
the matrix \ref{eqn:block-permutation}, the self-duality implies that
$A_{iq^{j+a_o/2}}=(A_{iq^j})^\vee$ for $0\le j<a_o/2$ and
$A_i=-\sgn(\rho)(A_{iq^{a_o/2}})^\vee$ where $A^\vee$ denotes the
inverse transpose of $A$.

Thus, the matrix of $Fr^{na_o}$ on
$H^1(\Cbar,\FF_{\rho\tensor\chi^i_{f}})$ is
\begin{equation*}
-\sgn(\rho)(A_{iq^{a_o/2-1}}^\vee\cdots A_{i}^\vee)
(A_{iq^{a_o/2-1}}\cdots A_i)=-\sgn(\rho)B^\vee B
\end{equation*} 
where $B=A_{iq^{a_o/2-1}}\cdots A_i$.  The first part of the following
lemma then finishes the proof.

\begin{lemma}\label{lemma:BcheckB}
  Consider invertible $N\times N$ matrices $B$ over an infinite field
  and let $B^\vee$ denote the inverse transpose of $B$.  If $N$ is odd
  then for every $B$, $B^\vee B$ has 1 as an eigenvalue; moreover,
  given $\alpha\neq1$ in the ground field, there exists a $B$ such
  that the multiplicity of 1 as an eigenvalue of $B^\vee B$ is 1 and
  $\alpha$ is not an eigenvalue of $B^\vee B$.  If $N$ is even, for
  any $\alpha$ there exists a $B$ such that $\alpha$ is not an
  eigenvalue of $B^\vee B$.  All of the above remains true of we
  restrict to matrices $B$ having any fixed non-zero determinant.
\end{lemma}

\begin{proof}
  First, note that $(B^\vee B)^t=B^tB^{-1}$ and $(B^\vee
  B)^{-1}=B^{-1}B^t$.  This implies that $(B^\vee B)^{-1}$ is
  conjugate to $(B^\vee B)^t$, which, by the Jordan form, is conjugate
  to $B^\vee B$.  Thus the set of eigenvalues of $B^\vee B$ is
  invariant under $\lambda\mapsto\lambda^{-1}$.  On the other hand,
  the product of the eigenvalues of $B^\vee B$ is $\det(B^\vee B)=1$.
  If $N$ is odd, this implies that at least one of the eigenvalues
  must be 1.  

For the existence assertions, we may build up a suitable $B$ using
$2\times2$ blocks of the
  form
\begin{equation*}
\pmat{a&b\\0&1}.
\end{equation*}
Indeed, these matrices have
determinant $a$ and the eigenvalues of
\begin{equation*}
\pmat{a&b\\0&1}^\vee
\pmat{a&b\\0&1}=\pmat{1&b/a\\-b&1-b^2/a}
\end{equation*}
vary with $b$ and avoid 1 and $\alpha$ for suitable $b$.  
\end{proof}

\subsection{Globalization}
Our next task is to define for each orbit $o$ of multiplication by $q$
on $(\Z/d\Z)^\times$ a sheaf $\GG_o$ on $X$ whose stalk at a geometric
point over $f\in X(\Fqn)$ is the cohomology group
$H^1(\Cbar,\FF_{\rho\tensor\sigma_{o,f}})$.  We will use
several constructions and results from \cite[Chaps.~5 and 6]{Katz}.  There
are some errors in Chapter~5, which Katz has addressed.  We refer
to his web site ({\tt http://www.math.princeton.edu/\~{}nmk}) for a
corrected version.

\subsubsection{}
Consider the product $X\times\Curve$ with its two projections $\pi_1$
and $\pi_2$ to $X$ and $\Curve$ respectively.
On the product $X\times\Curve$ we have a ``universal rational
function'' $F_{univ}$, characterized by the formula
$F_{univ}(f,p)=f(p)$.  The divisor of poles of $F_{univ}$ is
$D\times\Curve$, its divisor of zeroes is finite \'etale over $X$ (via
$\pi_1$) of degree equal to $\deg D$ and the divisor of zeroes of
$F_{univ}$ is disjoint from its divisor of poles.  

Let $\DD\subset X\times\Curve$ be the reduced divisor whose support is
the union of the divisor of $F_{univ}$ and $X\times (\Curve\setminus
U)$ where $j:U\into \Curve$ is a Zariski open subset over which $\rho$
is unramified.  Also let $\tilde j:V=(X\times\Curve)\setminus\DD\into
X\times\Curve$ be the inclusion.

\subsubsection{}
Let $\lambda:\XX\to X\times\Curve$ be the normalization of
$X\times\Curve$ in the field extension
$\Fq(X\times\Curve)(F_{univ}^{1/d})$ of $\Fq(X\times\Curve)$.  Clearly
$\lambda$ has degree $d$ and is \'etale over $V\subset
X\times\Curve$.

\subsubsection{}
Let $\uE$ denote the constant sheaf on $\XX$ with stalk $E$ and
consider $\lambda_*\uE$ and its restriction $\tilde j^*\lambda_*\uE$
to $V$.  Since $\lambda$ is \'etale of degree $d$ over $V$, $\tilde
j^*\lambda_*\uE$ is lisse of rank $d$.  The argument of
Lemma~\ref{lemma:reps} applies in this situation and we have a
factorization
\begin{equation*}
\tilde j^*\lambda_*\uE\cong
\bigoplus_{o\subset\Z/d\Z}\Sigma_o
\end{equation*}
where the sum is over orbits of multiplication by $q$ on $\Z/d\Z$ and
$\Sigma_o$ is a lisse sheaf of $E$-vector spaces on $V$ of rank $\#o$.

\subsubsection{}
After a small base extension, $\Sigma_o$ becomes isomorphic to a sum
of rank 1 lisse sheaves.  More precisely, the base change of $\lambda$
to $\Fq(\mu_d)$, i.e., 
\begin{equation*}
\XX\times\Fq(\mu_d)\to
X\times\Curve\times\Fq(\mu_d), 
\end{equation*}
is Galois with Galois group naturally identified with $\mu_d(\overline
F)$ by $\sigma\mapsto\sigma(F_{univ}^{1/d})/F_{univ}^{1/d}$.
Composing with an isomorphism $\mu_d(\overline F)\to\mu_d(E)$ (the
same one we used in \ref{sss:base-change}), we get a character
$\chi_{F_{univ}}$ which is unramified over $V$.  We let $\LL_{univ}^i$
be the rank 1 lisse sheaf on $V\times\Fq(\mu_d)$ corresponding to
$\chi^i_{F_{univ}}$.  We note that $\LL_{univ}^i$ in fact descends to
$V\times\Fq(\mu_{d_o})$ where as before $d_o=d/\gcd(d,i)$ for any
$i\in o$.

With these notations, we have a factorization
\begin{equation*}
\Sigma_o|_{V\times\Fq(\mu_{d_o})}\cong\bigoplus_{i\in o}\LL_{univ}^i
\end{equation*}
of lisse sheaves of $E$-vector spaces on $V\times\Fq(\mu_{d_o})$.
(This is the global version of the factorization at the end of
Lemma~\ref{lemma:reps}.)  

Similarly, globalizing Lemma~\ref{lemma:induced}, we have
\begin{equation*}
  \Sigma_o=b_*(\LL_{univ}^i)
\end{equation*}
for any $i\in o$, where $b$ is the projection
$X\times\Fq(\mu_{d_o})\to X$.

Globalizing Lemma~\ref{lemma:Frob-twist} and
\ref{sss:middle-ext-pullback}, we have
$Fr^*(\LL^i_{univ})\cong\LL^{iq}_{univ}$.

Globalizing Lemma~\ref{lemma:sigma-self-dual}, we have that $\Sigma_o$
is self-dual if and only if $-o=o$, in which case it is orthogonally
self-dual of weight 0.

\subsubsection{}
Recall that $\pi_1$ and $\pi_2$ denote the projections from
$X\times\Curve$ to $X$ and $\Curve$ respectively. Let $\mu$ be the
restriction of $\pi_1$ to $V$.  We define
\begin{align*}
\GG_{o,*}&=R^1\mu_*
\left((\tilde j^*\pi_2^*\FF_\rho)\tensor\Sigma_o\right)\\
\intertext{and}
\GG_{o,!}&=R^1\mu_!
\left((\tilde j^*\pi_2^*\FF_\rho)\tensor\Sigma_o\right).
\end{align*}
Since $\deg(D)>2g+1$, the arguments of \cite[5.2.1 and 6.2.10]{Katz}
show that these are lisse sheaves on $X$ whose formation is compatible
with arbitrary change of base.

\subsubsection{}
There is a natural ``forget supports'' morphism
$\GG_{o,!}\to\GG_{o,*}$ and we define $\GG_o$ to be the image of this
morphism.  Again by \cite[5.2.1 and 6.2.10]{Katz}, $\GG_o$ is lisse of
formation compatible with arbitrary change of base and by Deligne
\cite[3.2.3]{WeilII}, it is $\iota$-pure of weight 0.

By standard base change results, the stalk of $\GG_o$ at a geometric
point $\overline{f}$ over $f\in X(\Fqn)$ is
\begin{equation*}
\GG_{o,\overline{f}}\cong\im\left(H^1_c(\Ubar,j^*\FF_\rho\tensor j^*\FF_{\sigma_{o,f}})
\to H^1(\Ubar,j^*\FF_\rho\tensor j^*\FF_{\sigma_{o,f}})\right)
\end{equation*}
where $j:U\into\Curve$ is a Zariski open over which $\rho$ is
unramified and $f$ is regular and non-zero. 
By Lemma~\ref{lemma:IH}, this is
\begin{equation*}
H^1(\Cbar,j_*(j^*\FF_\rho\tensor j^*\FF_{\sigma_{o,f}}))
=H^1(\Cbar,\FF_{\rho\tensor\sigma_{o,f}}).
\end{equation*}

\subsubsection{}
By \ref{sss:coh-and-Ls},
\begin{align*}
L(\rho\tensor\sigma_{o,f},F_n,T)
&=\det\left(\left.1-Fr^{n}T\right|
H^1(\Cbar,\FF_{\rho\tensor\sigma_{o,f}})\right)\\
&=\det\left(\left.1-Fr_{n,f}T\right|\GG_{o,\overline{f}}\right).
\end{align*}
Thus we may study the $L$-functions $L(\rho\tensor\sigma_{o,f},F_n,T)$
for every $f\in X(\Fqn)$ by studying the sheaf $\GG_o$.

\subsubsection{}\label{sss:Gi-def}
We have variants of $\GG_o$ over a small base extension of $X$.  More
precisely, if $i\in o$ and we work over $\Fq(\mu_{d_o})$ then we can
define
\begin{align*}
\GG_{i,*}&=R^1\mu_* \left((\tilde
j^*\pi_2^*\FF_\rho)\tensor\LL^i_{univ}\right),\\ 
\GG_{i,!}&=R^1\mu_!
\left((\tilde j^*\pi_2^*\FF_\rho)\tensor\LL^i_{univ}\right),\\
\intertext{and} 
\GG_i&=\im\left(\GG_{o,!}\to\GG_{o,*}\right).
\end{align*}
(Here we are abusing abusing notation slightly by using $\mu$, $\tilde
j$ and $\pi_2$ to denote various maps to and from
$X\times\Curve\times\spec\Fq(\mu_{d_o})$.) 
By arguments similar to those mentioned above, we have that $\GG_i$ is
lisse and its stalk at a geometric point over $f\in X(\Fqn)$ is
$H^1(\Curve\times
\Fqbar,\FF_{\rho\tensor\chi_f^i})$.

The Grothendieck-Ogg-Shafarevitch formula says that the rank of
$\GG_i$ is 
\begin{equation*}
(2g_\Curve-2)(\deg\rho)+\deg\cond(\rho\tensor\chi_f^i).
\end{equation*}
Because we assumed $\deg(D)$ is large (cf.~\ref{ss:typical-local}
(e)), \cite[5.3.6]{Katz} says that $\GG_i$ is irreducible and 
\cite[5.5.1 and 5.7.1]{Katz} say
that $\GG_i$ is self-dual on $\overline{X}$ if and only if $\FF$ is
self dual and $i=d/2$, in which case its sign is the opposite of that
of $\FF$.

\subsubsection{}
Over $X\times\Fq(\mu_d)$, $\GG_o$ factors.  
More precisely, we have
\begin{equation*}
  b^*\GG_o=\oplus_{i\in o}\GG_i
\end{equation*}
and
\begin{equation*}
  \GG_o=b_*\GG_i
\end{equation*}
for any $i\in o$, where $b:X\times\Fq(\mu_{d_o})\to X$ is the
projection.  Also, $Fr^*(\GG_i)\cong\GG_{iq}$.

\begin{prop}\label{prop:GGis-distinct}
 Write $\Xbar$ for $X\times\Fqbar$.  Then  we have
 \begin{enumerate}
   \item $\GG_i\cong\GG_j$ on $\Xbar$ if and only if $i\equiv j\pmod d$
   \item $\GG_i\cong\GG_j^\vee$ on $\Xbar$ if and only if $\FF$ is
   self-dual (of weight $w=-1$) and $i\equiv -j\pmod d$.
 \end{enumerate}
 More generally, if $f:Y\to\Xbar$ is a connected, finite, \'etale
 cover, then $f^*\GG_i$ and $f^*\GG_j$ are isomorphic (resp.~dual) if
 and only if $i\equiv j\pmod d$ (resp.~$\FF$ is self-dual (of
 weight $w=-1$) and $i\equiv -j\pmod d$).
\end{prop}

\begin{cor}\label{cor:irred/self-dual}
  The lisse sheaf $\GG_o$ on $X$ is irreducible.  It is self-dual if
  and only if $o=-o$ and $\FF$ is self-dual (of weight $w=-1$) on
  $\Xbar$ (and thus by our assumptions self-dual on $X$).  In this
  case $\GG_o$ is orthogonally self-dual if $\FF$ is symplectically
  self-dual and $\GG_o$ is symplectically self-dual if $\FF$ is
  orthogonally self-dual.
\end{cor}

\begin{proof}[Proof of Corollary~\ref{cor:irred/self-dual}]
  We have $\GG_o=b_*\GG_i$ where $b:X\times\Fq(\mu_{d_o})\to X$ is the
  natural projection.  But $\GG_i$ is irreducible
  (see~\ref{sss:Gi-def}) and $Fr^{j*}(\GG_i)\not\cong\GG_i$ unless
  $iq^j\cong i$ (by Proposition~\ref{prop:GGis-distinct}), so it
  follows from Mackey's criterion that $\GG_o$ is irreducible.

It is also clear that if $\FF$ is self-dual (of weight $w=-1$)
and $-o=o$, then $\GG_o$ is self-dual on $X$ with
the asserted sign.

Suppose then that $\GG_o$ is self-dual on $X$.  On
$\Xbar$ we have $\GG_o\cong\oplus_{i\in o}\GG_i$ and
$\GG_o^\vee\cong\oplus_{j\in o}\GG_j^\vee$.  Since each $\GG_i$ is
irreducible we must have $\GG_i\cong\GG_j^\vee$ for some $j\in o$.
Then by Proposition~\ref{prop:GGis-distinct}, $\FF$ is self-dual and
$j=-i$.  This holds for every $i\in o$, so $-o=o$.
\end{proof}

\begin{proof}[Proof of Proposition~\ref{prop:GGis-distinct}]
  The ``if'' parts of both statements are trivial.  The proofs of
  the converses rely heavily on the details of the proofs in
  \cite{Katz}, especially those in Chapter~5, not just the results
  themselves.  
  
  We work throughout on $\Cbar=\Curve\times\spec\Fqbar$ and
  $\Xbar=X\times\spec\Fqbar$.  Since we have assumed that the degree
  of $D$ is large (cf.~hypothesis (e) on $\deg(D)$ in
  \ref{ss:typical-local}), by \cite[5.4.8]{Katz}, we may write
  $D=D_1+D_2$ where the $D_i$ satisfy several conditions.  If $p>2$, the
  conditions are:
\begin{itemize}
\item $\deg(D_1)\ge2g+2$ 
\item $\deg(D_2)\ge2g+1$
\item the coefficients of $D_2$ are invertible modulo $p$
\item If $D=\sum a_iP_i$ (where the $P_i$ are distinct $\Fqbar$ points
 of $\Curve$) and $a_i>2$, then $P_i\in|D_2|$
\item if $4|d$ then
  $2\deg(D_1)<2g_\Curve-2+\deg(D)-2\deg(|\n|\setminus|D_2|)$
\item if $4|d$ and $g_\Curve=0$ then in addition $\deg(D_2)\ge2$.
\end{itemize}
If $p=2$, the conditions are:
\begin{itemize}
\item $\deg(D_i)\ge 6g_\Curve+3$
\item the coefficients of $D_2$ are odd
\item If $D=\sum a_iP_i$ (where the $P_i$ are distinct $\Fqbar$ points
 of $\Curve$) and $a_i>2$, then $P_i\in|D_2|$.
\end{itemize}
  
We write $L(D_i)$ for $H^0(\Cbar,\O_{\Cbar}(D_i))$.  Fix a function
$f_1\in L(D_1)$ which has distinct zeroes, all of which are disjoint
from $|\n|\cup|D|$.  (It is elementary that the set of such functions
$f_1$ is dense in $L(D_1)$; cf.~\cite[5.0.6]{Katz}.)  Consider
functions $f_2\in L(D_2)$ which satisfy the following conditions: (i)
$f_2$ has distinct zeroes, all of which are disjoint from
$|\n|\cup|D|\cup f_1^{-1}(0)$; (ii) the ramification of $f_2$ is
minimal in the following strong sense: if $p>2$ then all of the zeroes
of the differential $df$ are simple zeroes and if $p=2$, then all of
the zeroes of $df$ have multiplicity exactly 2; and (iii) $f_2$
separates the points in
\begin{equation*}
S=\left(\{\hbox{zeroes of $df_2$}\}\cup f_1^{-1}(0)\cup|D|
\cup|\n|\right)\setminus|D_2|
\end{equation*}
i.e., each $s\in S$ is the only element of $S$ in its fiber
$f_2^{-1}(f_2(s))$.  Theorems 2.2.6 and 2.4.2 of \cite{Katz} guarantee
that the set of $f_2$ satisfying these restrictions is a dense open
subset of $L(D_2)$.

The map $F$ defined by $F(t)=f_1(t-f_2)$ defines a morphism from the
open subset $U=\A^1\setminus S$ of the affine line over $\Fqbar$ (with
coordinate $t$) to $\Xbar$.  Proposition 5.3.7 of \cite{Katz} says
that we can almost recover $\FF$ from $\GG_i$ via $F$ and
\cite[Thm.~5.4.9]{Katz} gives a reasonably complete description of the
ramification of $F^*(\GG_i)$ on $\P^1\setminus U$.  More precisely, we
have an isomorphism of perverse sheaves on $U$
\begin{equation*}
F^*(\GG_i)[1]\cong
\left(f_{2*}j_{2*}j_1^*(\FF\tensor\LL_{\chi^i(f_1)})\right)[1]
*_{mid,+}j_*\LL_{\chi^i}[1]
\end{equation*}
where $j_i:\Cbar\setminus|D|\into\Cbar\setminus|D_i|$ and
$j:\G_m=\A^1\setminus\{0\}\into \A^1$ are the natural
inclusions and $*_{mid,+}$ is the middle additive convolution
(for which we refer to \cite[Chapter~4]{Katz}).  If $\GG_i\cong\GG_j$
and $i\not\equiv j\pmod d$ we deduce
an isomorphism
\begin{equation*}
\left(f_{2*}j_{2*}j_1^*(\FF\tensor\LL_{\chi^i(f_1)})\right)[1]\cong
\left(f_{2*}j_{2*}j_1^*(\FF\tensor\LL_{\chi^j(f_1)})\right)[1]
*_{mid,+}j_*\LL_{\chi^{j-i}}[1].
\end{equation*}

Now if $p>2$ there is a point $t$ of $\A^1$ so that $f_2:\Cbar\to\A^1$
is ramified, with ramification index $e=2$ at exactly one point over
$t$ and is unramified at the others and so that
$\FF\tensor\LL_{\chi^i(f_1)}$ and $\FF\tensor\LL_{\chi^j(f_1)}$ are
unramified at all points over $t$.  Let $\HH_i\cong
f_{2*}j_{2*}j_1^*(\FF\tensor\LL_{\chi^i(f_1)})$, viewed as a
representation of $I(t)$, the inertia group at $t$, and similarly for
$\HH_j$.  Then, using a superscript to denote invariants,
$\HH_i/\HH_i^{I(t)}\cong\HH_j/\HH_j^{I(t)}$ and these representations
are spaces of dimension $\rk\FF$ on which $I(t)$ acts by a non-trivial
character of order 2.  But $\HH_j$ is in the class ${\mathcal P}_{conv}$
(see \cite[4.0]{Katz}) and so by \cite[4.1.10(1a)]{Katz} we have
\begin{equation*}
  \HH_i/\HH_i^{I(t)}\cong\HH_j/\HH_j^{I(t)}\tensor\LL_{\chi^{j-i}(x-t)}
\end{equation*}
as $I(t)$ representations.  This obviously contradicts the assumption
$i\not\equiv j\pmod d$ which
concludes the proof of part (1) when $p>2$.

The argument for $p=2$ is similar, but we have to contend with wild
ramification.  In this case, there is a point $t$ of $\A^1$ so that
$f_2:\Cbar\to\A^1$ is ramified, with ramification index $e=2$ and
$df_2$ vanishing to order exactly 2 at exactly one point over $t$ and
is unramified at the others and so that $\FF\tensor\LL_{\chi^i(f_1)}$
and $\FF\tensor\LL_{\chi^j(f_1)}$ are unramified at all points over
$t$.  Then, with $\HH_i$ and $\HH_j$ defined as before, we have that
$\HH_i/\HH_i^{I(t)}$ has dimension $\rk\FF$ and $I(t)$ acts through a
character of Swan conductor 1 (see \cite[2.7.1]{Katz}); moreover, the
character only depends on $f_2$, not on $i$.  Applying \cite[4.1.10
and 4.2.1]{Katz}, we have that
$F^*(\GG_i)/F^*(\GG_i)^{I(t)}\cong\chi^{2i}\rho$ as $I(t)$
representations, where $\rho$ is a character of 2-power order and Swan
conductor 1.  Thus if $\GG_i\cong\GG_j$ we have
$\chi^{2i}\rho\cong\chi^{2j}\rho'$ where $\rho$ and $\rho'$ have
2-power order.  Since $d$ is prime to $p=2$, we conclude that $i\equiv
j\pmod d$.

We now turn to the proof of part (2) of the proposition.  Let us
temporarily denote the sheaf $\GG_i$ constructed from $\FF$ as
$\GG(\FF,i)$, so that our hypothesis is that
$\GG(\FF,i)\cong\GG(\FF,j)^\vee$.  Since
$\GG(\FF,j)^\vee\cong\GG(\FF^\vee,-j)$, our hypothesis is equivalent
to $\GG(\FF,i)\cong\GG(\FF^\vee,-j)$.  The argument proving the first
part of the proposition does not use much about $\FF$; more precisely,
the only information about $\FF$ we use is the support of its Artin
conductor.  Since $\FF$ and $\FF^\vee$ have the same Artin conductor,
the argument generalizes immediately to prove that $i\cong-j\pmod d$.
Thus it remains to show that $\GG(\FF,i)\cong\GG(\FF^\vee,i)$ implies that
$\FF\cong\FF^\vee$. To that end, we choose functions $f_1$ and $f_2$
satisfying the same hypotheses as before.  Let
\begin{equation*}
  \HH=j_{2*}j_1^*(\FF\tensor\LL_{\chi^i(f_1)})
\end{equation*}
and
\begin{equation*}
  \HH'=j_{2*}j_1^*(\FF^\vee\tensor\LL_{\chi^i(f_1)}).
\end{equation*}
As representations $\gal(\overline F/F)$, $\HH$ and $\HH'$ are
irreducible and by assumption we have
\begin{equation*}
f_{2*}\HH[1]*_{mid,+}j_*\LL_{\chi^i}[1]\cong
f_{2*}\HH'[1]*_{mid,+}j_*\LL_{\chi^i}[1]
\end{equation*}
which implies $f_{2*}\HH=f_{2*}\HH'$.  

Choose a point $t\in\A^1$ such that $f_2:\Curve\to\A^1$ is unramified 
at every point over $t$, $\FF$ is unramified at every point over $t$,
and exactly one point of $f_1^{-1}(0)$, call it $s_0$, lies over $t$.
Then $\HH$ is ramified at $s_0$ and unramified at the other points
over $t$; more precisely, as a representation of $I(s_0)$, the inertia
group at $s_0$, $\HH$ is isomorphic to a direct sum of $\rk\FF$ copies
of $\LL_{\chi^i(f_1)}$, on which $I(s_0)$ acts by a non-trivial
character of finite order.  The same is true of $\HH'$.  Now we have
inclusions
\begin{equation*}
  \HH\into f_2^*f_{2*}\HH\cong f_2^*f_{2*}\HH'\hookleftarrow\HH'
\end{equation*}
Since the sheaves $\HH$ and $\HH'$ are irreducible, their images in
the middle either coincide or are linearly independent.  But we can
see that the latter is impossible by noting that as a representation
of $I(s_o)$, $f_2^*f_{2*}\HH$ has an unramified subspace of
codimension $\rk\FF$.  Since, as representations of $I(s_0)$, $\HH$
and $\HH'$ are both totally ramified of dimension $\rk\FF$, there is
not enough room in $f_2^*f_{2*}\HH$ for them to be linearly independent.
Thus we have an isomorphism $\HH\cong\HH'$.  Since $j_1$ and $j_2$ are
open immersions and $\FF$ is a middle extension, it follows
immediately that $\FF\cong\FF^\vee$.  This completes the proof of part
(2) of the proposition.

For the ``more generally,'' suppose $f:Y\to \Xbar$ is a connected,
finite, \'etale cover such that $f^*\GG_i$ is isomorphic to
$f^*\GG_j$.  Choose functions $f_1$ and $f_2$ and define $S$ as above, let
$F:U=\A^1\setminus S\to\Xbar$ be defined by $F(t)=f_2(t-f_1)$, and let
$g:V\to U$ be the pull back of $f:Y\to \Xbar$ to $U$.  We may choose
the $f_i$ so that $V$ is connected.  The proofs of \cite[1.5.1 and
1.7.1]{Katz} (applied in the context of \cite[5.4.9 or 5.6.2]{Katz})
show that each $F^*\GG_i$ is ``Lie irreducible'' i.e., it remains
irreducible when restricted to any connected, finite, \'etale cover of
$U$.  Considering the action of $\pi_1(U)$ on
$\Hom_V(g^*F^*\GG_i,g^*F^*\GG_j)$, which is 1-dimensional by Schur's
lemma, we see that there exists a rank 1 lisse sheaf $\LL_\psi$ (with
associated character $\psi$ of $\pi_1(U)$) such that $F^*\GG_i\cong
F^*\GG_j\tensor\LL_\psi$.  

We are going to use the nature of the ramification of $\GG_i$ and
$\GG_j$ to show that such a $\psi$ must be trivial.  First of all,
$\psi$ is unramified on $U=\A^1\setminus S$.  Since we assumed that
$\FF$ is tame at all places in $|D|$, $F^*\GG_i$ and $F^*\GG_j$ are
tame at $\infty\in\P^1$ and so $\psi$ must be tame there as well.  At
each place in $S$, the stalk of $F^*\GG_i$, viewed as representation
of the local inertia group, is the direct sum of a ramified
representation of some dimension $e$ and some copies of the trivial
representation and we always have the inequality $e\le r=\rk\FF$.  But
$N_i=\rk\GG_i$ is large (at least $(2g-2+\deg(D))\rk\FF$) and so
$N_i>2e$.  Similarly for $F^*\GG_j$.  This implies that $\psi$ must be
unramified at every place in $S$.  Thus $\psi$ is a character of
$\pi_1(U)$ which is unramified at every place of $S=\P^1\setminus
\left(U\cup\{\infty\}\right)$ and which is tame at $\infty$.  Since
$\A^1$ is ``tamely simply connected'' (i.e.,
$\pi_1^\text{tame}(\A^1)=0$), we must have that $\psi$ is trivial.
This means that $F^*\GG_i$ and $F^*\GG_j$ are already isomorphic on
$U$ which implies, by the argument of part (1), that $i\cong j\pmod
d$.

The argument when $f^*\GG_i$ is dual to $f^*\GG_j$ is quite similar
and will be omitted.  This completes the proof of the proposition.
\end{proof}

                \section{Monodromy groups}\label{s:monodromy}

\subsection{Definitions}
As usual, we write $\Xbar$ for $X\times\Fqbar$.  If $i\in\Z/d\Z$, set
$X_i=X\times\Fq(\mu_{d/(d,i)})$.  In the previous section we defined
sheaves $\GG_o$ on $X$ for each orbit $o\subset\Z/d\Z$ of
multiplication by $q$ and $\GG_i$ on $X_i$ for each $i\in\Z/d\Z$ and
we proved that $\GG_o\cong\oplus_{i\in o}\GG_i$ on $X_i$ and that
$Fr^*(\GG_i)\cong\GG_{iq}$.

These sheaves can be viewed as representations of various fundamental
groups.  More precisely, fix a geometric generic point
$\overline{\eta}$ of $\Xbar$; we also write $\overline{\eta}$ for the
induced geometric generic points of $X_i$ and $X$.  Consider the
arithmetic and geometric fundamental groups
\begin{equation*}
\pi_1(\Xbar,\overline{\eta})\subset
\pi_1(X_i,\overline{\eta})\subset
\pi_1(X,\overline{\eta}).
\end{equation*}
All three groups act on the stalk at $\overline{\eta}$ of $\GG_o$ and
the two smaller groups act on the stalk at $\overline{\eta}$ of
$\GG_i$, so we have homomorphisms
\begin{equation*}
  \tau_o:\pi_1(X,\overline{\eta})\to\aut(\GG_{o,\overline{\eta}})
\end{equation*}
and
\begin{equation*}
  \tau_i:\pi_1(X_i,\overline{\eta})\to\aut(\GG_{i,\overline{\eta}}).
\end{equation*}
Here $\aut(\GG_{o,\overline{\eta}})$ is viewed as the set of $E$
points of an algebraic group over $E$, isomorphic of course to
$\GL_{\rk\GG_o}$, and similarly with $\aut(\GG_{i,\overline{\eta}})$.
The isomorphism $Fr^*(\GG_i)\cong\GG_{iq}$ implies that if
$\Phi\in\pi_1(X,\overline{\eta})$ is an element inducing the geometric
($q^{-1}$-power) Frobenius automorphism of $\Fqbar$, then
$\tau_i^{\Phi}\cong\tau_{iq}$.

We define the arithmetic monodromy group $G_o^{\text{arith}}$ to be
the Zariski closure of the image of $\tau_o$ and the geometric
monodromy group $G_o^{\text{geom}}$ to be the Zariski closure of
$\tau_0(\pi(\Xbar,\overline{\eta}))$.  Similarly, $G_i^{\text{arith}}$
is by definition the Zariski closure of the image of $\tau_i$ and
$G_i^{\text{geom}}$ is by defintion the Zariski closure of
$\tau_i(\pi(\Xbar,\overline{\eta}))$.  Deligne proved
\cite[1.3.9]{WeilII} that $G_o^{\text{geom}}$ and $G_i^{\text{geom}}$
are (not necessarily connected) semisimple algebraic groups over $E$.
We will prove below that (after a suitable twist) the indices of
$G_o^{\text{geom}}\subset G_o^{\text{arith}}$ and
$G_i^{\text{geom}}\subset G_i^{\text{arith}}$ are finite, so the
arithmetic groups are also semisimple.

\subsection{Katz' monodromy calculation}\label{ss:KatzCalc}
The main theorem of \cite{Katz} is a calculation of the groups
$G_i^{\text{geom}}$.  Under the hypotheses of Sections~\ref{ss:data}
and \ref{ss:typical-local}-\ref{ss:special-local} and Theorem~\ref{thm:main}
(in particular, $\rho$ is everywhere tame or tame at places in $S_r$ and 
$p>\deg\rho+2$, and  $\deg(D)$ is large), we have that $G_i^{\text{geom}}$ 
is isomorphic to:
\begin{equation*}
\begin{cases}
\Sp(N_i)&\text{if $d/i=2$ and $\FF$ is orthogonally self-dual}\\
\Or(N_i)\text{ or }\SO(N_i)&
      \text{if $d/i=2$ and $\FF$ is symplectically self-dual}\\
\SL^{(\nu_i)}(N_i)&
      \text{if $d/i\neq2$ or $\FF$ is not self-dual.}
\end{cases}
\end{equation*}
(\cite[5.5.1 case (1b) and 5.7.1]{Katz} 
Here $N_i$ is the rank of $\GG_i$, $\GL(N_i)$, $\Sp(N_i)$, $\Or(N_i)$, and $\SO(N_i)$ refer to the
standard general linear, symplectic, orthogonal, and special
orthogonal groups over $E$, and
\begin{equation*}
\SL^{(\nu_i)}(N_i)=\{g\in \GL(N_i)|(\det g)^{\nu_i}=1\}.
\end{equation*}
In the second case, if $N_i$ is odd, then $G_i^{\text{geom}}=\Or(N_i)$.
Note that the connected component of $G_i^{\text{geom}}$ is either
$\Sp(N_i)$, $\SO(N_i)$, or $\SL(N_i)$.

\subsection{Structure of $G_o^{\text{geom},0}$}
In this subsection we apply the results of Katz to determine the
connected component of the algebraic group $G_o^{\text{geom}}$.  If
$o=\{i\}$, then $G_o^{\text{geom}}$ was already determined by Katz, as
in the previous subsection.  So for the rest of this subsection we
assume that $\#o>1$ and thus $d_o>2$.  Because of the decomposition
$\GG_o\cong\prod_{i\in o}\GG_i$ on $\Xbar$, we have
\begin{equation*}
G_o^{\text{geom}}\subset\prod_{i\in o}\aut(\GG_i)
\cong\prod_{i\in o}\GL(N_i).
\end{equation*}
Let $p_i:G_o^{\text{geom}}\to\aut(\GG_i)$ be the projection onto the
$i$-th factor.  It is elementary from the definitions that
$p_i(G_o^{\text{geom}})$ is contained in $G_i^{\text{geom}}$. Since
this image is Zariski dense and the image of a morphism of algebraic
groups is closed \cite[I.1.4a]{Borel}, we have
$p_i(G_o^{\text{geom}})=G_i^{\text{geom}}$.  It follows
from \cite[I.1.4b]{Borel} that
$p_i(G_o^{\text{geom},0})=G_i^{\text{geom},0}\cong \SL(N_i)$ where the
superscript 0 indicates the connected component of the identity.

\begin{prop}\label{prop:o-geom-monodromy}
Let $o$ be an orbit of multiplication by $q$ on $\Z/d\Z$.
\begin{enumerate}
\item If $\FF$ is not self-dual on $\Xbar$ or if $o\neq-o$ then the
projections $p_i$ induce an isomorphism 
\begin{equation*}
G_o^{\text{geom},0}\cong\prod_{i\in o}G_i^{\text{geom},0}
\cong \prod_{i\in o}\SL(N_i).
\end{equation*}
\item If $\FF$ is self dual on $\Xbar$ and $o=-o$, let $S\subset o$ be
a set of representatives for $o$ modulo $\pm1$.  Then the projections
$p_i$ induce an isomorphism
\begin{equation*}
G_o^{\text{geom},0}\cong\prod_{i\in S}G_i^{\text{geom},0}
\cong \prod_{i\in S}\SL(N_i).
\end{equation*}
If $j\not\in S$  then in terms of suitable bases, the projection
$p_j:G_o^{\text{geom},0}\to G_i^{\text{geom},0}$  sends a tuple of
matrices $(A_i)_{i\in S}$ to $A_{-j}^\vee={}^t(A_{-j})^{-1}$.
\end{enumerate}
\end{prop}

\begin{proof}
  Let $\ggg_o$ and $\ggg_i$ denote the Lie algebras of
  $G_o^{\text{geom},0}$ and $G_i^{\text{geom},0}$, which are
  semisimple.  The projections $p_i$ induce surjections
  $dp_i:\ggg_o\to\ggg_i\cong sl(N_i)$.  Let $\hh'_i=\ker dp_i$.  Since
  $\ggg_o$ is semisimple, we have a decomposition
  $\ggg_o=\hh_i\oplus\hh'_i$ where $\hh_i$ is an ideal mapping
  isomorphically onto $\ggg_i$.  Now take $j\in o$, $j\not\equiv i\pmod
  d$ and consider $dp_j$ restricted to $\hh_i$, so that
  $dp_j|_{\hh_i}:\hh_i\to\ggg_j$.  The source and target of this
  homomorphism are both simple (they are both isomorphic to $sl(N_i)$)
  so $p_j|_{\hh_i}$ is either 0 or an isomorphism.  Let us suppose for
  a moment that it is an isomorphism and define
  $d\phi_{ji}=dp_j\compose(dp_i|_{\hh_i})^{-1}:\ggg_i\isoto\ggg_j$.  Since
  $\SL(N_i)$ is simply connected we may integrate $d\phi_{ji}$ to an
  isomorphism $\phi_{ji}:G_i^{\text{geom},0}\to G_j^{\text{geom},0}$.
  Let $Y\to\Xbar$ be the finite \'etale cover which trivializes
  $\det\GG_i$ for all $i\in o$ and let $\overline{\eta}_Y$ be a geometric
  generic point of $Y$.  Then we have a commutative diagram
\begin{equation*}
\xymatrix{&&G_i^{\text{geom},0}\ar[dd]^{\phi_{ji}}\\
\pi_1(Y,\overline{\eta}_Y)\ar[r]&G_o^{\text{geom},0}
\ar[ur]^{p_i}\ar[dr]_{p_j}\\
&&G_j^{\text{geom},0}}
\end{equation*}
Now it is well known that the only automorphisms of $\SL$ are inner or
inner composed with $A\mapsto A^\vee={}^tA^{-1}$. (This follows easily
from the Lie algebra version, which is \cite[Chap.~IX, Thm.~5, p.
283]{Jacobson}.) Thus if $dp_j|_{\hh_i}$ is an isomorphism, then
$\tau_i$ and $\tau_j$ become isomorphic or contragredient over $Y$;
equivalently, $\GG_i$ and $\GG_j$ become isomorphic or dual on $Y$.
But Proposition~\ref{prop:GGis-distinct}, the first case is impossible
($j\not\equiv i$) and the second is impossible unless $j\equiv-i$ and
$\FF$ is self-dual.  Thus, under the hypotheses of (1), $p_j|_{\hh_i}$
must be zero for all $i\not\equiv j$.  From this we easily conclude
that $\ggg_o\cong\prod_{i\in o}\ggg_i$.  This implies that the
projections $p_i$ induce a local isomorphism
$G_o^{\text{geom},0}\to\prod_{i\in o}G_i^{\text{geom},0}$ and since
the target is simply connected, they in fact induce an isomorphism.
This concludes the proof of (1).

Under the hypotheses of (2), $dp_j|_{\hh_i}$ is zero if $j\not\equiv-i$
and we know (by the trivial part of
Proposition~\ref{prop:GGis-distinct}) that if $j\equiv-i$ then
$dp_j|_{\hh_i}$ is an isomorphism.  As in part (1), we easily conclude
that the $dp_i$ induce an isomorphism $\ggg_o\cong\prod_{i\in S}\ggg_i$
and thus the $p_i$ induce an isomorphism
$G_o^{\text{geom},0}\cong\prod_{i\in S}G_i^{\text{geom},0}$.  Also,
there is an isomorphism $\phi_{ji}$ as in the displayed equation
above, and since $i\not\equiv j$, this isomorphism is not inner, so in
terms of suitable bases it is $A_i\mapsto A_i^\vee$.  This completes
the proof of the proposition.
\end{proof}

The last sentence of the proposition can also be deduced by explicit
matrix calculations, as in \ref{ss:reduced-cp} below.

\subsection{Structure of $G_o^{\text{geom}}$}
Let $\Phi^\text{geom}_o$ and $\Phi^\text{geom}_i$ denote the groups of
connected components of $G_o^{\text{geom}}$ and $G_i^{\text{geom}}$
respectively.  By Katz' monodromy calculation, we have
$\Phi^\text{geom}_i\cong\mu_{\nu_i}$, the roots of unity of order
$\nu_i$ for some integer $\nu_i$.  The isomorphism
$\tau_i^{\Phi}\cong\tau_{iq}$ and the fact that
$\pi_1(\Xbar,\overline{\eta})$ is a normal subgroup of
$\pi_1(X,\overline{\eta})$ imply that $\nu_i$ is independent of $i$
for $i$ running through a fixed orbit $o$; let $\nu_o$ denote the
common value of the $\nu_i$. Thus $\Phi^\text{geom}_o$ is a subgroup
of $\prod_{i\in o}\Phi^\text{geom}_i=\mu_{\nu_o}^{a_o}$; the
isomorphism $\tau_i^{\Phi}\cong\tau_{iq}$ implies that this subgroup
is invariant under cyclic permutation ($i\mapsto iq$) of the factors.
Also, since $p_i(G_o^{\text{geom}})=G_i^{\text{geom}}$, the projection
$p_i$ induces a surjection $\Phi^\text{geom}_o\to\Phi^\text{geom}_i$
for each $i$.  So in all, we have that $\Phi^\text{geom}_o$ is a
subgroup of $\prod_{i\in o}\mu_{\nu_o}^{a_o}$ which maps surjectively
onto each factor and which is invariant under cyclic permutation of
the factors.

\subsection{Arithmetic monodromy groups}
Our next goal is to determine the structure of the arithmetic
monodromy group $G_o^\text{arith}$ or rather of a twisted version of
it.  This will amount to determining its component group.

\subsubsection{}
Given an $\ell$-adic unit $\beta\in \O_E^\times$ there is a continuous
Galois representation $\gal(\Fqbar/\Fq)\to E^\times$ which sends $Fr$
to $\beta$.  We denote the corresponding lisse sheaf on $\spec\Fq$, as
well as its pull back to various schemes over $\Fq$, by
$\beta^\text{deg}$.  

If we write $\GG(\FF,o)$ for the sheaf on $X$ defined above using the
sheaf $\FF$ on $\Curve$ and the orbit $o\subset\Z/d\Z$, then the
projection formula implies that we have a canonical isomorphism
$\GG(\FF\tensor\beta^\text{deg},o)\cong\GG(\FF,o)\tensor\beta^\text{deg}$
of sheaves on $X$.  Define $G^\text{geom}(\beta)$ and
$G^\text{arith}(\beta)$ to be the geometric and arithmetic monodromy
groups associated to $\GG(\FF,o)\tensor\beta^\text{deg}$.  Since
$\beta^\text{deg}$ is trivial on $\Xbar$ we have
$G^\text{geom}(\beta)=G^\text{geom}$.  On the other hand,
$G^\text{arith}(\beta)$ will in general differ from $G^\text{arith}$.

The connection with $L$-functions also changes: we have
\begin{align*}
\det\left(1-T\,Fr_{n,f}\left|(\GG(\FF,o)\tensor\beta^\text{deg})_{\overline
f}\right.\right)
&=\det\left(1-\beta T\,Fr_{n,f}\left|\GG(\FF,o)_{\overline f}\right.\right)\\
&=L(F_n,\rho\tensor\sigma_{o,f},\beta T)
\end{align*}
for all $f\in X(\Fqn)$.

\subsubsection{}
For the rest of this section, we view $\FF$ and $o$ as being fixed and
we drop them from the notation.  Let $\Gamma(\beta)$ be defined as
$G^\text{arith}(\beta)/G^\text{geom}(\beta)$ and consider the
following commutative diagram, where the columns and rows are exact by
definition.
\begin{equation*}
  \xymatrix{&0\ar[d]&0\ar[d]\\
0\ar[r]&G^{\text{geom},0}\ar[r]\ar[d]&G^{\text{arith},0}(\beta)\ar[d]\\
0\ar[r]&G^\text{geom}\ar[r]\ar[d]&G^{\text{arith}}(\beta)\ar[r]\ar[d]&
    \Gamma(\beta)\ar[r]&0\\
&\Phi^\text{geom}\ar[r]\ar[d]&\Phi^\text{arith}(\beta)\ar[d]\\
&0&0\\}
\end{equation*}

The next proposition says that for a suitable $\beta$,
$G^{\text{geom}}$ has finite index in $G^{\text{arith}}(\beta)$ and so
$G^{\text{geom},0}=G^{\text{arith},0}(\beta)$. Thus for such a $\beta$
we can complete the diagram into the following, where all rows and
columns are exact:
\begin{equation*}
  \xymatrix{&0\ar[d]&0\ar[d]\\
0\ar[r]&G^{\text{geom},0}\ar@{=}[r]\ar[d]&
    G^{\text{arith},0}(\beta)\ar[d]\ar[r]&0\ar[d]\\
0\ar[r]&G^\text{geom}\ar[r]\ar[d]&G^\text{arith}(\beta)\ar[r]\ar[d]&
    \Gamma(\beta)\ar[r]\ar@{=}[d]&0\\
0\ar[r]&\Phi^\text{geom}\ar[r]\ar[d]&\Phi^\text{arith}(\beta)\ar[r]\ar[d]&
    \Gamma(\beta)\ar[r]\ar[d]&0\\
&0&0&0\\}
\end{equation*}

\begin{prop}\label{prop:good-beta}
 Expanding $E$ if necessary, there exists a
  $\beta\in\O_E^\times$ such that the conditions below hold.  
\begin{enumerate}
\item[(a)] The arithmetic monodromy group associated to
  $\GG(\FF,o)\tensor\beta^\text{deg}$ contains $G_o^{\text{geom}}$ as
  a finite index subgroup, i.e., $\Gamma(\beta)$ is finite.
\item[(b)] $\Gamma(\beta)$ is cyclic of order $a_o$ or $2a_o$.  Its order
  is $2a_o$ if and only if either (i) $a_o>1$, $o=-o$, $\FF$ is
  orthogonally self-dual, $N_i$ is odd and $\nu_i$ is odd; or (ii)
  $o=\{d/2\}$, $\FF$ is symplectically self-dual,
  $G^\text{geom}=\SO(N_{d/2})$ and $G^\text{arith}=\Or(N_{d/2})$.
\item[(c)] $\Phi^\text{arith}(\beta)$ is the semi-direct product
  $\Phi^\text{geom}\sdp\Gamma(\beta)$ where the action of
  $\Gamma(\beta)$ on $\Phi^\text{geom}\subset\prod\mu_{\nu_o}$ is by
  cyclic permutation of the factors.
\end{enumerate}
If $\FF$ is self-dual (of weight $w=-1$) on $X$ and $o=-o$, then we
may take $\beta=1$.
\end{prop}

\begin{proof}
  First suppose that $\FF$ is self-dual (of weight $w=-1$) and
  $o=\{d/2\}$.  Then $\GG(\FF,o)$ self-dual and so
  $G_o^{\text{arith}}$ is {\it a priori\/} contained in an orthogonal
  or symplectic group.  But as we have seen, $G_o^{\text{geom}}$ is
  the full symplectic group or contains the special orthogonal group,
  so (a), (b), and (c) are clear in this case.
  
  Next, we make an observation about determinants.  Let $\Phi$ be an
  element of $\pi_1(X,\overline\eta)$ inducing the geometric Frobenius
  on $\Fqbar$.  Then since
  $\tau_i^\Phi\cong\tau_{iq}$, we have that
  $\det\tau_i(\Phi^{a_o})=\det\tau_{iq}(\Phi^{a_o})$ and so
  $\det\tau_i(\Phi^{a_o})$ is independent of $i\in o$.  This means
  that there is a $\beta\in \O_E^\times$ such that
  $\det\left(\tau_i\tensor\beta^\text{deg}\right)(\Phi^{a_o})=1$ for
  all $i\in o$.

Now assuming that $o\neq\{d/2\}$ or $\FF$ is not self-dual, we have
seen that the groups $G^{\text{geom},0}$ are all $\SL(N_i)$ and so
$\left(\tau_o\tensor\beta^\text{deg}\right)(\Phi^{a_o})$ lies in
$G^\text{geom}$.  Since
$\left(\tau_o\tensor\beta^\text{deg}\right)(\Phi)$ generates
$\Gamma(\beta)$, this proves that $\Gamma(\beta)$ is finite cyclic of
order dividing $a_o$, indeed of order exactly $a_o$ since
$\tau_o(\Phi)$ permutes the factors of $\GG_o\cong\oplus_{i\in
  o}\GG_i$ cyclically.  It also shows that $\Phi^\text{arith}(\beta)$
is a semi-direct product, i.e., the lower row of our diagram is split
exact.  That the action of $\Gamma(\beta)$ on $\Phi^\text{geom}$ is as
asserted follows easily from the formula $\tau_i^\Phi\cong\tau_{iq}$.

This completes the proof of the proposition except in the case where
$\FF$ is self-dual (of weight $w=-1$), $a_o>1$, and $o=-o$, in which
case we insist that $\beta=1$ and we have to show that $\Gamma$ is
finite of order $a_o$ or $2a_o$.  But under these hypotheses,
$\GG(\FF,i)$ and $\GG(\FF,-i)$ are dual on $X_i$ and so
$\det\tau_i(\Phi^{a_o}) =\left(\det\tau_{-i}(\Phi^{a_o})\right)^{-1}$.
Since $\det\tau_i(\Phi^{a_o})$ is independent of $i\in o$, this
implies that these determinants are $\pm1$. A matrix calculation (see
\ref{ss:reduced-cp} below) shows that this determinant is in fact 1 if
$\FF$ is symplectically self-dual, and it is $(-1)^{N_i}$ if $\FF$ is
orthogonally self-dual.  In light of
Proposition~\ref{prop:o-geom-monodromy}, this implies that
$\tau_o(\Phi^{2a_o})$ lies in $G_o^{\text{geom}}$ and
$\tau_o(\Phi^{a_o})$ lies in $G_o^{\text{geom}}$ except in the cases
mentioned in part (2).  This completes the proof of the proposition.
\end{proof}

Note that for $\beta$ as in the proposition, the twisted sheaf
$\GG(\FF\tensor\beta^\text{deg},o)$ is again $\iota$-pure of weight 0.

\subsection{Reduced characteristic polynomials}\label{ss:reduced-cp}
Let $\GL(N)$ denote the general linear group over some field and let
$G\subset \GL(N)$ be a closed algebraic subgroup.  We define the
reduced characteristic polynomial function as follows.  For each
irreducible component of $G$, let $P_0(T)$ be the gcd of the
(reversed) characteristic polynomials of the elements of that
component.  Then define $P^\text{red}_g(T)$, the reduced
characteristic polynomial of $g\in G$, to be the usual (reversed)
characteristic polynomial, divided by the gcd $P_0$ for the component
in which $g$ lies.  The key property of the reduced characteristic
polynomial is that if $\alpha$ is any element of the ground field,
then the set of $g\in G$ such that $P^\text{red}_g(\alpha)=0$ is a
Zariski closed subset which contains no irreducible components of $G$.
In particular, if the field is $\C$, this set has Haar measure zero.

Now we compute the reduced characteristic polynomials (or rather the gcd's
$P_0$) for various components of the groups
$G^\text{arith}_o(\beta)\subset\GL(\GG_o\tensor\beta^\text{deg})$.
 
\subsubsection{}\label{sss:quad-orth}
If $o=\{d/2\}$ and $\FF$ is orthogonally self-dual, then
$G^\text{arith}=\Sp(N_{d/2})$ and $P_0(T)=1$.

\subsubsection{}\label{sss:quad-symp}
If $o=\{d/2\}$ and $\FF$ is symplectically self-dual, then
$G^\text{arith}$ is either $\SO(N_{d/2})$ or $\Or(N_{d/2})$.  In the
former case $N_{d/2}$ is necessarily even (see \ref{ss:KatzCalc}) and
so $P_0(T)=1$.  In the latter, there are two cases depending on the
parity of $N_{d/2}$.  If $N_{d/2}$ is even, $P_0(T)$ is 1 on
$\SO(N_{d/2})$ and $1-T^2$ on $\Or_-(N_{d/2})$.  If $N_{d/2}$ is odd,
$P_0(T)$ is $1-T$ on $\SO(N_{d/2})$ and $1+T$ on $\Or_-(N_{d/2})$.

\subsubsection{}\label{sss:non-self-dual}
Next we consider the case where $o\neq-o$ or $\FF$ is not self-dual.
Here we claim that $P_0(T)=1$ on every component of
$G^\text{arith}_o$.  Recall that components of $G^\text{arith}_o$ are
indexed by tuples $((\zeta_i)_{i\in o},b)$ where
$(\zeta_i)\in\prod_{i\in S}\mu_{\nu_o}$ and $b\in\Z/a_o\Z$.  For
convenience, we prove the assertion only for components where $b$ is a
generator of $a_o$; the other cases are similar but would require more
notational complexity.  Let us fix $j\in o$ and a basis of
$\GG_{j,\overline{\eta}}$.  We extend this to a basis of
$\GG_{o,\overline{\eta}}$ by applying
$\tau_o(\Phi^b),\tau_o(\Phi^{2b}), \dots$ to the original basis.  In
terms of this basis, the matrix of an element $g$ of the component
indexed by $((\zeta_i)_{i\in o},b)$ is a ``block cyclic permutation
matrix," i.e., it has the form
\begin{equation*}
  \begin{pmatrix}
0&0&0&\cdots&A_{jq^{a_o-1}}\\
A_j&0&0&\cdots&0\\
0&A_{jq}&0&\cdots&0\\
\vdots&\vdots&\vdots&\ddots&\vdots\\
0&0&0&\cdots&0
\end{pmatrix}
\end{equation*}
where the blocks are $N_o\times N_o$ and $\det A_i=\zeta_i$ for all
$i\in o$.  Moreover, as $g$ varies through the component, the matrices
$A_i$ vary (independently) over all matrices with these determinants.
(This comes from Proposition~\ref{prop:o-geom-monodromy} above.)  It
follows easily that $P_0(T)=1$.

\subsubsection{}\label{sss:non-Galois-self-dual}
Finally, we consider the case where $o=-o$, $a_o>1$, and $\FF$ is
self-dual (of weight $w=-1$) and we restrict to components indexed by
$((\zeta_i)_{i\in o},b)$ where $b$ is prime to $a_o$.  In this case,
we claim that if $N_o$ is even, then $P_0(T)=1$ on every such
component whereas if $N_o$ is odd, $P_0(T)=(1-T^{a_o})$ if $\FF$ is
symplectically self-dual and $P_0(T)=(1+T^{a_o})$ if $\FF$ is
orthogonally self-dual.

It will be convenient to argue with matrices.  (This is
essentially the same argument as in \ref{sss:induced-zeroes-proof}.)
Choosing a basis as above, elements $g$ of the component indexed by
$((\zeta_i)_{i\in o},b)$ are block cyclic permutation matrices, as
above, but as we will see, there are relations among the $A_i$.  To
see this, note that the matrix of the form on
$\GG_{o,\overline{\eta}})$ (which is orthogonal resp.~symplectic if
$\FF$ is symplectically resp.~orthogonally self-dual) is
\begin{equation*}
\begin{pmatrix}
0&I_{N_oa_0/2}\\
\epsilon I_{N_oa_0/2}&0
\end{pmatrix}
\end{equation*}
where $I_{N_oa_0/2}$ denotes the identity matrix of size $N_oa_o/2$
and $\epsilon=-\sgn(\rho)$, i.e., $\epsilon=1$ if $\rho$ is
symplectically self-dual and $-1$ if it is orthogonally self-dual.
Writing out the condition that $g$ respects the form, we find that
$A_{-jq^{kb}}=A_{jq^{kb}}^\vee$ for $k=0,\dots,a_o/2-2$ and
$A_{-q^{(a_o/2-1)b}}=\epsilon A_{-q^{(a_o/2-1)b}}^\vee$.  By
Proposition~\ref{prop:o-geom-monodromy}, other than these
restrictions, the matrices vary freely among those with determinants
$(\zeta_i)$.  

Now since the matrix of $g$ is block cyclic permutation, its
eigenvalues are all of the $a_o$-th roots of those of $g^{a_o}$.  The
matrix calculation above shows that the matrix of $g^{a_o}$ is block
diagonal with blocks of the form
\begin{equation*}
\epsilon A_{jq^{(a_o/2-1)b}}^\vee A_{jq^{(a_o/2-2)b}}^\vee\cdots A_{j}^\vee
A_{jq^{(a_o/2-1)b}}A_{jq^{(a_o/2-2)b}}\cdots A_{j} 
\end{equation*}
which is of the form $\epsilon B^\vee B$.  (To tie up a loose end in
Proposition~\ref{prop:good-beta}, note that these blocks have
determinant $\epsilon^{N_o}$.)  By Lemma~\ref{lemma:BcheckB}, if $N$
is odd, all the matrices $B^\vee B$ have 1 as an eigenvalue,
generically of multiplicity 1 and have no other shared eigenvalues.
If $N$ is even then there are no shared eigenvalues.  This completes
the proof of our claims about $P_0(T)$.

                    \section{Equidistribution}\label{s:equidistribution}

In this section we fix the sheaf $\FF$ and the orbit $o$ and then
choose a $\beta$ as in Proposition~\ref{prop:good-beta}.  We will drop
this data from the notation and so just write $G^\text{arith}$ and
$G^\text{geom}$ for the arithmetic and geometric monodromy groups
attached to $\GG(\FF,o)\tensor\beta^\text{deg}$.  Also, $\Gamma$ will
denote $G^\text{arith}/G^\text{geom}$.

\subsection{Maximal compact subgroups}
Using the embedding $E\into\Qlbar\cong\C$ we may extend scalars and
define semisimple algebraic groups $G^\text{arith}_{/\C}$ and
$G^\text{geom}_{/\C}$ over $\C$.  Let $G^\text{arith}(\C)$ and
$G^\text{geom}(\C)$ denote their complex points, which we regard as
complex semisimple Lie groups.

We will denote by $K^\text{arith}$ and $K^\text{geom}$ maximal compact
subgroups of $G^\text{arith}(\C)$ and $G^\text{geom}(\C)$.  By Weyl's
``unitarian trick,'' $K^\text{arith}$ is Zariski dense in
$G^\text{arith}_{/\C}$ and so $K^\text{arith}/K^\text{geom}\cong
G^\text{arith}/G^\text{geom}\cong\Gamma$ is a finite cyclic group.
Also, the group of components of $K^\text{arith}$ and $K^\text{geom}$
are the same as those of $G^\text{arith}$ and $G^\text{geom}$.

We define the reduced characteristic polynomials $P^\text{red}_k(T)$
for $k\in K^\text{arith}$ as in \ref{ss:reduced-cp} above (dividing
the usual reversed characteristic polynomial by the gcd of the
characteristic polynomials over each connected component).  Again
because $K^\text{arith}$ is Zariski dense in $G^\text{arith}_{/\C}$,
the reduced characteristic polynomials for $K^\text{arith}$ are just
the restrictions of the reduced characteristic polynomials from
$G^\text{arith}$ (via the embedding $\iota$).

\subsection{Haar measures}
Fix an element $\gamma\in\Gamma\cong K^\text{arith}/K^\text{geom}$ and
let $K^\text{arith}_\gamma$ denote the inverse image of $\gamma$ in
$K^\text{arith}$.  We denote the set of conjugacy classes of
$K^\text{arith}$ which meet $K^\text{arith}_\gamma$ by
$K^{\text{arith},\#}_\gamma$; since $\Gamma$ is abelian, this is just
the quotient of $K^\text{arith}_\gamma$ by the conjugation action of
$K^\text{geom}$.

Let $d\mu_{\text{Haar},\gamma}$ be the $K^\text{geom}$-translation
invariant measure on $K^\text{arith}_\gamma$ of total mass 1.  (We may
take the left or right invariant measure as either is bi-invariant.)
Let $d\mu_{\text{Haar},\gamma}^\#$ be its push-forward onto
$K^{\text{arith},\#}_\gamma$.  The main equidistribution statement
will be that a suitably normalized sum of point masses corresponding
to Frobenius elements converges to the measure
$d\mu_{\text{Haar},\gamma}^\#$.

\subsection{Frobenius classes}
Let $f$ be an element of $X(\Fqn)$ and denote as usual a corresponding
Frobenius element (defined up to conjugacy) by
$Fr_{n,f}\in\pi_1(X,\overline{\eta})$.  The monodromy representation
$\tau_o$ gives us an element (up to conjugacy) $\tau_o(Fr_{n,f})\in
G^\text{arith}(E)\into G^\text{arith}(\C)$ and we denote its
``semi-simple part'' (obtained from a Jordan form by throwing away the
off-diagonal terms) by $\tau_o(Fr_{n,f})^{ss}$.  Because
$\GG(\FF,o)\tensor\beta^\text{deg}$ is $\iota$-pure of weight 0, the
eigenvalues of $\tau_o(Fr_{n,f})^{ss}$ lie on the unit circle, and so
$\tau_o(Fr_{n,f})^{ss}$ is conjugate to an element of
$K^\text{arith}$.  The $K^\text{arith}$-conjugacy class of this
element is well-defined and we denote it by $\theta(f,n)$.

Note that the image of $\tau_o(Fr_{n,f})$ in
$\Gamma=G^\text{arith}/G^\text{geom}$ (which we have seen is
$\Z/a_o\Z$ or $\Z/2a_o\Z$) is just the class $\gamma$ of $n$.  Thus as $f$
varies through $X(\Fqn)$, the classes $\theta(f,n)$ all lie in the set
of classes of $K^\text{arith}$ over a fixed element $\gamma\in\Gamma$, i.e.,
in $K^{\text{arith},\#}_\gamma$.

\numberwithin{equation}{subsection}
\subsection{Equidistribution}
For each integer $n$ we have the finite
set of points $X(\Fqn)$ and the corresponding conjugacy classes
$\theta(f,n)$ in $K^{\text{arith},\#}$.  We define a measure $d\mu_n$ on
the set of conjugacy class $K^{\text{arith},\#}$ by
averaging the point masses at the various classes $\theta(f,n)$ for
$f\in X(\Fqn)$.  Thus, if $\phi$ is a class function on $K^\text{arith}$,
\begin{equation*}
\int_{K^{\text{arith},\#}}\phi\,d\mu_n
=\frac1{\#X(\Fqn)}\sum_{f\in X(\Fqn)}\phi(\theta(f,n)).
\end{equation*}
Note that this measure is supported on
$K^{\text{arith},\#}_\gamma$ where $\gamma$ is the class of $n$ in
$\Gamma$.

The basic equidistribution statement is that the measures $d\mu_n$
converge weakly to $d\mu^\#_{\text{Haar},\gamma}$ as $n\to\infty$
through a fixed class in $\Gamma$.  In other words, if $\phi$ is a
continuous class function on $K^\text{arith}$, we have
\begin{equation}\label{eqn:weak-conv}
\lim_{\substack{n\to\infty\\ [n]=\gamma}}\int_{K^{\text{arith},\#}}\phi\,d\mu_n=
\int_{K^{\text{arith},\#}}\phi\,d\mu^\#_{\text{Haar},\gamma}
\end{equation}

This result is \cite[9.7.10]{KatzSarnak} (with $S=\spec\Fq$) which is a
mild generalization of \cite[3.5.3]{WeilII}.
\numberwithin{equation}{subsubsection}

\subsection{Good test functions}
We will apply the equidistribution statement \ref{eqn:weak-conv} to a
well-chosen test function to conclude that for large enough $n$, there
are many $f\in X(\Fqn)$ such that a given $\alpha$ is not a root of
the reduced characteristic polynomial $P^\text{red}_{\theta(f,n)}$.

Let $K^\text{arith}_\alpha$ denote the subset of elements $k\in K^\text{arith}$
where $P^\text{red}_k(\alpha)=0$.  This is a Zariski closed subset
which is a proper subset of each component of $K^\text{arith}$.

\begin{prop} 
For every $\epsilon>0$ there exist smooth class functions
$f_\alpha:K^\text{arith}\to\R$ indexed by $\alpha\in S^1$ such that 
\begin{enumerate}
\item[(a)] $0\le f_\alpha(k)\le 1$ for all $k\in K^\text{arith}$, all
  $\alpha\in S^1$.
\item[(b)] For all $\alpha\in S^1$, $f_\alpha(k)=1$ for all $k\in
K^\text{arith}_\alpha$.
\item[(c)]  There exists $n_0$ such that for each $\gamma\in\Gamma$ and all
  $n>n_0$ in the class of $\gamma$, $\int_{K^{\text{arith},\#}_\gamma}
  f_\alpha\,d\mu_n<\epsilon$ for all $\alpha\in S^1$.
\end{enumerate}
\end{prop}

\begin{proof}
Let $f_\alpha$ be defined by the formula
\begin{equation*}
f_\alpha(k)=e^{-C\left|P^\text{red}_k(\alpha)\right|^2}
\end{equation*}
where $P^\text{red}_k$ is the reduced characteristic polynomial of $k$
and $C$ is a positive real number.  Clearly $f_\alpha$ is a smooth
class function of $k$ which satisfies the first two requirements of
the proposition.  

Because $f_\alpha$ vanishes on a proper Zariski closed subset of each
component of $K^\text{arith}$ (i.e., on a set of Haar measure zero)
and $S^1$ is compact, we can choose one $C$ so that
\begin{equation*}
\int_{K^{\text{arith},\#}_\gamma}f_\alpha\,d\mu^\#_{\text{Haar},\gamma}
<\epsilon/2
\end{equation*}
for all $\alpha\in S^1$.

Next, we claim that for sufficiently large $n$,
\begin{equation*}
  \left|
\int_{K^{\text{arith},\#}_\gamma}f_\alpha\,d\mu_n
-\int_{K^{\text{arith},\#}_\gamma}f_\alpha\,d\mu^\#_{\text{Haar},\gamma}
\right|<\epsilon/2
\end{equation*}
for all $\alpha\in S^1$.  For a fixed $\alpha$, this is just our
equidistribution statement~\ref{eqn:weak-conv}.  Again by the
compactness of $S^1$, there is one $n_0$ so that the displayed
inequaltiy holds for all $n>n_0$ in the class of $\gamma$ and all
$\alpha\in S^1$.  Since $\Gamma$ is finite, there is one $n_0$ that
works for all $\gamma$.

Combining the two displayed inequalities shows that the functions
$f_\alpha$ also satisfy the third requirement of the proposition.
\end{proof}

\begin{cor}\label{cor:non-vanishing-density}
Let $X(\Fqn)_{\alpha}$ be the set of elements $f\in X(\Fqn)$ where
$P^\text{red}_{\theta(f,n)}(\alpha)$ vanishes.  Then for every $\epsilon>0$
there exists an integer $n_0$ such that for $n>n_0$
\begin{equation*}
\frac{\# X(\Fqn)_{\alpha}}{\# X(\Fqn)}<\epsilon
\end{equation*}
\end{cor}

\begin{proof}
Indeed, the fraction on the left hand side is bounded above by
$\int_{K^\#}f_\alpha\,d\mu_n$ where $f_\alpha$ is the function
appearing in the proposition.
\end{proof}

\section{End of the proof of the main theorem}\label{s:end-of-proof}
              
We are now in a position to prove the main technical theorem,
Theorem~\ref{thm:main}.  We first give the basic structure of the
argument, then adapt it to the various cases, considering one orbit
$o$ at a time (i.e., part (1) of the theorem).  Then we discuss the
case of several orbits at once (i.e., part (2) of the theorem).

\subsection{The basic argument}
We are given data $\Curve$, $\rho$, $d$, $S_s$, $S_i$,
$S_r$, and $(\alpha_n)_{n\ge1}$ satisfying the hypotheses of
\ref{ss:data}.  By twisting, we may assume that $\rho$ has weight
$w=-1$ and that the $\alpha_n$ all have $\iota$-weight 0.
The representation $\rho$ gives rise to a middle extension
sheaf $\FF$ on $\Curve$.  Fix $o\subset\Z/d\Z$, an orbit for
multiplication by $q$.  Then we choose a divisor $D$ and local
conditions $(S_n,C_{n,w})$ as described in
\ref{ss:typical-local}-\ref{ss:special-local}.  Using $\Curve$ and
$D$, we construct the space $X$ parameterizing certain degree $d$
covers of $\Curve$ and the sheaf $\GG(\FF,o)$.  Then we choose an
$\ell$-adic unit $\beta$ as in \ref{prop:good-beta} and consider
$\GG(\FF\tensor\beta^\text{deg},o)$, as well as its arithmetic
monodromy group $G^\text{arith}$ and its compact form
$K^\text{arith}$. 

Proposition~\ref{prop:good-pts-density} guarantees that for all
sufficiently large $n$, the density of points $f\in X(\Fqn)$
satisfying the local conditions imposed by $(S_n,C_{n,w})$ is bounded
below by some positive constant $C$ independent of $n$.  Applying
Corollary~\ref{cor:non-vanishing-density} with
$\alpha=(\beta\alpha_n)^{-1}$ guarantees that for any $\epsilon>0$,
for all sufficiently large $n$ relatively prime to $a_o$, the density
of points $f\in X(\Fqn)$ such that
$P^\text{red}_{\theta(f,n)}((\beta\alpha_n)^{-1})\neq0$ is
at least $1-\epsilon$.  Since
\begin{equation*}
L(\rho\tensor\sigma_{f,o},F_n,T)=
P_{\theta(f,n)}(\beta^{-1}T)=
P^\text{red}_{\theta(f,n)}(\beta^{-1}T)
\left(\frac{P_{\theta(f,n)}(\beta^{-1}T)}
{P^\text{red}_{\theta(f,n)}(\beta^{-1}T)}\right)
\end{equation*}
the remainder of the argument consists of relating the exceptional
situations to the specific choices of local conditions and the
``forced zeroes,'' i.e., the inverse roots of
$P_k(T)/P^\text{red}_k(T)$ for $k=\theta(f,n)\in K^\text{arith}$.

\subsection{The case where $o\neq-o$ or $\rho$ is not self-dual}
In this case, by \ref{sss:non-self-dual}, $P_k(T)=P^\text{red}_k(T)$
for all $k\in K^\text{arith}$.  Thus there are no ``forced zeroes''
and so the basic argument already proves part (1) of the theorem in
this case.

\subsection{The case where $o=-o$, $a_o>1$, and $\rho$ is self-dual}
In this case, by Proposition~\ref{prop:good-beta} we may take
$\beta=1$.  Let $N_o$ be the rank of $\GG_i$ for any $i\in o$.  Then
we have seen in \ref{sss:non-Galois-self-dual} that if $N_o$ is even
then $P_k(T)=P^\text{red}_k(T)$ for all $k\in K^\text{arith}$, whereas
if $N_o$ is odd, then $P_k(T)=P^\text{red}_k(T)(1+\sgn(\rho) T^{a_o})$
where $\sgn(\rho)$ is $-1$ if $\rho$ is symplectic and $1$ if it is
orthogonal.  In particular, if $N_o$ is even or if $\alpha_n^{a_o}\neq
-\sgn(\rho)$ then the basic argument already suffices.

If hypothesis~\ref{eqn:cond2} fails or if $\rho$
has odd degree, then we have chosen $D$ and $(S_n,C_w)$ 
so that $N_o$ is even.  (These are the choices we
made in \ref{ss:special-local}.)  

So let us assume that hypothesis~\ref{eqn:cond2} holds, $\rho$ has
even degree, and that $\alpha_n^{a_o}=-\sgn(\rho)$, i.e., that we are in
the exceptional situation of type (iii) or (iv).  In these cases,
$P_{k}(T)/P^\text{red}_k(T)=(1+\sgn(\rho) T^{a_o})$ has
$\alpha_n$ as inverse root to order exactly one.

This completes the proof of the theorem in the case
appearing in the section title.

\subsection{The case where $o=\{d/2\}$ and $\FF$ is self-dual}
The argument is quite similar to that in the previous subsection, with
different adjustments for the exceptional cases.

If $\rho$ is orthogonally self-dual, then the monodromy group
$G^\text{arith}$ is symplectic and so by \ref{sss:quad-orth}, the
ratio $P_k(T)/P^\text{red}_k(T)$ is 1.

From now on we assume that $\rho$ is symplectically self-dual so that
the monodromy group is an orthogonal group.  By
Proposition~\ref{prop:good-beta} we may assume $\beta=1$.  If
hypothesis~\ref{eqn:cond1} fails, then by \ref{lemma:global-signs} the
signs in the functional equation vary as $f$ varies.  This implies
that the arithmetic and geometric monodromy groups are both
$\Or(N_o)$.  But then our choice of local conditions in
\ref{ss:special-local} forces $k=\theta(f,n)$ into the component
($\SO(N_o)$ or $\Or_-(N_o)$) where $\alpha_n$ is not an inverse root
of $P_{k}(T)/P^\text{red}_{k}(T)$.

From now on, we also assume that hypothesis~\ref{eqn:cond1} holds, so
that the sign in the functional equation is fixed for a fixed $n$ and
all $f\in X(\Fqn)$ satisfying the local conditions.
Then there
are four cases, depending on the parity of $N=N_o$ and the sign
$W=W(\rho\tensor\chi_f,F_n)$ in the functional equation.  More
precisely, if $N$ is even and $W=1$, then
$P_k(T)/P_k^\text{red}(T)=1$.
If $N$ is even and $W=-1$, then $P_k(T)/P_k^\text{red}(T)=(1-T^2)$ and
so if $\alpha_n=\pm1$ (i.e., we are in an exceptional situation of
type (i)), then $\alpha_n$ is a simple inverse root of
$P_k(T)/P_k^\text{red}(T)$.  If $N$ is odd then
$P_k(T)/P_k^\text{red}(T)=(1+WT)$ and so if $\alpha_n\neq-W$, then
$\alpha_n$ is not an inverse root of $P_k(T)/P_k^\text{red}(T)$,
whereas if $\alpha_n=-W$ (i.e., we are in an exceptional situation of
type (ii)), then $\alpha_n$ is a simple inverse root of
$P_k(T)/P_k^\text{red}(T)$.

This completes the proof of the theorem in the case appearing in the
section title, and thus the proof of all of part (1) of the theorem.

\subsection{Part (2) of Theorem~\ref{thm:main}}
The argument is similar to that for part (1).  We choose $D$ and local
conditions $(S_n,C_{n,w})$ according to the recipe in
\ref{ss:typical-local} and the second paragraph of
\ref{ss:special-local} and construct $X$ and a sheaf
$\GG_o=\GG(\FF,o)$ on $X$ for each orbit $o\subset\Z/d\Z$.  Then we
choose $\ell$-adic units $\beta_o$ as in \ref{prop:good-beta} and
consider $\GG(\FF\tensor\beta_o^\text{deg},o)$, its arithmetic
monodromy group $G_o^\text{arith}$, and its compact form
$K_o^\text{arith}$.

Applying Proposition~\ref{prop:good-pts-density} and
Corollary~\ref{cor:non-vanishing-density}, we find that for all
sufficiently large $n$, there exists an element $f\in X(\Fqn)$
satisfying the local conditions imposed by $(S_n,C_{n,w})$ and such
that for all $o$, $\beta_o\alpha_n$ is not an inverse root of
$P^\text{red}_{\theta(f,n)}(T)$.  Thus we are reduced to considering
the zeroes of $P_k(T)/P^\text{red}_k(T)$ where $k={\theta(f,n)}$.

In the exceptional situations of type (i)-(iv) and in the
non-exceptional situations, the analysis is exactly as for part (1).
The exceptional situations of type (v) and (vi) are like those of type
(iii) and (iv), except that in the former, we assume that
hypothesis~\ref{eqn:cond2} fails.  We used this hypothesis to show, in
Proposition~\ref{prop:avoiding-induced-zeroes}, that {\it for one
  orbit\/} $o$, if \ref{eqn:cond2} fails, we can choose local
conditions so that the rank of $\GG_o$ is even, and so
$P_k(T)/P^\text{red}_k(T)=1$.  But as we
already remarked after \ref{prop:avoiding-induced-zeroes}, it is not
possible in general to do this for several orbits $o$ at once.  In 
\ref{ss:special-local} we chose local conditions to handle possible 
trouble with the orbit $o=\{d/2\}$ (when $d$ is even) and so we 
have no control over the orbits appearing in exceptional situations 
of types (v) and (vi).  So in
these situations, the rank of $\GG_o$
may be odd or even, and
$P_k(T)/P^\text{red}_k(T)$ may be 1, so that
we have non-vanishing of the $L$-function, or it may be
$1+\sgn(\rho)T^{a_o}$, so that we have simple vanishing 
(if $\alpha_n^{a_o}=-\sgn(\rho)$) or non-vanishing (if not) of the
$L$-function.  Thus the conclusion is that we have vanishing to 
order at most 1, as desired.  This completes the proof of part (2) 
of Theorem~\ref{thm:main}.

\subsection{Proof of Theorem~\ref{thm:intro-main}}
We want to apply Theorem~\ref{thm:main}, part (2), to the data $F$, $\rho$,
and $d$, setting $S_s=S_i=S_r=\emptyset$ and $\alpha_n=q^{-ns_0}$.
The hypotheses of \ref{ss:data} are satisfied, except possibly
\ref{sss:rho-duality-hyp}.  But if $\chi$ is a character of
$G/G_{\infty}$, then the truth of Theorem~\ref{thm:intro-main} for
$\rho\tensor\chi$ and all $s_0$ implies the truth of
Theorem~\ref{thm:intro-main} for $\rho$ and all $s_0$.  Thus we may
legitimately apply Theorem~\ref{thm:main}.

Since we assume that $d|q-1$, all the orbits $o\subset\Z/d\Z$ are
singletons.  In particular, the exceptional situations of types
(iii)-(vi) do not occur.  Exceptional situations of types (i) or (ii)
can occur only if $d$ is even, $\rho$ is symplectically self-dual and
the exponent of the local Artin conductor $\cond_v(\rho)$ is even for
all $v$.  (This is what \ref{eqn:cond1} says when
$S_s=S_i=S_r=\emptyset$.)  If no exceptional situations occur, then
Theorem~\ref{thm:intro-main} follows immediately from
Theorem~\ref{thm:main}, and we may even replace ``infinitely many
$n$'' with ``all sufficiently large $n$.''

So let us assume that $d$ is even, $\rho$ is symplectically self-dual,
and \ref{eqn:cond1} is satisfied.  Then
$\deg(\cond(\rho\tensor\chi_f^{d/2}))$ is even for all $f$ satisfying
the local conditions and so exceptional situation (ii) is in fact
impossible.  Exceptional situation (i) occurs only if the root number
$W(\rho\tensor\chi_f^{d/2},F_n)=-1$.  But if this happens then for any
even multiple $m$ of $n$, $W(\rho\tensor\chi_f^{d/2},F_m)=1$
(cf.~\ref{sss:funl-eqn}) and so we avoid all exceptional situations.
Thus there are infinitely many values of $n$ for which there exists a
good $f$.  This completes the proof of Theorem~\ref{thm:intro-main}.

With slightly more work, one can prove that
Theorem~\ref{thm:intro-main} holds with ``infinitely many $n$''
replaced by ``all sufficiently large even $n$,'' and in many cases by
``all sufficiently large $n$.''

           \section{Application to elliptic curves}\label{s:elliptic-app}

The goal of this section is to prove the following two theorems.

\begin{thm}\label{thm:semistable}
  Let $\Curve$ be a geometrically irreducible curve over a finite
  field $\Fq$ of characteristic $p>3$ and let $F=\Fq(\Curve)$.  Let
  $E$ be a non-isotrivial elliptic curve over $F$.  Then there exists
  a finite separable extension $F'/F$ such that:
\begin{enumerate}
\item[(a)] $E$ has split multiplicative reduction at some place of $F'$
\item[(b)] $E$ is semistable over $F'$, i.e., it has good or multiplicative
      reduction at every place of $F'$
\item[(c)] $\ord_{s=1}L(E/F',s)=\ord_{s=1}L(E/F,s)$
\end{enumerate}
\end{thm}

\begin{thm}\label{thm:Heegner}
  Let $\Curve$ be a geometrically irreducible curve over a finite
  field $\Fq$ of characteristic $p>3$ and let $F=\Fq(\Curve)$ and
  $F_n=\Fqn(\Curve)$.  Let $E$ be a non-isotrivial elliptic curve over
  $F$ of conductor $\n$.
\begin{enumerate}
\item Fix three finite, pairwise disjoint sets of places $S_s$, $S_i$,
  $S_r$ of $F$.  Then for all sufficiently large $n$ relatively prime
  to some integer $B$, there is a quadratic extension $K/F_n$ such
  that
\begin{equation*}
\ord_{s=1}L(E/K,s)\le\ord_{s=1}L(E/F,s)+1
\end{equation*}
and such that the places of $F_n$ over $S_s$ (resp.~$S_i$, $S_r$) are
split (resp.~inert, ramified).
\item If $E$ has split multiplicative reduction at some place $\infty$
  of $F$ and we let $S_s=|\n|\setminus\{\infty\}$, $S_i=\emptyset$ and
  $S_r=\{\infty\}$, then for all sufficiently large $n$ prime to $B$
  there exists $K$ as above so that $\ord_{s=1}L(E/K,s)$ is odd.  In
  particular, if $\ord_{s=1}L(E/F,s)=1$, then $\ord_{s=1}L(E/K,s)=1$.
  The same conclusion holds if we take $S_s=|\n|\setminus\{\infty\}$,
  $S_i=\{\infty\}$ and $S_r=\emptyset$.
\end{enumerate}
\end{thm}

Theorem~\ref{thm:intro-elliptic} of the introduction is an immediate
consequence.  Indeed, we first apply \ref{thm:semistable} to find a
suitable $F'$, then apply the second part of \ref{thm:Heegner}, with
$F'$ playing the role of $F$, to find $K$.

To prove these two theorems, we will apply Theorem~\ref{thm:main} to
the representation $\rho$ of $\gal(\overline{F}/F)$ on the Tate module
$V_\ell(E)$ for some $\ell\neq p$.  Note that $\rho$ is symplectically
self-dual of weight 1 and it satisfies the hypotheses of
Subsection~\ref{ss:data}.  We have $L(E/F,s)=L(\rho,F,q^{-s})$.

\subsection{Proof of \ref{thm:Heegner}}
We begin with an easy lemma.

\begin{lemma}\label{lemma:constant-exts}
If $E$ is an elliptic curve over
$F=\Fq(\Curve)$ and if $F_n=\Fqn(\Curve)$, then there exists an integer
$b$ such that $\ord_{s=1}L(E/F_n,s)=\ord_{s=1}L(E/F,s)$ for all $n$
relatively prime to $b$.
\end{lemma}

\begin{proof}
First assume that $E$ is non-constant, so that $L(E/F,s)$ is a
polynomial in $q^{-s}$.  Writing
$L(E/F,s)=\prod_{i=1}^N(1-\alpha_iq^{-s})$ we have that
$\ord_{s=1}L(E/F,s)$ is the number of $\alpha_i$ which are equal to
$q$.  On the other hand,
$L(E/F_n,s)=\prod_{i=1}^N(1-\alpha_i^nq^{-ns})$ and so
$\ord_{s=1}(L(E/F_n,s)$ is equal to the number of $\alpha_i$
satisfying $\alpha_i^n=q^n$.  Thus we may take $b$ to be the least
common multiple of the orders of all roots of unity appearing in the
set $\{\alpha_i/q|i=1,\dots,N\}$.

If $E$ is constant, the argument is similar, except that $L(E/F,s)$ is
now a polynomial in $q^{-s}$ divided by $(1-q^{-s})(1-q^{2-s})$.
\end{proof}

\subsubsection{}
The first part of Theorem~\ref{thm:Heegner} is an easy consequence of
the main Theorem~\ref{thm:main}.  Indeed, 
Lemma~\ref{lemma:constant-exts} says that for all $n$ prime to $b$,
$\ord_{s=1}L(E/F_n,s)=\ord_{s=1}L(E/F,s)$.  On the other hand,
Theorem~\ref{thm:main}, applied with $d=2$, the given $S_s$, $S_i$,
and $S_r$, and $\alpha_n=q^n$, says that for all sufficiently large $n$ (prime
to $a_o=1$) there exists an $f\in F_n^\times$ such that the quadratic
extension $K=F_n(\sqrt{f})$ satisfies the local conditions imposed by
$S_s$, $S_i$, and $S_r$ and with 
\begin{equation*}
\ord_{s=1}\frac{L(E/K,s)}{L(E/F_n,s)}=
L(E/F_n,\chi_f,s)\le1
\end{equation*}
where $\chi_f$ is the quadratic character of $F_n$ associated to $K$.
Moreover, we can conclude that $\ord_{s=1}L(E/F_n,\chi_f,s)=0$ unless we are
in one of the exceptional situations (i) or (ii).

For the second part of Theorem~\ref{thm:Heegner}, we take
$B=b\deg\infty$, so that if $n$ is prime to $B$, then there is a
unique place of $F_n$ over $\infty$.  The assertion is that the sign
in the functional equation of $L(E/K,s)$ is $-1$, and this follows
easily from the factorization of the sign into a product of local
factors.  Indeed, over each place of $F_n$ in $|\n|\setminus\infty$
there are two places of $K$ and the local root number there are equal
and so cancel.  The only remaining contribution is at the unique place
of $K$ over $\infty$ (which is unique because we have assumed $\infty$
is inert or ramified in $K$).  There $E$ is split multiplicative and
the local contribution is $-1$.  This means that the sign in the
functional equation of $L(E/K,s)$ is $-1$, i.e., the $L$-function
vanishes to odd order.

This completes the proof of Theorem~\ref{thm:Heegner}.
\qed

\subsection{Proof of Theorem~\ref{thm:semistable}}
  For brevity, we say that an extension $F'$ of $F$ is ``good'' if
  $\ord_{s=1}L(E/F',s)=\ord_{s=1}L(E/F,s)$.  Theorem~\ref{thm:main}
  guarantees the existence of good extensions $F'=F_n(f^{1/d})$
  satisfying various local conditions for $n$ sufficiently large and
  prime to $a=[\Fq(\mu_d):\Fq]$ and the $b$ of
  Lemma~\ref{lemma:constant-exts}.
  
  We proceed in three main steps.  First we find a good extension $F'$
  of $F$ such that $E$ has a place of split multiplicative reduction
  over $F'$.  Then we replace $F$ with $F'$ and eliminate places of
  reduction types $II$, $II^*$, $IV$ and $IV^*$ (i.e., we replace $F$
  with a good extension such that there are no places of these types).
  Lastly we eliminate places of reduction types $III$, $III^*$, and
  $I_0^*$.

\subsubsection{Step 1:}
Since $E$ is assumed to be non-isotrivial, its $j$-invariant is
non-constant and thus has a pole at some place $v_0$ of $F$.  Thus $E$
is potentially multiplicative at this place.  There are three
possibilities: (i) $E$ is split multiplicative at $v_0$; (ii) $E$ is
non-split multiplicative at $v_0$; (iii) $E$ has reduction type
$I_n^*$ for some $n>0$.  

In case (i) there is nothing to do for the
first step and we set $F'=F$.  

In case (ii) we need a quadratic extension in which $v_0$ is inert.
If the integer $b$ appearing in Lemma~\ref{lemma:constant-exts} is
odd, we may take $F'=F_2=\F_{q^2} F$.  If $b$ is even, we need a
geometric extension.  For the rest of step 1, we set $d=2$ and
$\alpha_n=q^n$.  Set $S_s=S_r=\emptyset$ and $S_i=\{v_0\}$.  If the
hypothesis \ref{eqn:cond1} fails, or if it holds and the signs
appearing in Lemma~\ref{lemma:global-signs} are $+1$ then we are not
in an exceptional situation and so for large enough $n$ prime to $b$
Theorem~\ref{thm:main} supplies us with a good quadratic extension
$F'$ of $F_n$ such that $v_0$ is inert.  In the case where hypothesis
\ref{eqn:cond1} holds and the signs appearing in
Lemma~\ref{lemma:global-signs} are $-1$ then we are in an exceptional
situation and we proceed in two substeps.  First we set
$S_r=S_i=\emptyset$ and $S_s=\{v_0\}$.  By Lemma~\ref{lemma:signs}(2)
the signs appearing in Lemma~\ref{lemma:global-signs} are now $+1$ and
we can find a good quadratic extension $F'$ of $F_n$ for some large
$n$ prime to $b$ where $v_0$ is split.  Replacing $F$ with $F'$ we now
have two places of multiplicative reduction, call them $v_0$ and
$v_1$.  Setting $S_s=S_r=\emptyset$ and $S_i=\{v_0\}$ we see that
hypothesis \ref{eqn:cond1} fails (because of $v_1$) and so we are not
in an exceptional situation.  The argument in the first part of case
(ii) gives us a quadratic extension $F'$ of $F$ where $v_0$ is inert
and so $E$ has split multiplicative reduction at the place of $F'$
over $v_0$.  This completes the analysis in case (ii).

In case (iii) we will find a quadratic extension in which $v_0$ is
ramified.  We set $S_s=S_i=\emptyset$ and $S_r=\{v_0\}$.  For any
ramified local character $\chi_{v_o}$ at $v_0$, we have
$\cond_{v_0}(\rho\tensor\chi_{v_0})=1$ which is odd, so the hypothesis
\ref{eqn:cond1} fails and we are not in an exceptional situation.
Then for $n$ large and prime to $b$ Theorem~\ref{thm:main} supplies a
good quadratic extension $F'$ of $F_n$ in which every place over $v_0$
is ramified.  Then $E$ will have multiplicative reduction at each
place of $F'$ over $v_0$.  If necessary, i.e., if the reduction is not
split multiplicative, then we replace $F$ with $F'$ and apply the
argument of case (ii) again to find a good extension over which $E$ is
split multiplicative.

We now replace $F$ with $F'$ and so we may assume that $E$ has a place
of split multiplicative reduction over $F$.  This property is
preserved in arbitrary finite extensions of $F$ so we may forget about
it for the rest of the proof.

\subsubsection{Step 2:}
For a finite extension $F'$ of $F$ and an integer $m$, we let
$S_m(F')$ be the set of places $v$ of $F'$ where $E$ has additive
reduction and $m=12/gcd(v(\Delta_v),12)$ where $\Delta_v$ is the
discriminant of a minimal model of $E$ at $v$.  Thus $S_m$ consists of
places of reduction type $I_0^*$ for $m=2$, types $IV$ and $IV^*$ for
$m=3$, types $III$ and $III^*$ for $m=4$, and types $II$ and $II^*$
for $m=6$, and $S_m$ is empty for other values of $m$.  We need to
find a good extension $F'$ of $F$ such that $S_m(F')$ is empty for all
$m$.  To do this we use the well-known fact that $E$ obtains good
reduction over any place of $S_m(F)$ which is ramified of index a
multiple of $m$.  (Here we use crucially that $p>3$.)

In step 2, we will find a good extension $F'$ so that $S_3(F')$ and
$S_6(F')$ are empty.  For the rest of this step (except the very end)
we let $d=3$ and $\alpha_n=q^n$.  Theorem~\ref{thm:main} will supply
us with good cubic extensions of $F_n$ in which the places over
$S_3(F)\cup S_6(F)$ are totally ramified.  If $F'$ is such an extension,
then places of $F'$ over $S_6(F)$ are in $S_2(F')$ and places of $F'$
over $S_3(F)$ are places of good reduction.  Thus replacing $F$ with
$F'$ we will reduce to the case where $S_3(F)$ and $S_6(F)$ are empty.

To start, let $S_r=S_3(F)\cup S_6(F)$, and $S_s=S_i=\emptyset$.  If
$q\equiv1\pmod3$, the extensions $F_n(f^{1/3})/F_n$ are Galois and we
are in a non-exceptional situation.  Theorem~\ref{thm:main} supplies
us with good cubic extensions in which the places of $S_3(F)$ and
$S_6(F)$ are totally ramified.

If $q\equiv2\pmod3$ but the integer $b$ of
Lemma~\ref{lemma:constant-exts} is odd, then we may replace $F$ with
$F_2$ and then proceed as in the previous paragraph.

If $q\equiv2\pmod3$ and $b$ is even, we again set $S_r=S_3(F)\cup
S_6(F)$, and $S_s=S_i=\emptyset$ and consider the integer $N_o$ where
$o\subset(\Z/3\Z)$ is the orbit of multiplication by $q$ not containing
0.  If $N_o$ is even, we are not in an exceptional situation and we
obtain a good cubic extension as above.  If $N_o$ is odd, then we
are in an exceptional situation of type (iii) and so we will modify
our input data.  Note that the parity of $N_o$ is the same as the
parity of
\begin{equation*}
\sum_{v\text{ over }|\n|\cap S_r}\cond_v(\rho\tensor\chi_f)\deg v 
+\sum_{v\text{ over }|\n|\setminus S_r}\cond_v(\rho)\deg v 
\end{equation*}
for any $f\in X(\Fqn)$ satisfying the local conditions.  Thus one of
these sums is odd.  If the second sum is odd, then $E$ must have a
place of multiplicative reduction of odd degree.  If $v$ is such a
place, then $\cond_v(\rho)=1$ but $\cond_v(\rho\tensor\chi_f)=2$, and
so if we change $S_r$ to $S_3(F)\cup S_6(F)\cup\{v\}$, then $N_o$ is
now even and we may proceed as in the first part of this paragraph.
If the first sum is odd, we make a preliminary quadratic extension
using Theorem~\ref{thm:main}.  More precisely, we set $d=2$,
$S_s=S_r=\emptyset$, $S_i=S_3(F)\cup S_6(F)$, and $\alpha_n=q^n$.
Because we have a place of split multiplicative reduction, this is not
an exceptional situation and we find a good quadratic extension $F'$
of $F_n$ for $n$ large and relatively prime to $b$.  Now every place
of $S_3(F')\cup S_6(F')$ has even degree.  Replacing $F$ with $F'$ we
return to the setup with $d=3$, $S_r=S_3(F)\cup S_6(F)$,
$S_s=S_i=\emptyset$, and $\alpha_n=q^n$.  Now we have that every place
in $S_r$ has even degree and so either $N_o$ is even or the second
displayed sum is odd and we may proceed as in the first part of this
paragraph. 

Applying step 2 iteratively, replacing $F$ with $F'$ at each
iteration, we may now aassume that $S_3(F)$ and $S_6(F)$ are empty.

\subsubsection{Step 3:}
Now we use quadratic extensions to eliminate $S_2(F)$ and $S_4(F)$.
Let $d=2$, $S_s=S_i=\emptyset$, $S_r=S_2(F)\cup S_4(F)$, and
$\alpha_n=q^n$.  Since we have a place of multiplicative reduction,
hypothesis~\ref{eqn:cond1} fails and so we are in a non-exceptional
situation.  Theorem~\ref{thm:main} gives us a good quadratic extension
$F'$ of $F_n$ for some large $n$ prime to $b$ in which every place of
$F_n$ over $S_r$ is ramified.  This means that $E$ acquires good
reduction at every place over $S_2(F)$, $S_4(F')$ is empty and
$S_2(F')$ consists of precisely the places over $S_4(F)$.  Replacing
$F$ with $F'$ and repeating this construction once more yields a good
extension $F'$ where $S_m(F')$ is empty for all $m$.  This $F'$ is an
extension of the original $F$ with all the required properties and
this completes the proof of Theorem~\ref{thm:semistable}.  \qed

\providecommand{\bysame}{\leavevmode\hbox to3em{\hrulefill}\thinspace}
\providecommand{\MR}{\relax\ifhmode\unskip\space\fi MR }
\providecommand{\MRhref}[2]{%
  \href{http://www.ams.org/mathscinet-getitem?mr=#1}{#2}
}
\providecommand{\href}[2]{#2}

\end{document}